\documentclass[preprint,12pt]{elsarticle}

\usepackage{amssymb}
\usepackage{amsmath}
\usepackage{subcaption}
\usepackage{physics}
\usepackage{hyperref}
\usepackage{siunitx}

\journal{Journal of Computational Physics}

\begin{document}

\begin{frontmatter}



\title{Fast Summation on the Sphere with Applications to the Barotropic Vorticity Equation}


\author[inst1]{Anthony Chen}

\affiliation[inst1]{organization={Department of Mathematics, University of Michigan},
            addressline={530 Church St}, 
            city={Ann Arbor},
            postcode={48109}, 
            state={MI},
            country={USA}}

\author[inst2]{Christiane Jablonowski}

\affiliation[inst2]{organization={Department of Climate and Space Sciences and Engineering, University of Michigan},
            addressline={2455 Hayward St}, 
            city={Ann Arbor},
            postcode={48109}, 
            state={MI},
            country={USA}}

\begin{abstract}
Fast summation refers to a family of techniques for approximating $O(N^2)$ sums in $O(N\log{N})$ or $O(N)$ time. These techniques have traditionally found wide use in astrophysics and electrostatics in calculating the forces in a $N$-body problem. In this work, we present a spherical tree code, and apply it to the problem of efficiently solving the barotropic vorticity equation. 
\end{abstract}



\end{frontmatter}


\section{Introduction}
\label{sec:intro}

Fast summation refers to a wide family of techniques for approximating sums that have $O(N^2)$ time complexity in $O(N\log{N})$ or $O(N)$ time. For example, when one has $N$ particles interacting electrostatically, the potential of each particle is 
\begin{equation}
    \phi(\mathbf{x}_i)=\sum_{j=1,j\neq i}^NG(\mathbf{x}_i,\mathbf{x}_j)q_j
\end{equation}
where $\{\mathbf{x}_i\}_{i=1}^N$ is a set of particles, $G$ is the Coulomb potential, and $q_i$ is the charge of particle $\mathbf{x}_i$. Similar sums appear in many other fields of physics, such as the astrophysical N-body problem, where $G$ is instead the gravitational potential and $q_j$ is replaced with $m_j$, the mass of each particle. There are a wide range of fast summation techniques that have been created, such as Appel's tree code \cite{appel1985efficient} which introduces the dual tree traversal idea, followed by the Barnes-Hut tree code \cite{barneshut} which is based on monopole expansions, and the Fast Multipole Method \cite{fmmgreengard} which uses higher order multipole expansions as well. Additional work has been in the direction of kernel independence, with the creation of the kernel independent Fast Multipole Method \cite{ying2004kernel} which generalizes beyond the Coulomb potential and the Blackbox Fast Multipole Method \cite{fong2009black} which uses polynomial interpolation for its kernel independence instead of multipole expansions. More recent techniques have used barycentric Lagrange interpolation, such as the Barycentric Lagrange Tree Code \cite{bltc}, and the Barycentric Lagrange Dual Tree Traversal Fast Multipole Method\cite{bldtt}. Past work has focused on fast summation in $\mathbb{R}^n$, as electrons and stars move in Cartesian space, to the neglect of other geometries. Many problems are naturally set on the sphere, and we can take advantage of this geometry in our fast summation. A number of problems in geophysics can naturally be formulated to take advantage of fast summation techniques.  

The use of Lagrangian particle methods for fluid dynamics has a long history. Chorin introduced the vortex method for incompressible fluid problems \cite{chorin1973numerical}, and work has continued in this area \cite{mimeau2021review}, with good results for challenging problems, as seen in \cite{xu2023dynamics}. On the sphere, there has been work focusing on point vortices~\cite{kidambi1998motion,sakajo1999motion}, vortex patches~\cite{surana2008vortex}, and vortex sheets~\cite{sakajo2004motion,sakajo2009extension}. More recent work has focused on general vorticity distributions~\cite{bosler2014lagrangian} and the shallow water equations~\cite{capecelatro2018purely}, with both Lagrangian particle methods and smoothed particle hydrodynamics demonstrating their power. Some of these Lagrangian solvers have applied fast summation techniques~\cite{sakajo2009extension,kropinski2014fast}, but did not set the problem on a sphere, instead setting the problem in $\mathbb{R}^3$ or taking a stereographic projection to a plane. Lagrangian methods have a number of advantages, which include not being constrained by a CFL condition for stability, and relatively easier adaptive mesh refinement. 

In this work, we design a spherical tree code that explicitly accounts for the spherical geometry found in the problems of geophysics. We demonstrate the utility of this tree code by applying it to the problem of the time evolution of the barotropic vorticity equation, going to higher resolutions and longer times than what has been previously done. The rest of this paper is organized as follows. In Sec.~\ref{sec:bve}, we first present an example of a problem in geophysical fluid dynamics in which we would like to perform fast summation. In Sec.~\ref{sec:fss}, we introduce our technique for performing fast summation on a sphere. Results of numerical tests are then presented in Sec.~\ref{sec:tests}. We conclude in Sec.~\ref{sec:conc}, and present some future research directions. 

\section{A Lagrangian Discretization of the Barotropic Vorticity Equation}
\label{sec:bve}

We briefly present a Lagrangian discretization of the Barotropic Vorticity Equation (BVE), which is a simplified equation set for two dimensional incompressible inviscid fluid motion in the presence of a Coriolis force. The BVE is a good model to capture large scale atmospheric motions. The BVE is given by \begin{equation}\frac{D(\zeta+f)}{Dt}=0\end{equation} where $\frac{D}{Dt}$ is the total or material derivative, $\zeta$ is the relative vorticity, $f$ is the Coriolis parameter, and the quantity $\zeta+f$ is known as the absolute vorticity. To close this, following the presentation in \cite{vallis2017atmospheric}, we first define the stream function $\psi$ as 
\begin{equation}\label{eq:stream}
\Delta\psi=\zeta
\end{equation} 
which we can solve for $\psi$ from $\zeta$, before we compute
\begin{equation}\label{eq:bve}
\mathbf{u}(\mathbf{x})=\nabla\psi\times\mathbf{x}
\end{equation}
where $\mathbf{x}$ is a point on the sphere $S$ with radius $1$. For the discretization, we follow the method in \cite{bosler2014lagrangian}, with some small differences. First, the spherical Laplacian has a Green's function 
\begin{equation}
    G(\mathbf{x},\mathbf{y})=-\frac{1}{4\pi}\log(1-\mathbf{x}\cdot\mathbf{y})
\end{equation}
\cite{bogomolov1977dynamics, kimura1987vortex}, where $\mathbf{x}$ and $\mathbf{y}$ are points in $\mathbb{R}^3$ that lie on the surface of the sphere. Using this, we find the velocity as \begin{equation}\label{eq:velbsint}\mathbf{u}(\mathbf{x},t)=-\frac{1}{4\pi}\int_S\frac{\mathbf{x}\times\mathbf{y}}{1-\mathbf{x}\cdot\mathbf{y}}\zeta(\mathbf{y})dS(\mathbf{y})\end{equation} We next introduce the flow map \begin{equation}\label{eq:flowmap}
    \frac{\partial\mathbf{x}(\mathbf{a},t)}{\partial t}=\mathbf{u}(\mathbf{x}(\mathbf{a},t),t)
\end{equation} where $\mathbf{a}$ is a set of Lagrangian coordinates and $\mathbf{x}(\mathbf{a},t)$ is the location of the fluid particle originally at $\mathbf{a}$ at time $t$. Then substituting Eq.~\ref{eq:velbsint} into Eq.~\ref{eq:flowmap} we have that \begin{equation}\label{eq:velbsint2}
    \frac{\partial\mathbf{x}}{\partial t}(\mathbf{a},t)=-\frac{1}{4\pi}\int_S\frac{\mathbf{x}(\mathbf{a},t)\times\mathbf{x}(\mathbf{b},t)}{1-\mathbf{x}(\mathbf{a},t)\cdot\mathbf{x}(\mathbf{b},t})\zeta(\mathbf{x}(\mathbf{b},t))dS(\mathbf{b})
\end{equation} and substituting Eq.~\ref{eq:flowmap} into Eq.~\ref{eq:bve} we get that \begin{equation}
    \frac{D\zeta}{Dt}(\mathbf{x}(\mathbf{a},t),t)=-2\Omega\frac{\partial z}{\partial t}(\mathbf{a},t)
\end{equation} where we used that $f=2\Omega z$ and $\Omega$ is the rotation rate of the Earth, $2\pi\,\mathrm{day}^{-1}$. Then, we discretize the integral in Eq.~\ref{eq:velbsint2} with the midpoint rule, giving us a system of particles with a time evolution given by 
\begin{align}\begin{split}
    \frac{d\mathbf{x}_i}{dt}&=-\frac{1}{4\pi}\sum_{j=1,j\neq i}^N\frac{\mathbf{x}_i\times\mathbf{x}_j}{1-\mathbf{x}_i\cdot\mathbf{x}_j}\zeta_jA_j\\ \frac{d\zeta_i}{dt}&=-2\Omega\frac{dz_i}{dt}
\end{split}\end{align} where $\zeta_i$ is the relative vorticity of particle $\mathbf{x}_i$ and $A_i$ is the area for the discretized quadrature, where we now have regular derivatives $\frac{d}{dt}$ as these represent the movement of particles in space. For this work, we choose a set of particles given by an iteratively refined icosahedron, with the areas representing the faces of the dual polyhedron. These areas are known in the literature as the node patch areas, which have been shown to have favorable error behavior as compared to centroid based areas \cite{wilson2022}. We can now apply fast summation to the particles. Note that unlike \cite{bosler2014lagrangian}, we do not use panels and only have particles. The tree structure of the particles is useful when performing adaptive mesh refinement (AMR). One special feature of the BVE is that they admit nontrivial exact solutions, the Rossby-Haurwitz waves \cite{haurwitz1940motion}, which we can use to measure the error of our numerical scheme exactly. We present more details of this in Sec.~\ref{sec:tests}. 

We note that the kernel as written here, for particles on the sphere, is only singular when $\mathbf{x}=\mathbf{y}$. However, if $\mathbf{x}$ or $\mathbf{y}$ lie outside of the sphere, then it is possible for $1-\mathbf{x}\cdot\mathbf{y}$ to equal $0$ even when $\mathbf{x}\neq\mathbf{y}$. It is possible to slightly rewrite the kernel, since \begin{equation}\frac12\norm{\mathbf{x}-\mathbf{y}}^2=\frac{\langle\mathbf{x}-\mathbf{y},\mathbf{x}-\mathbf{y}\rangle}{2}=\frac{\norm{\mathbf{x}}^2+\norm{\mathbf{y}}^2}{2}-\mathbf{x}\cdot\mathbf{y}\end{equation} where $\langle\cdot,\cdot\rangle$ is the inner product, so when $\norm{\mathbf{x}}=\norm{\mathbf{y}}=1$ then we have that \begin{equation}\label{eq:altkernel}\frac12\norm{\mathbf{x}-\mathbf{y}}^2=1-\mathbf{x}\cdot\mathbf{y}\end{equation} but this kernel is also well defined off the sphere. This will be important later. 

There are several reasons that one might choose a Lagrangian discretization over an Eulerian scheme \cite{lagrangeatmos}. Lagrangian discretizations automatically satisfy a CFL condition so time stepping is limited by accuracy instead of stability, and Lagrangian schemes have advantages in tracer advection and conservation, as seen in \cite{lauritzen2014standard}. However, there are challenges relating to particles becoming disordered with time, necessitating remeshing \cite{perlman1985accuracy, wee2006modified}. Remeshing allows us to interpolate from a disordered set of particles to an ordered set of particles, improving the error scaling. Additionally, in a Lagrangian scheme, we can perform AMR effectively because we do not need to vary our time step as spatial refinement occurs, as opposed to an Eulerian scheme, as seen in \cite{jablonowski2006block}. 

\section{Fast Summation on a Sphere}
\label{sec:fss} 

We wish to compute convolutions on the sphere, in the form
\begin{equation}
    \phi(\mathbf{x})=\int_SG(\mathbf{x},\mathbf{y})f(\mathbf{y})dS
\end{equation}
where $S$ is a sphere, $G$ is an interaction kernel or Green's function, $f$ is some scalar field defined on the surface of the sphere, $\mathbf{x}$ and $\mathbf{y}$ are points on the sphere, and $\phi$ is a potential. When discretized, one ends up with the sum
\begin{equation}  \phi(\mathbf{x}_i)=\sum_{i=1}^NG(\mathbf{x}_i,\mathbf{x}_j)f(\mathbf{x}_j)A_j
\end{equation}
where $A_j$ is the quadrature weight associated to point $\mathbf{x}_j$. In some cases $G$ may be singular, in which we can choose to regularize $G$. We refer to $\mathbf{x}_i$ as the target particle and the $\mathbf{x}_j$ as the source particles. We wish to compute this quickly. As previously seen, the problem of solving the barotropic vorticity equation can be put into this form. 

We now present an overview of the tree code, with details later in this section. We start with an icosahedron and refine its faces to give a tree structure, iteratively splitting each triangular face into four smaller triangles. Considering two triangles at a time, we designate one as the target containing the target particles, and one as the source containing the source particles. If $n_t$, the number of particles in the target triangle, or $n_s$, the number of particles in the source triangle, is greater than $N_T$, some threshold, we then make use of various approximations to reduce the amount of computation. If the two triangles are well separated in a technical sense to be defined later, we then perform the appropriate interaction. If the two triangles are not well separated, we refine each triangle as appropriate. We now go into more detail. 

\subsection{Grids on a Sphere}

There is no regular and uniform grid on a sphere, and each choice of a grid will pose its own issues \cite{spherepoints}. For this work, we have chosen to use a refined icosahedron for several reasons, the primary of which is the fact that the triangles at a given refinement level of the icosahedral mesh are relatively uniform in area. The latitude longitude grid suffers from non-uniformity where cells shrink dramatically at the poles, as does the cubed sphere grid, where cells shrink at the corners of the original cube. 

\subsection{Interpolation in Triangles}

Because we are using faces of a refined icosahedron to perform our fast summation approximations, we need to interpolate in triangles. Optimal interpolation points for the triangle are not known, but there one can compute Fekete points using an iterative optimization procedure~\cite{feketepoints}. However, for this work, we use the approximate Fekete points which can be directly computed, as presented in~\cite{liu2019}. For the interpolating functions, we use the Spherical Bernstein-B\'ezier (SBB) polynomials. For a spherical triangle with vertices $\mathbf{v}_1$, $\mathbf{v}_2$, and $\mathbf{v}_3$, any point inside of the spherical triangle can be written in the form $\mathbf{v}=\beta_1\mathbf{v}_1+\beta_2\mathbf{v}_2+\beta_3\mathbf{v}_3$ and $\beta_1,\beta_2,\beta_3\geq0$. These are the barycentric coordinates of the triangle. Then, a SBB basis polynomial is defined as \begin{equation}B_{ijk}(\mathbf{v})=\beta_1^i\beta_2^j\beta_3^k\end{equation} with $i+j+k=d$~\cite{alfeld1996bernstein}. For a fixed degree $d$, there will be $\frac12(d+1)(d+2)$ SBB basis polynomials, and general interpolating functions will be given as a sum of these basis functions. This basis has found considerable success in area of spherical spline interpolation~\cite{neamtu2004approximation, baramidze2006spherical, baramidze2018nonnegative}. 

\subsection{Treatment of Clusters}

\begin{figure}
    \centering
    \begin{subfigure}{0.45\textwidth}
        \includegraphics[width=\linewidth]{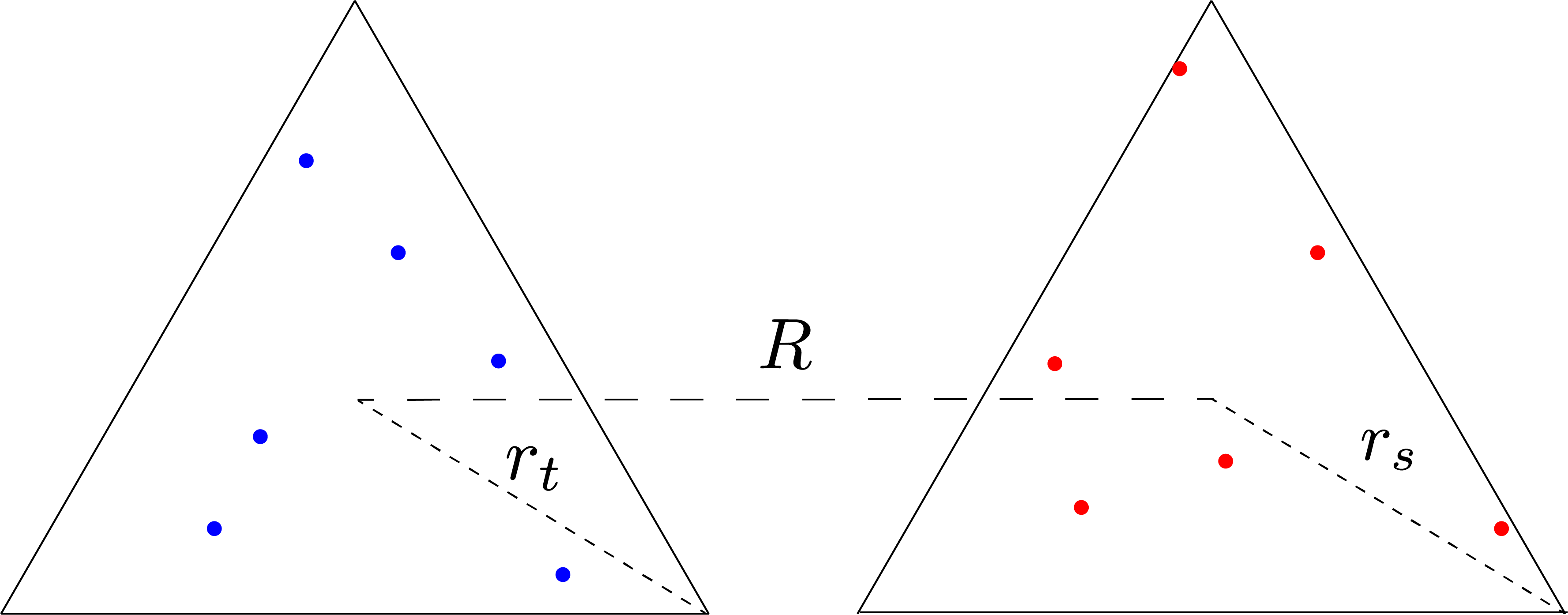}
        \caption{Particle-Particle}
        \label{fig:pptri}
    \end{subfigure}
    \hspace{0.05\textwidth}
    \begin{subfigure}{0.45\textwidth}
        \includegraphics[width=\linewidth]{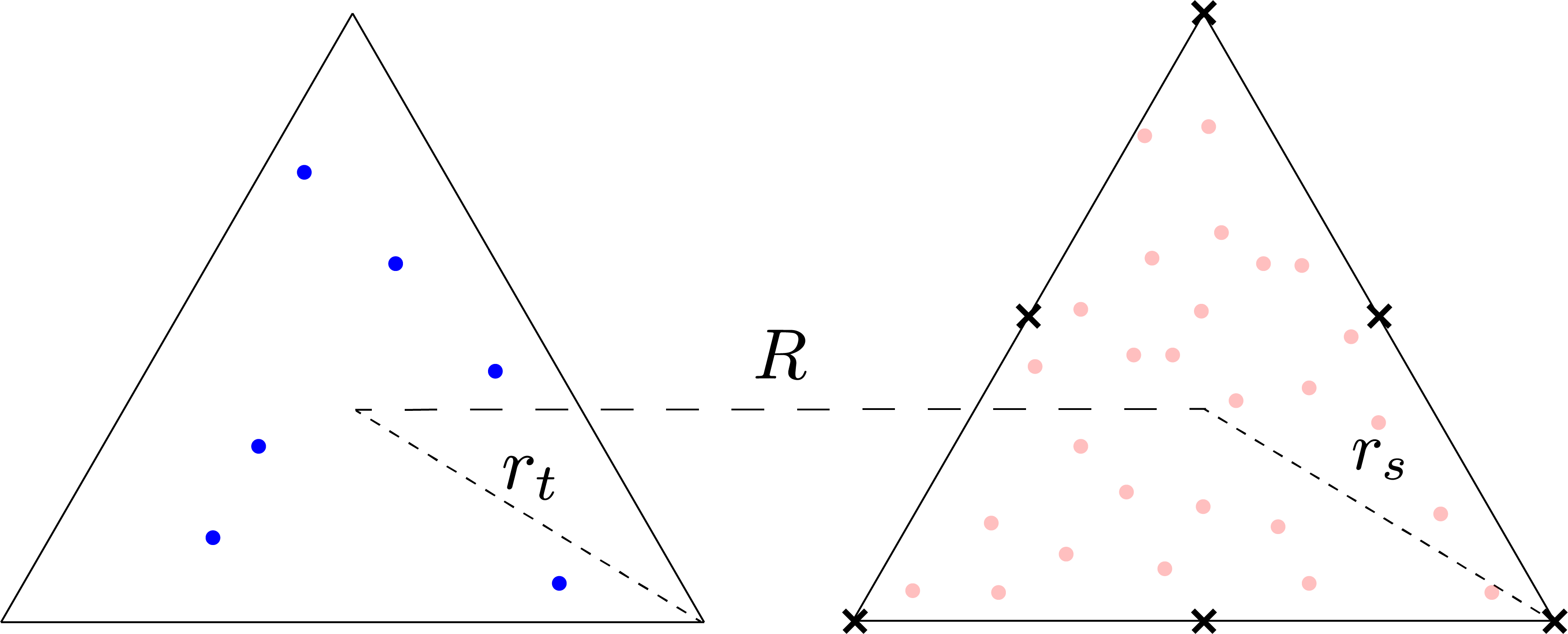}
        \caption{Particle-Cluster}
        \label{fig:pctri}
    \end{subfigure}
    \begin{subfigure}{0.45\textwidth}
        \includegraphics[width=\linewidth]{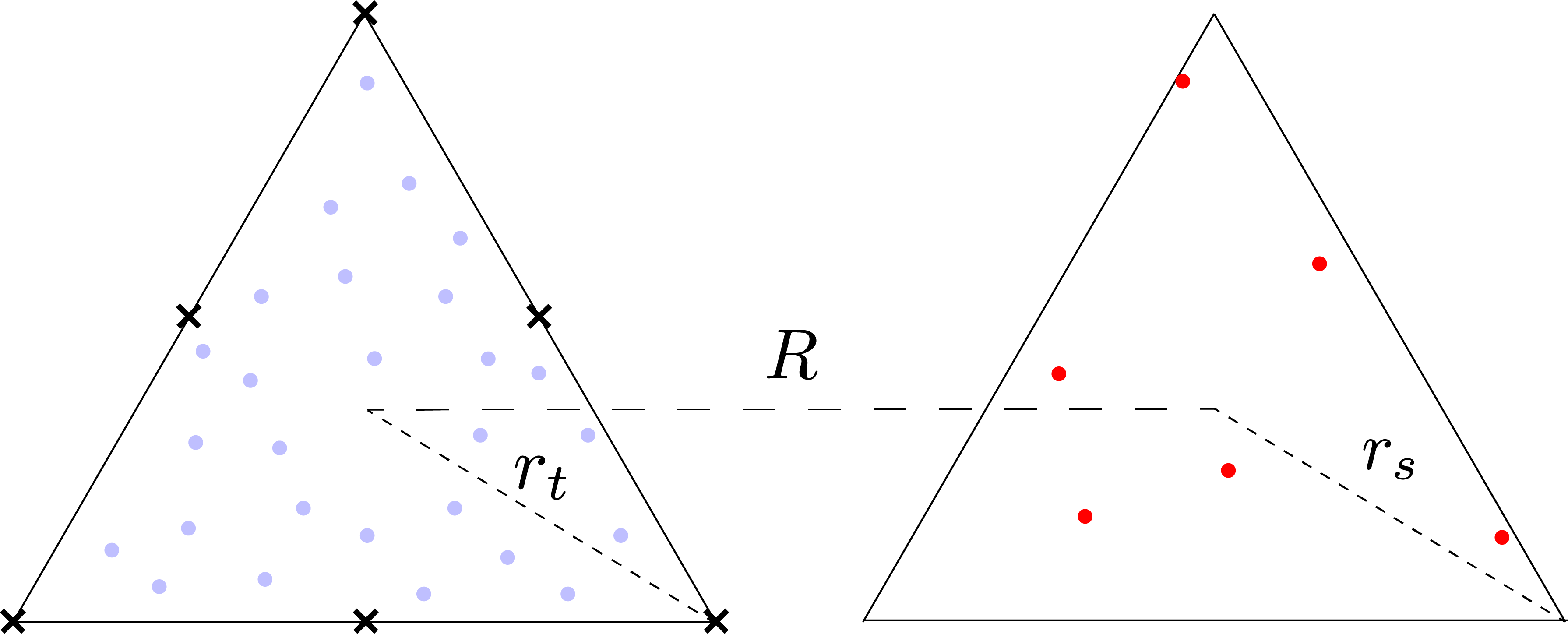}
        \caption{Cluster-Particle}
        \label{fig:cptri}
    \end{subfigure}
    \hspace{0.05\textwidth}
    \begin{subfigure}{0.45\textwidth}
        \includegraphics[width=\linewidth]{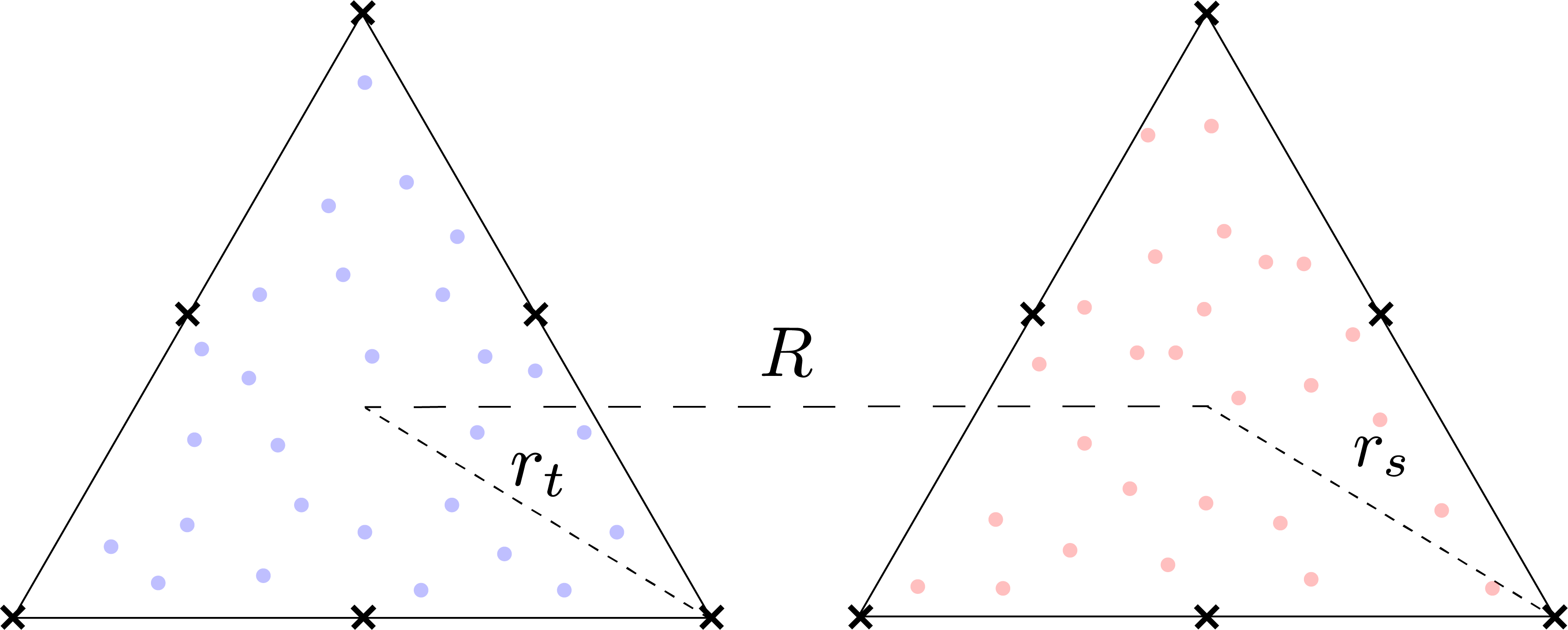}
        \caption{Cluster-Cluster}
        \label{fig:cctri}
    \end{subfigure}
    \caption{The four types of interactions are presented here. The blue particles are the target particles and the red particles are the source particles. When a triangle has many particles, as seen in \ref{fig:pctri}, \ref{fig:cptri}, and \ref{fig:cctri}, we use the proxy points marked with the crosses. $R$ is the distance between the center of the two triangles, $r_t$ is the radius of the target triangle, and $r_s$ is the radius of the source triangle.}
    \label{fig:tris}
\end{figure}

We can treat either the source or the target triangle as a cluster. This leads to four possibilities: particle-particle interactions, particle-cluster interactions, cluster-particle interactions, and cluster-cluster interactions. Let $T_t$ denote a target triangle, $T_s$ denote a source particle, and $\phi(\mathbf{x}_i,T_t,T_s)$ denote the potential at particle $\mathbf{x}_i$ in target $T_t$ due to the particles in $T_s$. These four cases are depicted in Fig.~\ref{fig:tris}. 

\subsubsection{Particle-Particle Interactions}

In a particle-particle (PP) interaction, as seen in Fig.~\ref{fig:pptri}, the two triangles both have few particles, or they are not well separated and the two triangles are leaves of the tree. In this case, we simply interact the particles directly as we would for direct summation, so for a particle $\mathbf{x}_i\in T_t$, we have that 
\begin{equation}\label{eq:pp}\phi(\mathbf{x}_i,T_t,T_s)=\sum_{\mathbf{y}_j\in T_s}G(\mathbf{x}_i,\mathbf{y}_j)q_jA_j\end{equation}
where $q_j$ is a charge or mass associated to each particle, and $A_j$ is an area. 

\subsubsection{Particle-Cluster Interactions}

In a particle-cluster (PC) interaction, as seen in Fig.~\ref{fig:pctri}, we are taking the source triangle as a cluster. In this case, we will interpolate the Green's function in the triangle. We have that \begin{equation}G(\mathbf{x}_i,\mathbf{y}_j)\approx\sum_{k}\alpha_{k}(\mathbf{x}_i)B_{k}(\mathbf{y})\end{equation} for some fixed $\mathbf{x}_i\in T_t$, an interpolation degree $d$, and where $k$ indexes over the SBB basis polynomials. Then, substituting this into Eq.~\ref{eq:pp} we get \begin{equation}\phi(\mathbf{x}_i,T_t,T_s)\approx\sum_{k}\alpha_{k}(\mathbf{x}_i)\sum_{\mathbf{y}_j\in T_s}B_{k}(\mathbf{y}_j)q_jA_j\end{equation} and we can denote \begin{equation}\omega_k=\sum_{\mathbf{y}_j\in T_s}B_k(\mathbf{y}_j)q_jA_j\end{equation} as a proxy weight that is independent of $\mathbf{x}_i$ and hence only needs to be computed once per PC interaction. To compute the $\alpha_k$, one solves a linear system. For the proxy source points $\widetilde{\mathbf{y}}_m\in T_s$, one has a linear system of the form \begin{equation}\begin{bmatrix}B_1(\widetilde{\mathbf{y}}_1)&B_2(\widetilde{\mathbf{y}}_1)&\cdots&B_N(\widetilde{\mathbf{y}}_1)\\B_1(\widetilde{\mathbf{y}}_2)&B_2(\widetilde{\mathbf{y}}_2)&\cdots&B_N(\widetilde{\mathbf{y}}_2)\\\vdots&\vdots&\ddots&\vdots\\B_1(\widetilde{\mathbf{y}}_N)&B_2(\widetilde{\mathbf{y}}_N)&\cdots&B_N(\widetilde{\mathbf{y}}_N)\end{bmatrix}\begin{bmatrix}\alpha_1(\mathbf{x}_i)\\\alpha_2(\mathbf{x}_i)\\\vdots\\\alpha_N(\mathbf{x}_i)\end{bmatrix}=\begin{bmatrix}G(\mathbf{x}_i,\widetilde{\mathbf{y}}_1)\\G(\mathbf{x}_i,\widetilde{\mathbf{y}}_2)\\\vdots\\G(\mathbf{x}_i,\widetilde{\mathbf{y}}_N)\end{bmatrix}\end{equation} which one solves to find the $\alpha_k$, then allowing the computation of the PC interaction as \begin{equation}\phi(\mathbf{x}_i,T_t,T_s)\approx\sum_k\alpha_k(\mathbf{x}_i)\omega_k\end{equation}

\subsubsection{Cluster-Particle Interactions}

In a cluster-particle (CP) interaction, as seen in Fig.~\ref{fig:cptri}, we are taking the target triangle as a cluster. Given a set of interpolation points $\widetilde{\mathbf{x}}_n\in T_t$, we first compute \begin{equation}\phi(\widetilde{\mathbf{x}}_n,T_t,T_s)=\sum_{\mathbf{y}_j\in T_s}G(\widetilde{\mathbf{x}}_n,\mathbf{y}_j)q_jA_j\end{equation} which we can then use to form an interpolating polynomial for the potential $f_{\phi}(\mathbf{x};T_s)$ for $\mathbf{x}\in T_t$ and we then simply have that 
\begin{equation}
    \phi(\mathbf{x}_i,T_t,T_s)=f_{\phi}(\mathbf{x}_i;T_s)
\end{equation}

\subsubsection{Cluster-Cluster Interactions}

In a cluster-cluster (CC) interaction, as seen in Fig.~\ref{fig:cctri}, we treat both target and source triangles as clusters, so we combine the approximations of the particle-cluster interaction and cluster-particle interaction. The interaction then proceeds as follows. We first compute the proxy weights \begin{equation}\omega_k=\sum_{\mathbf{y}_j\in T_s}B_k(\mathbf{y}_j)q_jA_j\end{equation} for the source cluster. Next, we use this to compute the potential at the proxy target particles \begin{equation}\phi(\widetilde{\mathbf{x}}_n,T_t,T_s)\approx\sum_k\alpha_k(\widetilde{\mathbf{x}}_i)\omega_k\end{equation} and we then use to form an interpolating function $f_{\phi}$ for the potential, which we then evaluate at the target particles \begin{equation}\phi(\mathbf{x}_i,T_t,T_s)=f_{\phi}(\mathbf{x}_i)\end{equation} thus combining the PC and CP approximations. 

\subsection{Dual Tree Traversal}

Now that we have the details of the various types of interaction, we can detail the tree traversal strategy, which is based on the dual tree traversal as presented in \cite{appel1985efficient}. We start with the 20 top level triangles, giving 400 pairs of source and target triangles. For each pair, we first determine if the two triangles are well separated. To do this, we use the radii of the triangle, $r_s$ and $r_t$, and the distance between the circumcenters of the triangle, $R$. If 
\begin{equation}
    \frac{r_t+r_s}{R}<\theta
\end{equation}
we say that that the two triangles are well separated. $\theta$ is known in the literature as the multipole acceptance criteria (MAC) and can be varied based on one's error tolerance. Reducing $\theta$ reduces error at the cost of increased computation. Then, if two triangles are well separated, we check the number of particles each contains against $N_T$ are perform the appropriate interaction. If the two triangles are not well separated and both contain fewer than $N_T$ triangles, we interact them directly. If the two triangles are not well separated and one or both contain many particles, we refine the triangle with more triangles, up until some limit on the tree depth. This gives us our spherical tree code, which is asymptotically $O(N\log{N})$ in the number of particles. This is reflected in the timing results in Sec.~\ref{sec:tests}. 

\section{Numerical Results}
\label{sec:tests}

We now present runtime results for both serial and parallel configurations, before moving on to discussions of error, a test of the AMR scheme, and some tests that are of particular interest to the atmospheric sciences. We do this with a series of test cases. For time stepping, we use the Runge Kutta 4th order scheme. Computations were done on the High Performance Computing system Derecho \cite{derecho} at the National Center for Atmospheric Research (NCAR). We also compare with the Barycentric Lagrange Tree Code (BLTC) \cite{bltc}, for which we use the kernel rewritten using Eq.~\ref{eq:altkernel}. The BLTC is a kernel independent tree code, with approximations based on barycentric Lagrange interpolation at Chebyshev points in $\mathbb{R}^3$. The BLTC has two parameters that are used to control the error, the interpolation degree $d$ and the MAC $\theta$. The code used for this is available at \href{https://github.com/cygnari/fast-bve}{https://github.com/cygnari/fast-bve}. 

\subsection{Runtime Results}

\begin{figure}
    \centering
    \includegraphics[scale=0.45]{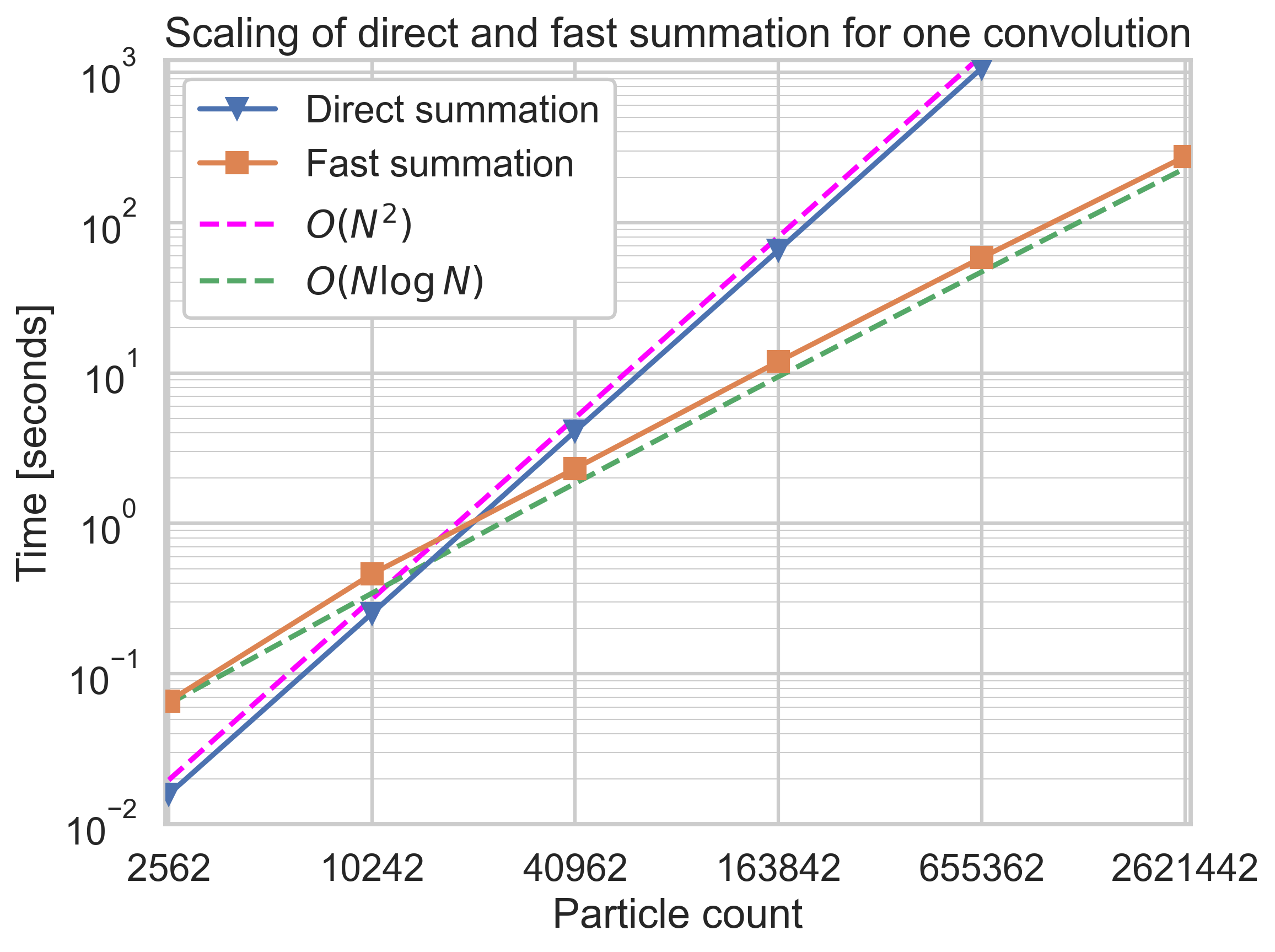}    
    \caption{We observe the speedup that the spherical tree code has over direct summation, with the expected scalings. These are run on one processor, and for the fast summation, we set $\theta=0.7$ and $d=6$, which gives between 3 and 4 digits of accuracy. }
    \label{fig:direct_fast_scaling}
\end{figure}

\begin{figure}
    \centering
    \begin{subfigure}{0.45\textwidth}
        \includegraphics[width=\linewidth]{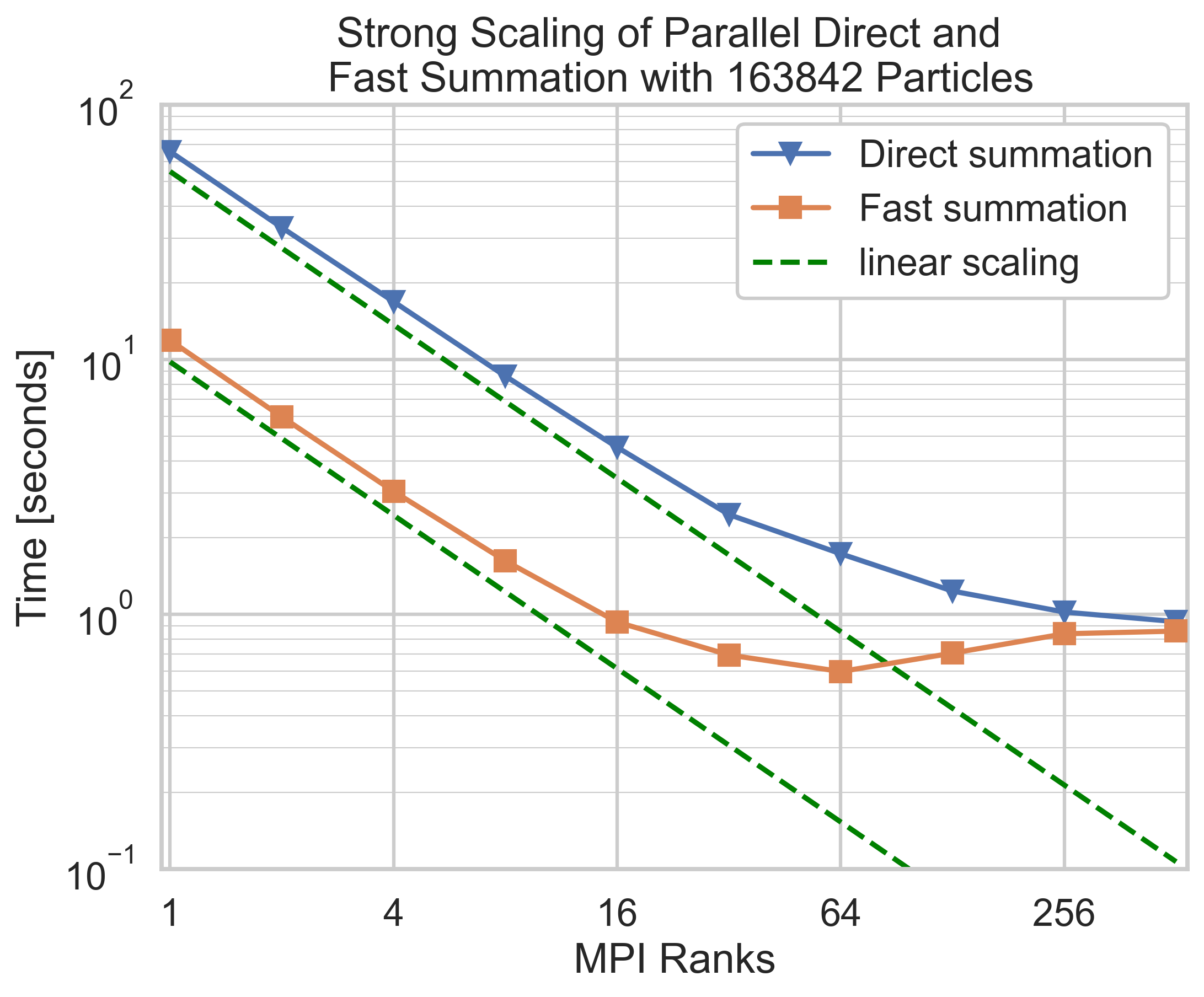}
        \caption{Strong scaling runtime}
    \end{subfigure}
    \begin{subfigure}{0.45\textwidth}
        \includegraphics[width=\linewidth]{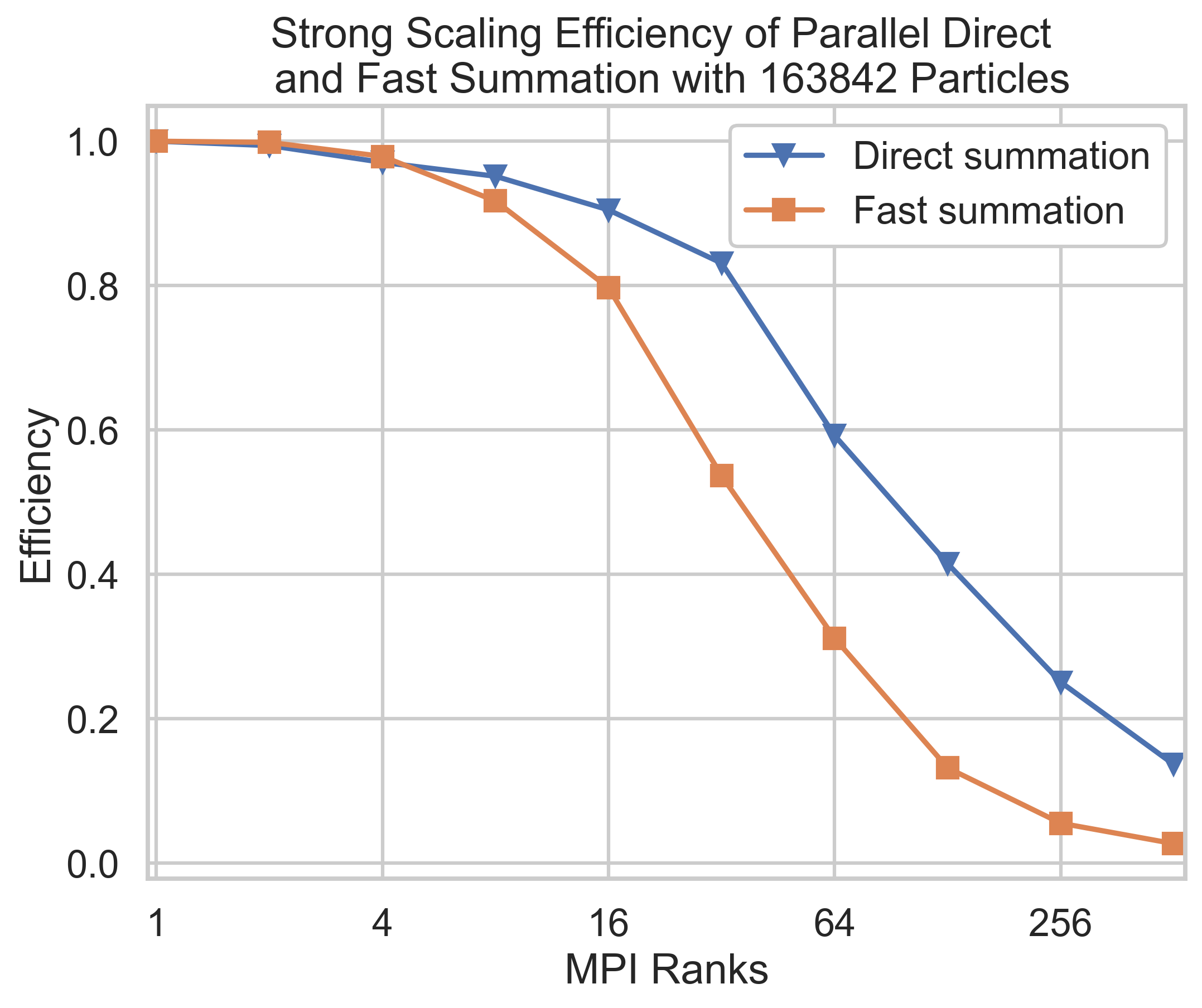}
        \caption{Strong scaling efficiency}
    \end{subfigure}
    \caption{On the left, we plot the strong scaling of both the direct and fast summation for a fixed problem size of 163842 particles. On the right hand side, we plot the parallel efficiency. For these tests, we are using $\theta=0.7$ and $d=6$. }
    \label{fig:parallelperformance}
\end{figure}

For these timing results, we test the time to perform a single convolution. These timing results were acquired on Derecho. For the serial results, we run with a single process. For the parallelization, we use MPI for both intra- and inter-node communication, making use of the Remote Memory Access paradigm \cite{mpirma}. For strong scaling results, we fix the problem size at 163842 points and increase the number of MPI ranks. First, we present the serial timings in Fig.~\ref{fig:direct_fast_scaling}, that demonstrate that the fast summation significantly improved the scaling in the problem size. At 655362 particles, the fast summation is 30 times faster than the direct summation. While we only present runtime results for only one value of $\theta$ and $d$, the $O(N\log{N})$ scaling still holds when $\theta$ and $d$ are varied, just with different constant factors. The parallel performance results are presented in Fig.~\ref{fig:parallelperformance}. We observe good scaling for the fast summation up until 64 processors, at which point the performance plateaus, as the gains from additional parallelism are outweighed by the increases in communication. 

\begin{figure}
    \centering
        \includegraphics[scale=0.45]{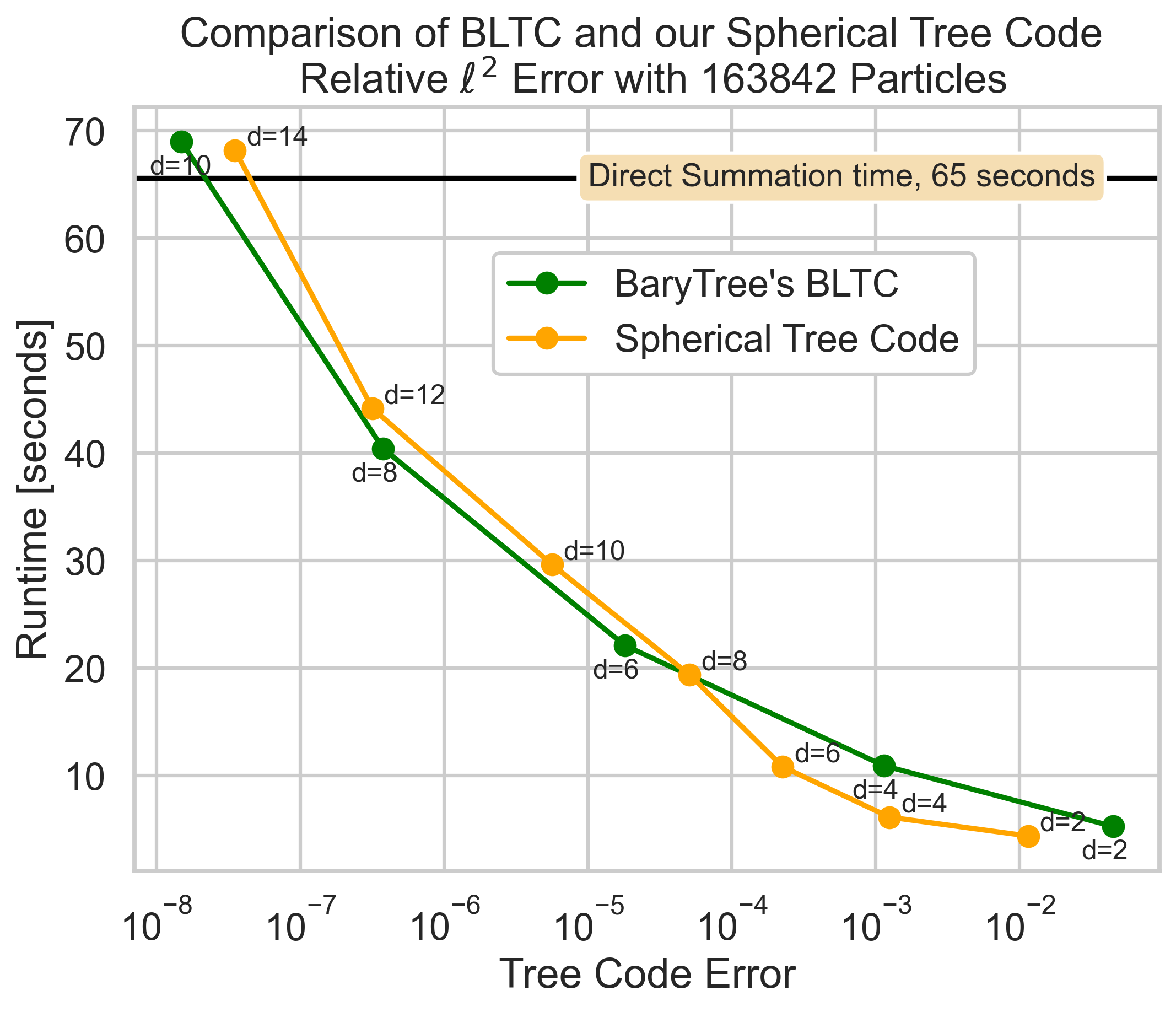}
    \caption{
	Here, we plot the error of the BLTC and our spherical tree code as we increase the interpolation degree. For the BLTC, we vary the interpolation degree from $d=2$ to $10$ in increments of $2$. For our spherical tree code, we vary the interpolation degree from $d=2$ to $d=14$ in increments of $2$. For the BLTC and the spherical tree code, we have $\theta=0.7$. The runtime for direct summation is given by the black line across the top. }
    \label{fig:treecodeerror_theta.png}
\end{figure}

As mentioned previously, performing the various fast summation approximations introduces error. We can measure this by computing the relative $\ell^2$ error between the output of direct summation and the output of the fast summation. Additionally, we can control the error of the fast summation technique by changing the interpolation degree $d$, at the cost of increasing the runtime. This is shown in Fig.~\ref{fig:treecodeerror_theta.png}. Here, comparing $\mathbf{v}_F$, the fast summation particle velocities, and $\mathbf{v}_D$, the direct summation particle velocities, we define 
\begin{equation}
    E=\sqrt{\frac{\sum_{i=1}^N\norm{\mathbf{v}_{F,i}-\mathbf{v}_{D,i}}^2A_i}{\sum_{i=1}^N\norm{\mathbf{v}_{D,i}}^2A_i}}
\end{equation}
as the relative $\ell^2$ error. For both fast summation techniques, we observe that as the interpolation degree is increased, the error decreases and the runtime increases. The spherical tree code is competitive with the BLTC across the error range presented here. 

\subsection{Rossby Haurwitz Wave}

\begin{figure}
    \centering
    \begin{subfigure}{0.45\textwidth}
        \includegraphics[width=\linewidth]{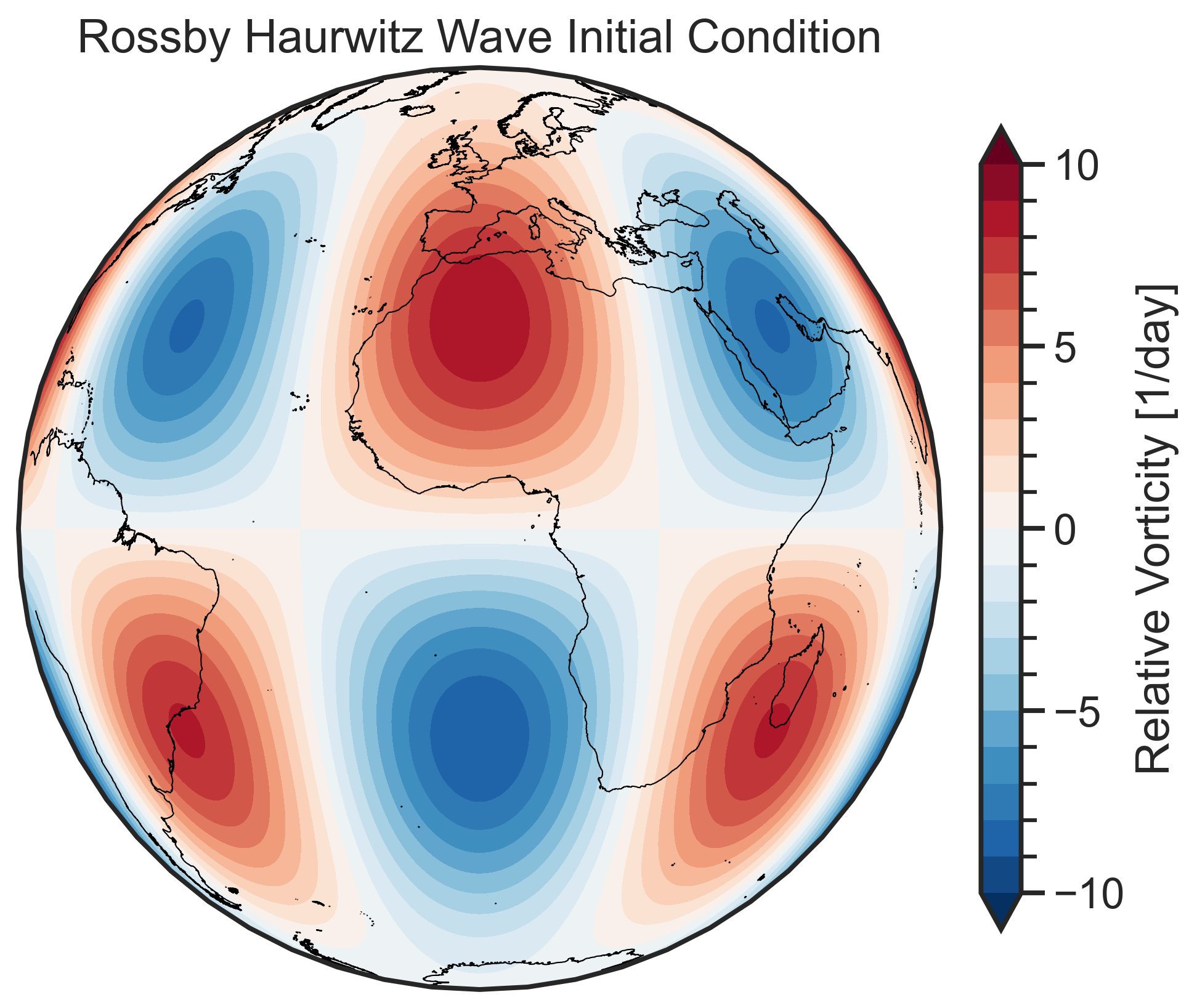}
    \end{subfigure}
    \begin{subfigure}{0.45\textwidth}
        \includegraphics[width=\linewidth]{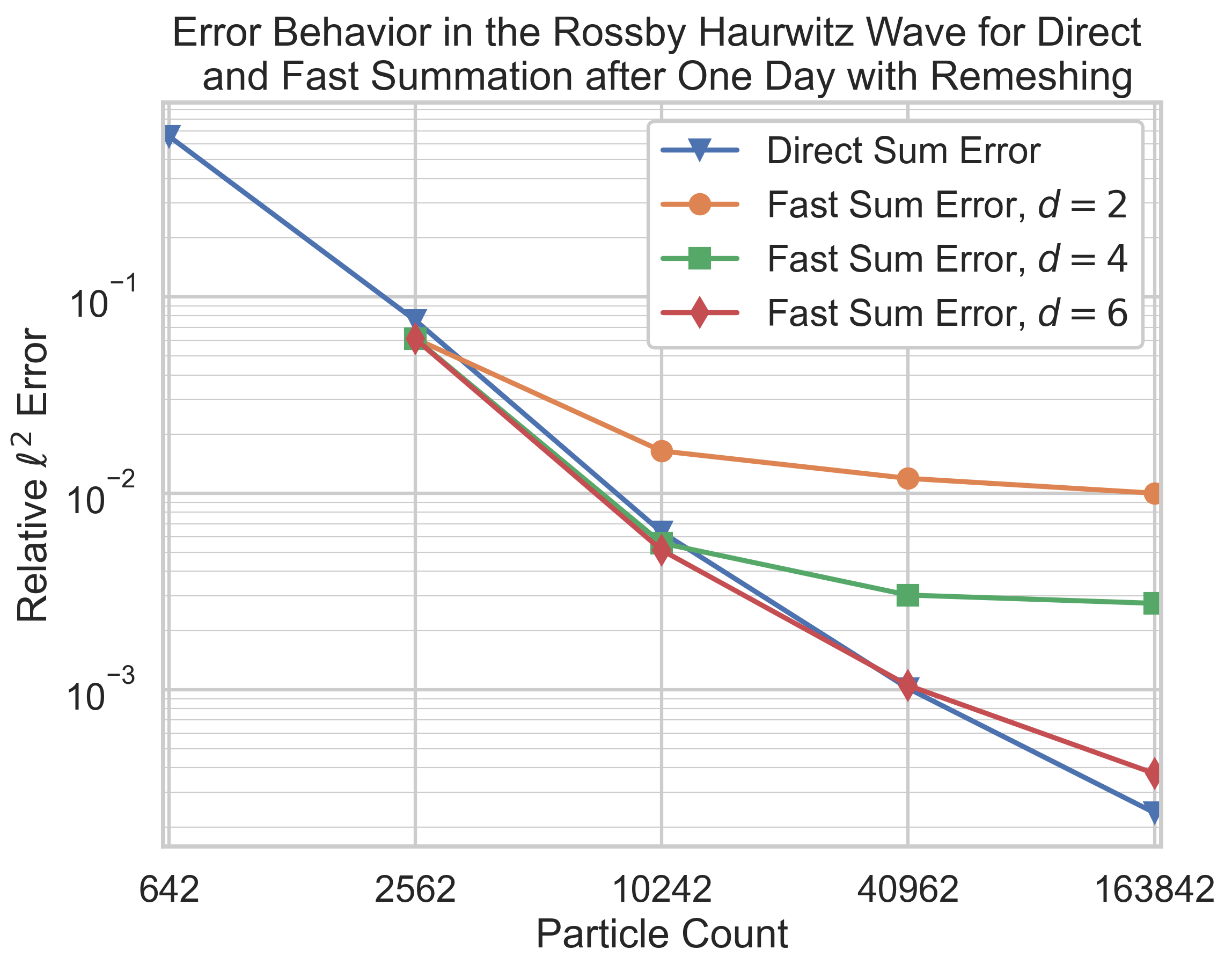}
    \end{subfigure}
    \caption{On the left, we plot the Rossby Haurwitz wave initial condition with a wave number of $4$. On the right, we plot the relative error behavior as we increase the particle count from 642, corresponding to 8 degree resolution, to 163842, corresponding to 0.5 degree resolution. We do not apply fast summation to the configuration with only 642 particles because there is very little room for speedup. In this case, error is coming from the discretization of the integral, interpolation in the remeshing, and fast summation. Here, the time step is small enough for time stepping error to be negligible. For these runs, we use $\theta=0.7$. }
    \label{fig:rhfig}
\end{figure}

\begin{figure}
	\centering
	\begin{subfigure}{0.23\textwidth}
		\includegraphics[width=\linewidth]{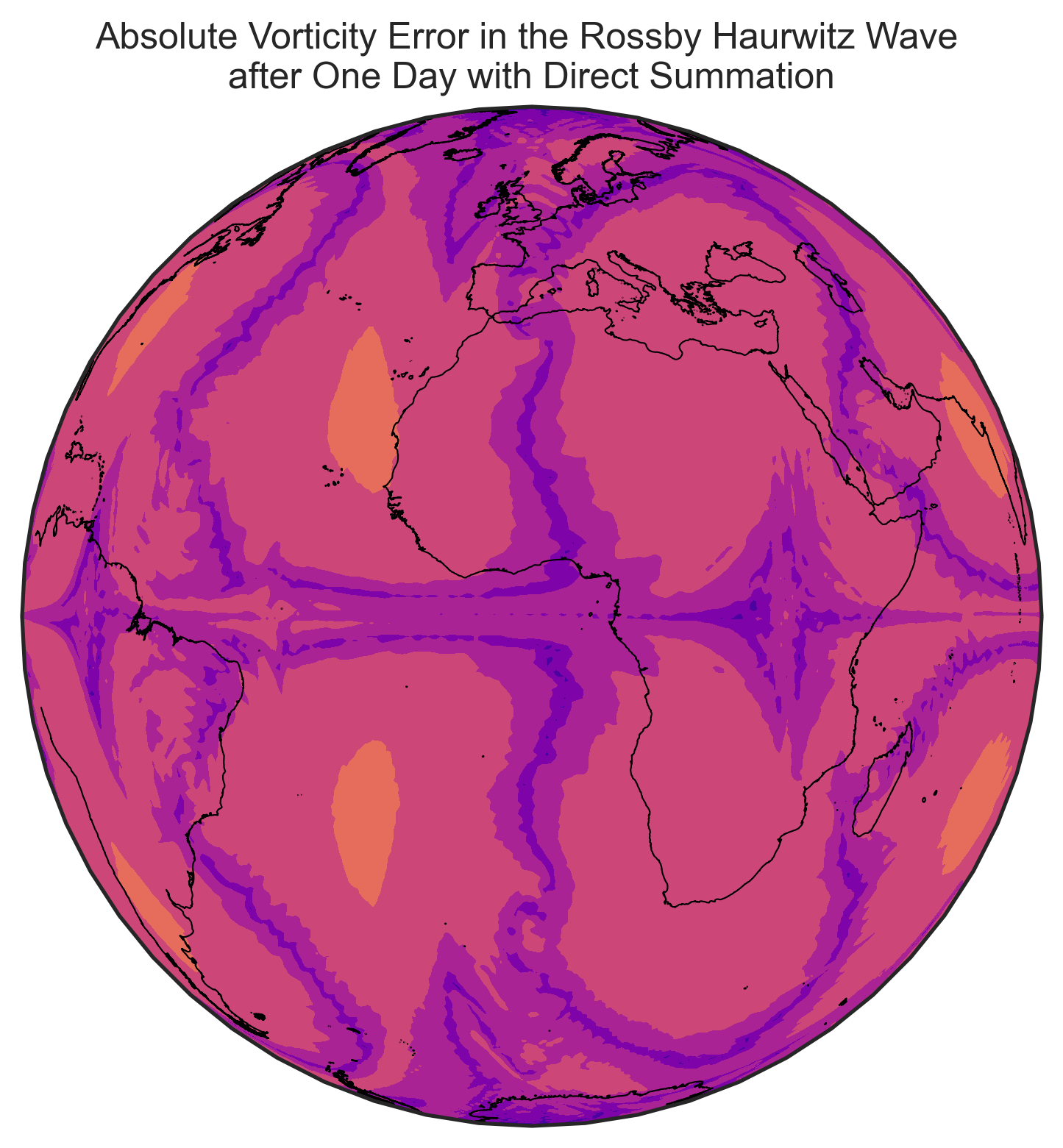}
	\end{subfigure}
	\begin{subfigure}{0.23\textwidth}
		\includegraphics[width=\linewidth]{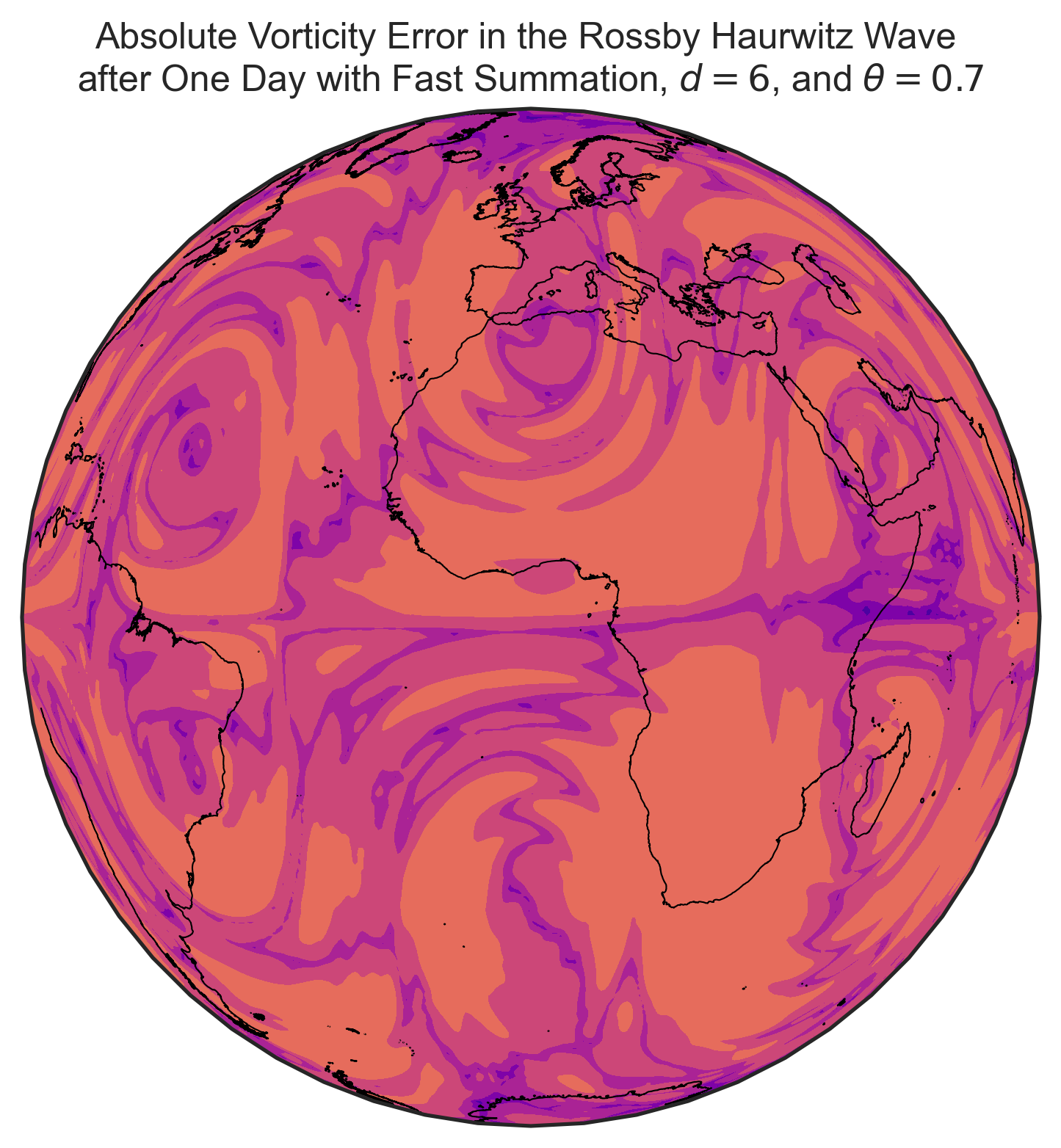}
	\end{subfigure}
	\begin{subfigure}{0.23\textwidth}
		\includegraphics[width=\linewidth]{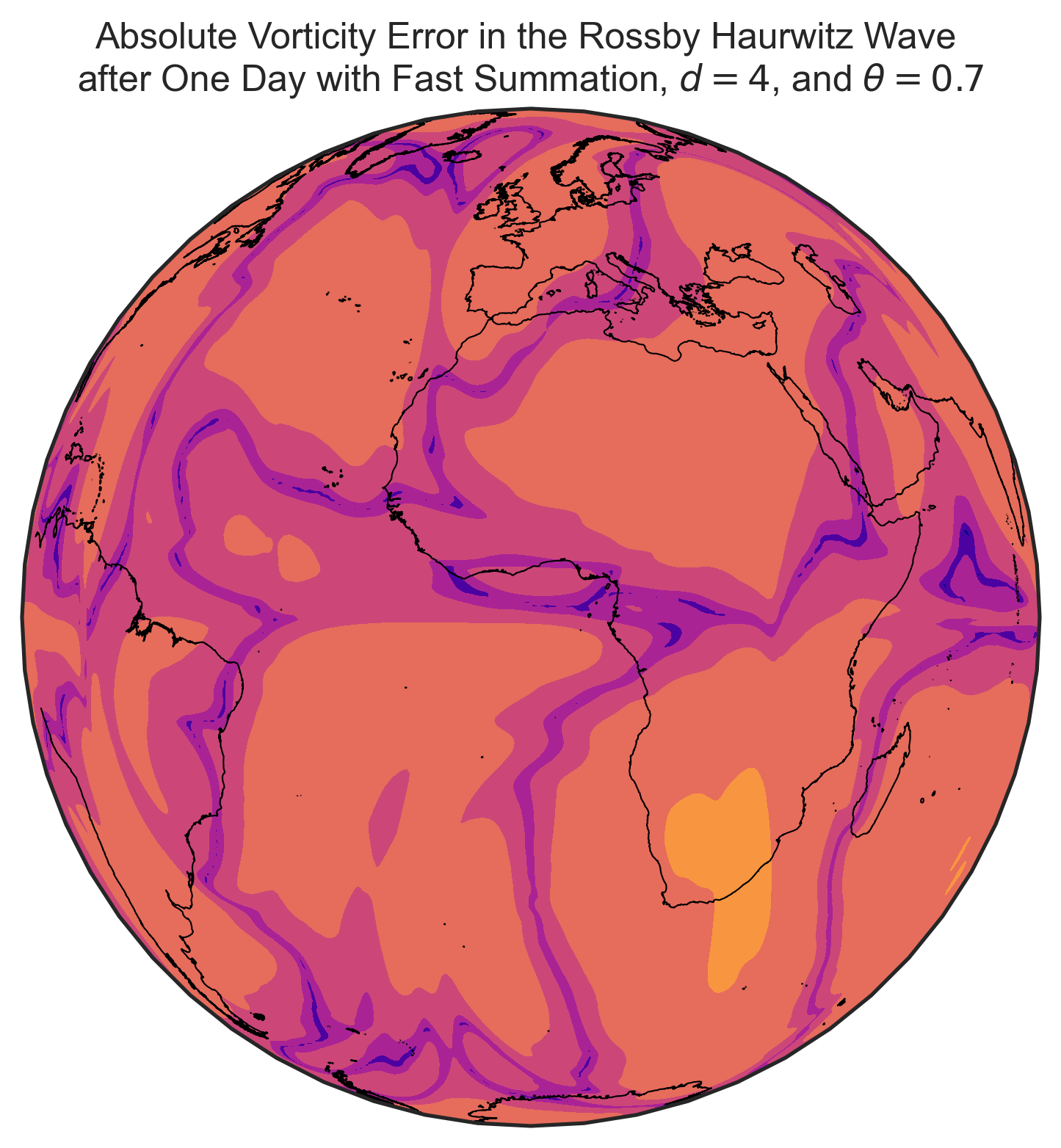}
	\end{subfigure}
	\begin{subfigure}{0.23\textwidth}
		\includegraphics[width=\linewidth]{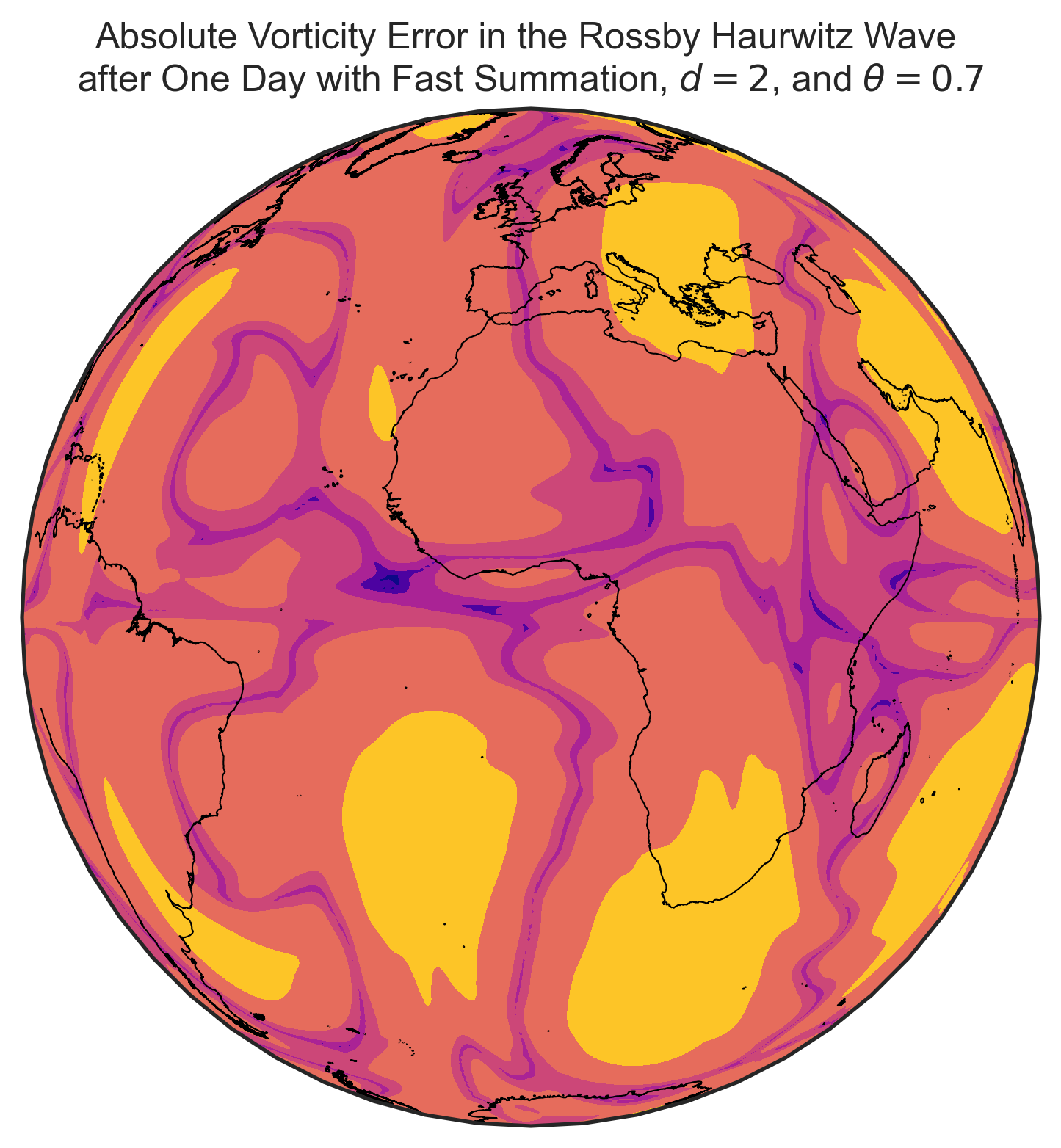}
	\end{subfigure}
	\begin{subfigure}{\textwidth}
		\centering
		\includegraphics[scale=0.25]{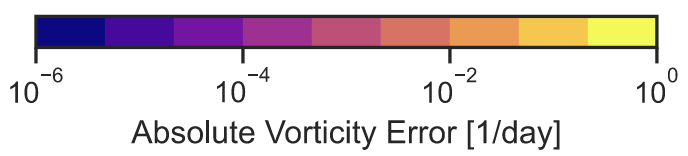}
	\end{subfigure}
	\caption{We plot the spatial structure of the vorticity error after one day. From left to right, we use direct summation, fast summation with interpolation degree $d=6$, fast summation with $d=4$, and fast summation with $d=2$. We plot the absolute error here because there are areas where the true vorticity is $0$ and hence the relative error is not defined everywhere. The MAC parameter used here is $\theta=0.7$.  }
	\label{fig:rherrors}
\end{figure}

The Rossby-Haurwitz waves are exact solutions of the BVE, letting us exactly measure errors. We follow the test case as detailed in \cite{williamson1992standard}. The initial vorticity is given as 
\begin{equation}
    \zeta_0(\theta,\phi)=2w\sin\theta+(k^2+3k+2)(\sin\theta)(\cos^k\theta)(\cos k\phi)
\end{equation} where $k$ is the zonal wavenumber, $\theta$ is the latitude, $\phi$ is the longitude, and $w$ is a constant related to the speed of the wave as \begin{equation}v=\frac{k(3+k)w-2\Omega}{(1+k)(2+k)}\end{equation} where $v$ is the speed of the wave and $\Omega$ is the rotation rate of the Earth. We take $k=4$ and $v=0$ so that the vorticity at all times are given by the initial vorticity. With $k=4$ the initial vorticity is then 
\begin{equation}\label{eq:rhwave}
    \zeta_0(\theta,\phi)=\frac{2\pi}{7}\sin\theta+30\sin\theta\cos^4\theta\cos(4\phi)
\end{equation}
To compute the error, we define a relative $\ell^2$ error as 
\begin{equation}
    E_N=\sqrt{\frac{\sum_{i=1}^N(\zeta_i-\zeta_0(\mathbf{x}_i))^2A_i}{\sum_{i=1}^N\zeta_0(\mathbf{x}_i)^2A_i}}
\end{equation}
and we note that we can define a relative $\ell^{\infty}$ error in an analogous way. In Fig.~\ref{fig:rhfig}, we present the Rossby Haurwitz wave initial condition, as well as the error behavior. Additionally, we plot the spatial structure of the error in Fig.~\ref{fig:rherrors}. We note that these errors are larger than the errors found in Table 1 of~\cite{bosler2014lagrangian}. The primary reason for this is that we are using a second order biquadratic remeshing scheme, as opposed to a sixth order cubic Hermite remeshing scheme. 

\begin{figure}
    \centering
    \begin{subfigure}{0.45\textwidth}
        \includegraphics[width=\linewidth]{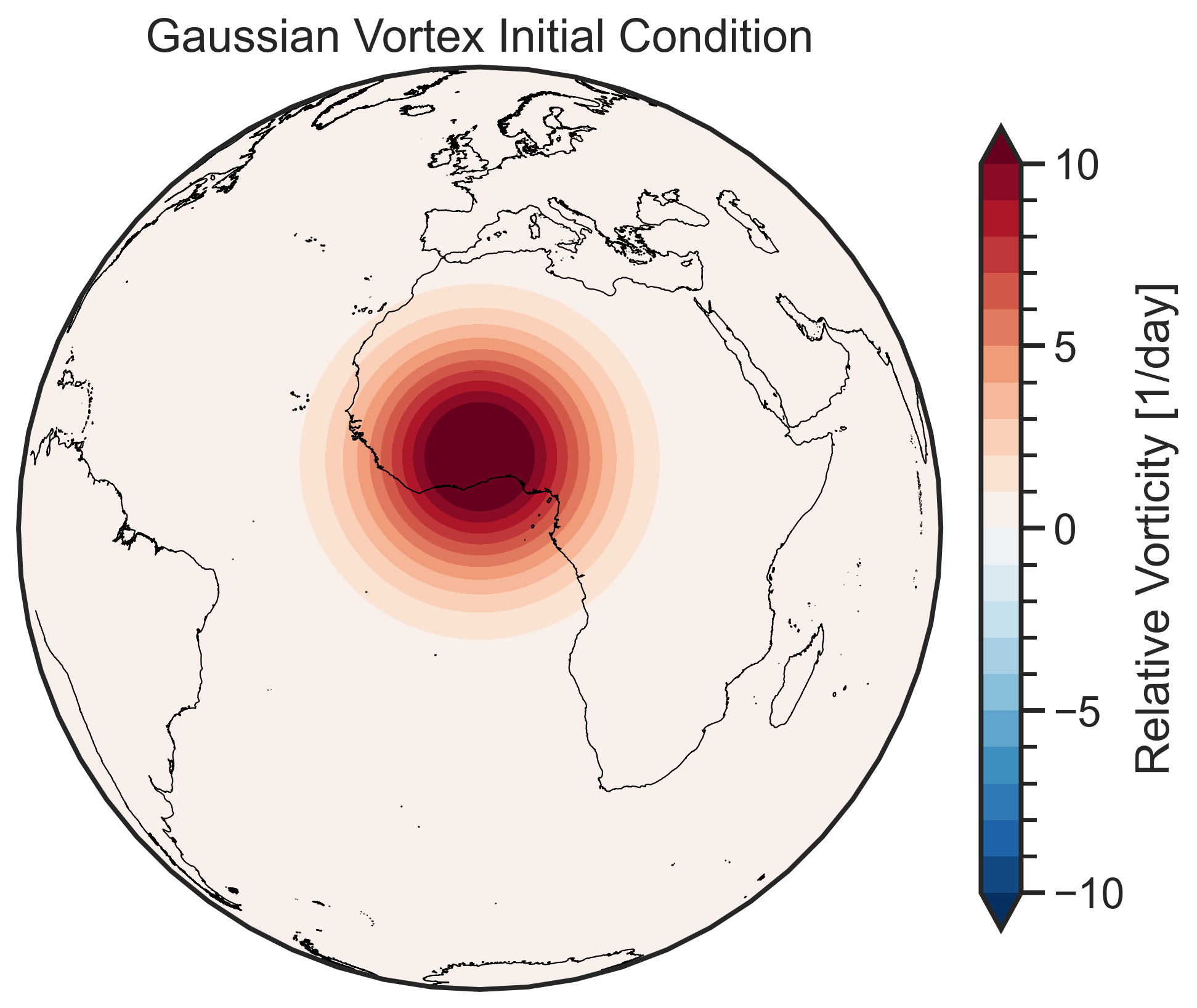}
        \caption{Gaussian Vortex initial condition}
        \label{fig:gv0}
    \end{subfigure}
    \begin{subfigure}{0.45\textwidth}
        \includegraphics[width=\linewidth]{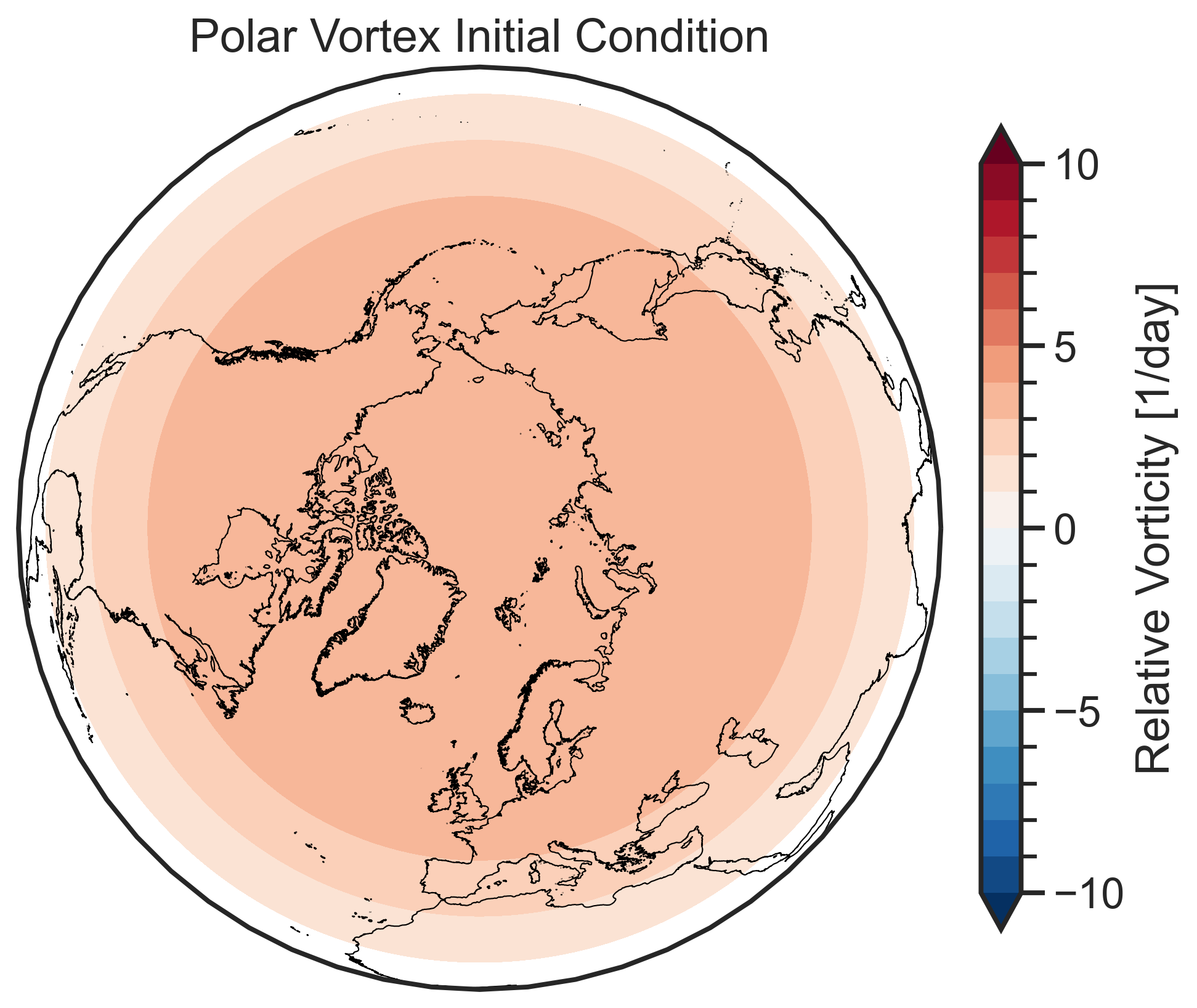}
        \caption{Polar Vortex initial condition}
        \label{fig:pv0}
    \end{subfigure}
    \caption{Gaussian and polar vortex initial conditions}
    \label{fig:pvgv0}
\end{figure}

\subsection{Gaussian Vortex}
\label{sec:results:gv}

We use the Gaussian vortex test case in order to assess the efficacy of our AMR scheme, as described in \ref{sec:amr:appendix}. In this case the initial vorticity is given as 
\begin{equation}
    \zeta_0=4\pi e^{-16\norm{\mathbf{x}-\mathbf{x}_c}^2}+C\,\mathrm{day}^{-1}
\end{equation} where $C$ is a constant chosen such that the total initial vorticity is $0$ and $\mathbf{x}_c$ is the center of the vortex, which we take to be at $\phi=0$ and $\theta=\frac{\pi}{20}$. We plot this initial condition in Fig.~\ref{fig:gv0}. In this case, we have no exact solution, but we can validate our results by comparing with those found in~\cite{bosler2014lagrangian}. 

\begin{figure}
    \centering
    \begin{subfigure}{0.20\textwidth}
        \includegraphics[width=\linewidth]{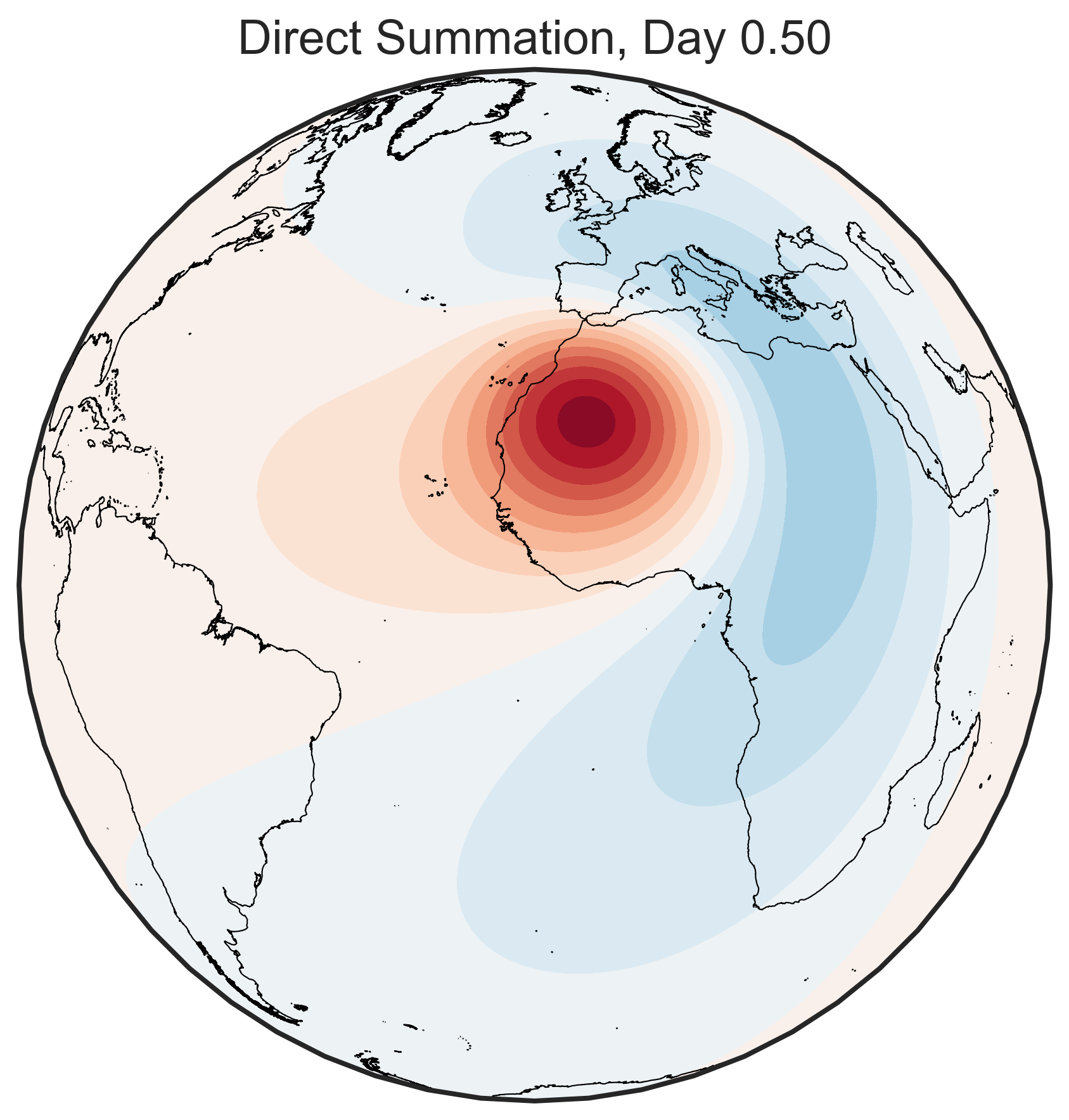}
    \end{subfigure}
    \hspace{0.025\textwidth}
    \begin{subfigure}{0.20\textwidth}
        \includegraphics[width=\linewidth]{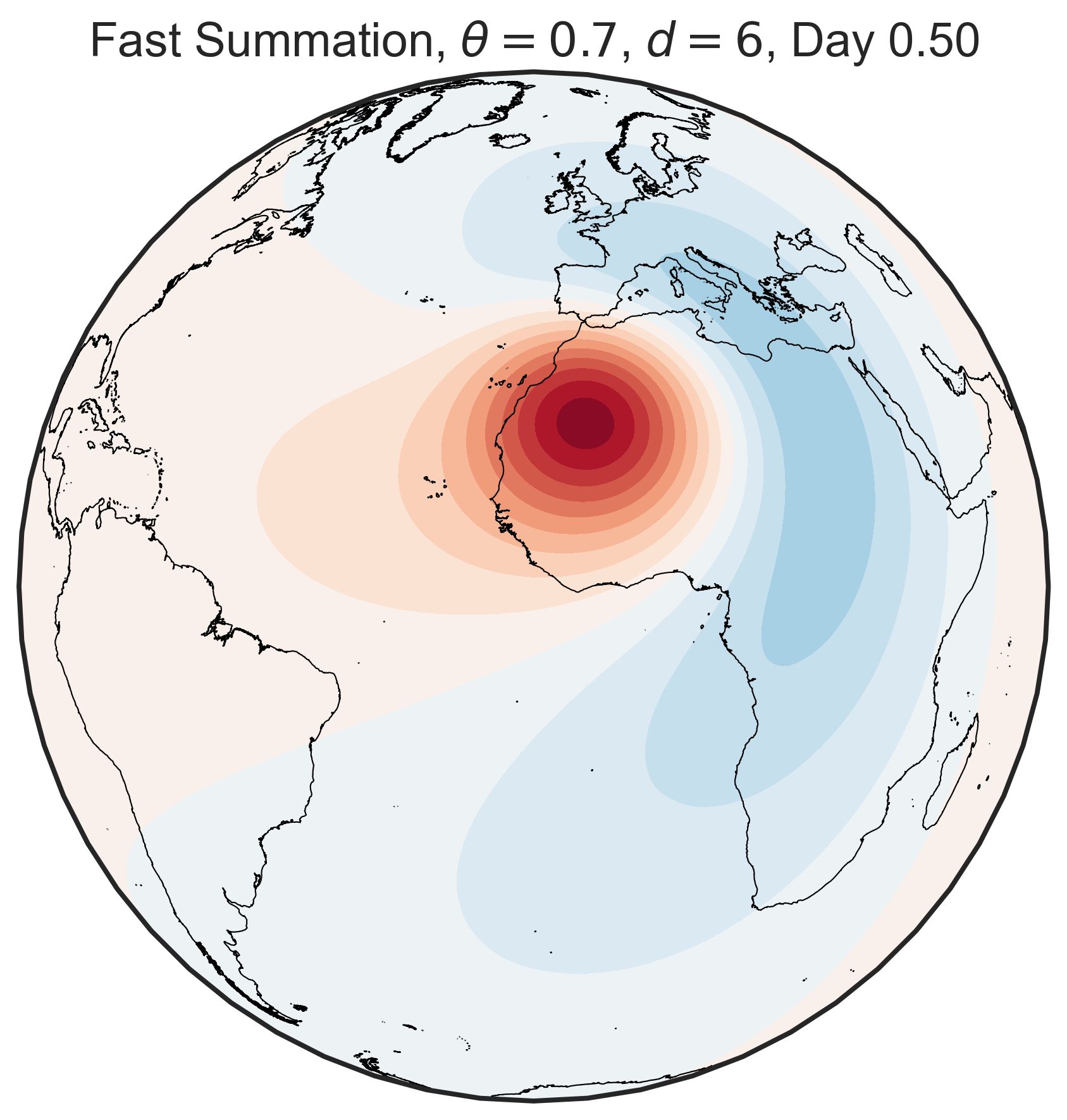}
    \end{subfigure}
    \hspace{0.025\textwidth}
    \begin{subfigure}{0.20\textwidth}
        \includegraphics[width=\linewidth]{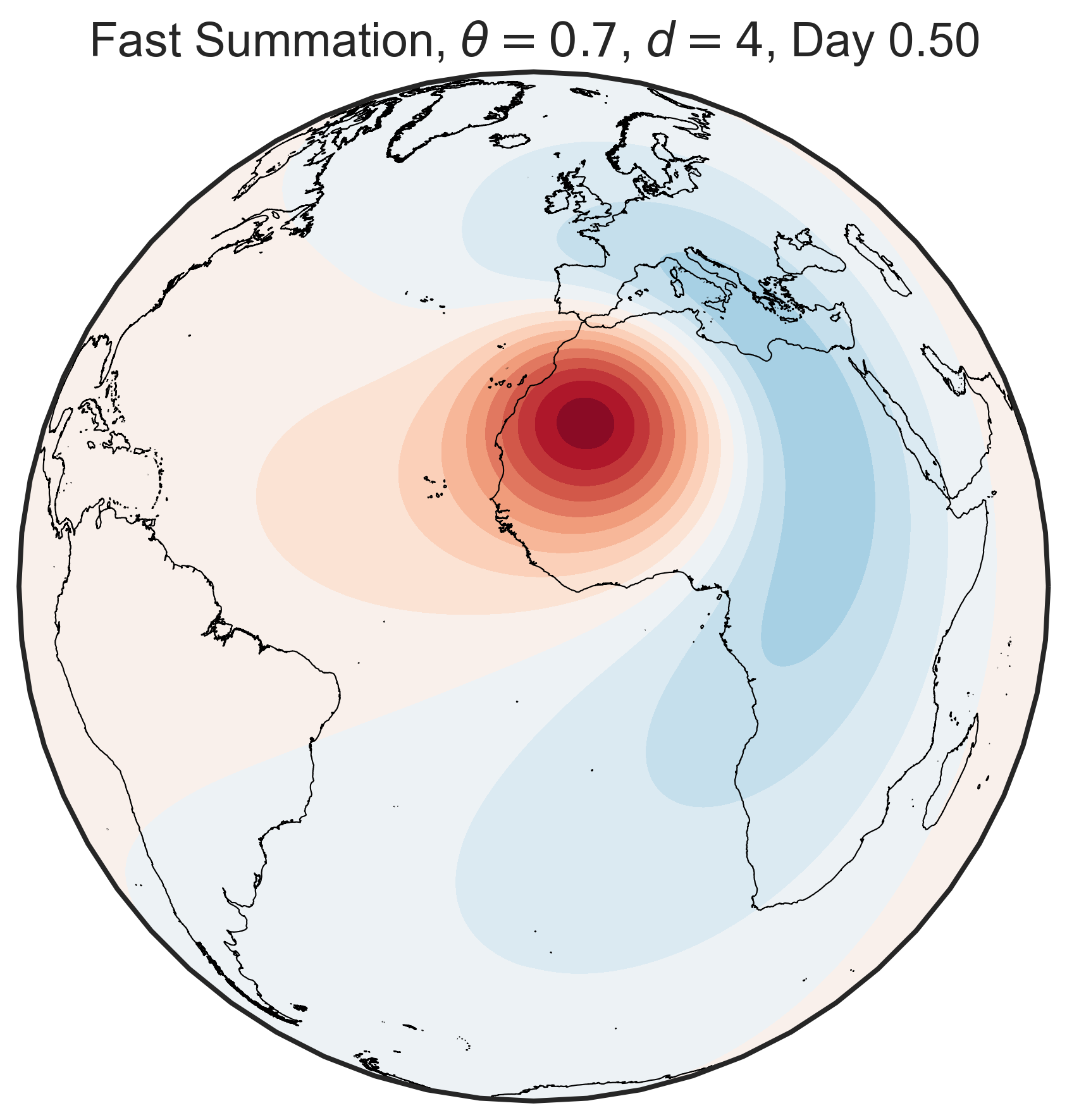}
    \end{subfigure}
    \hspace{0.025\textwidth}
    \begin{subfigure}{0.20\textwidth}
        \includegraphics[width=\linewidth]{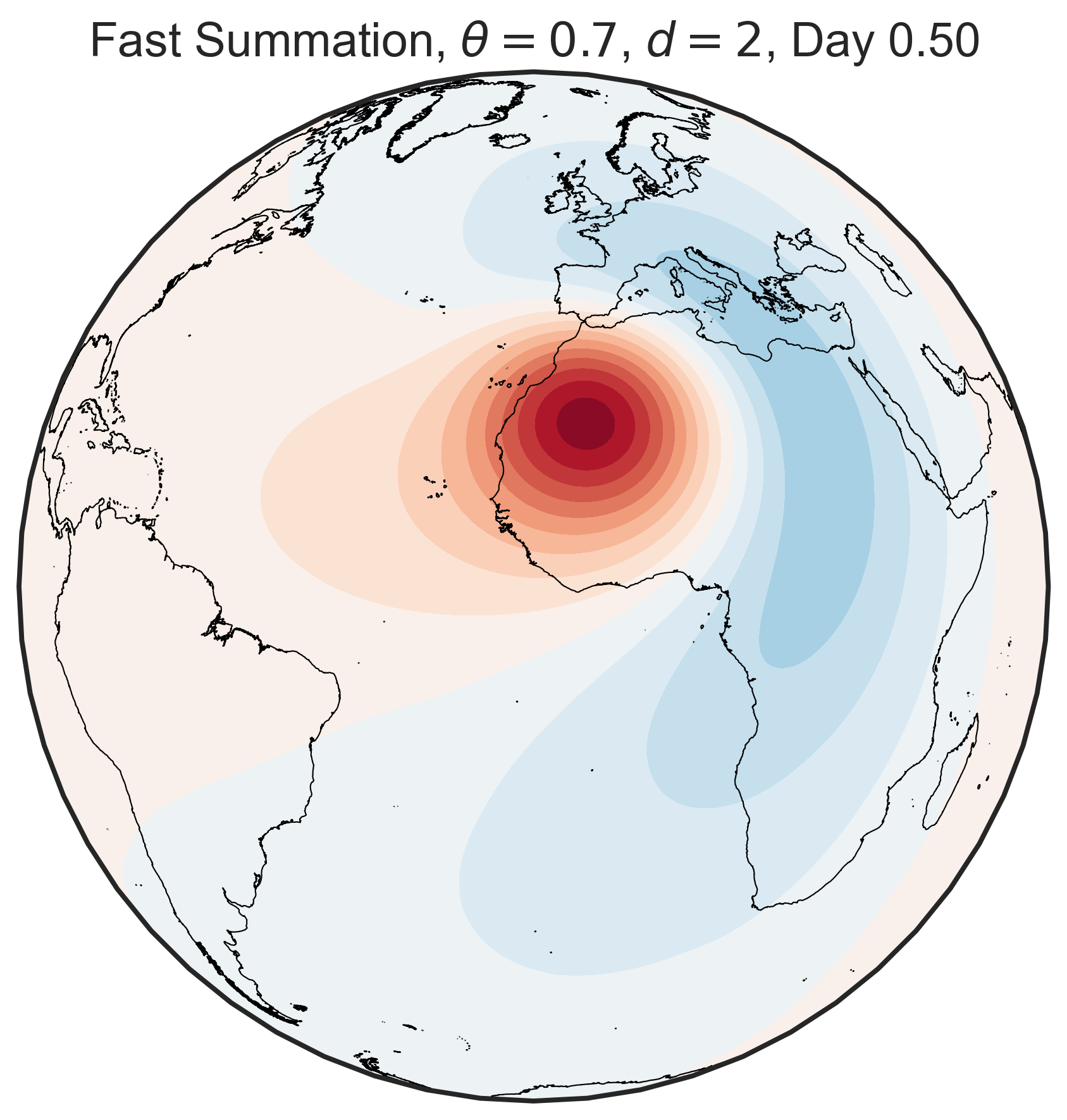}
    \end{subfigure}
    \begin{subfigure}{0.20\textwidth}
        \includegraphics[width=\linewidth]{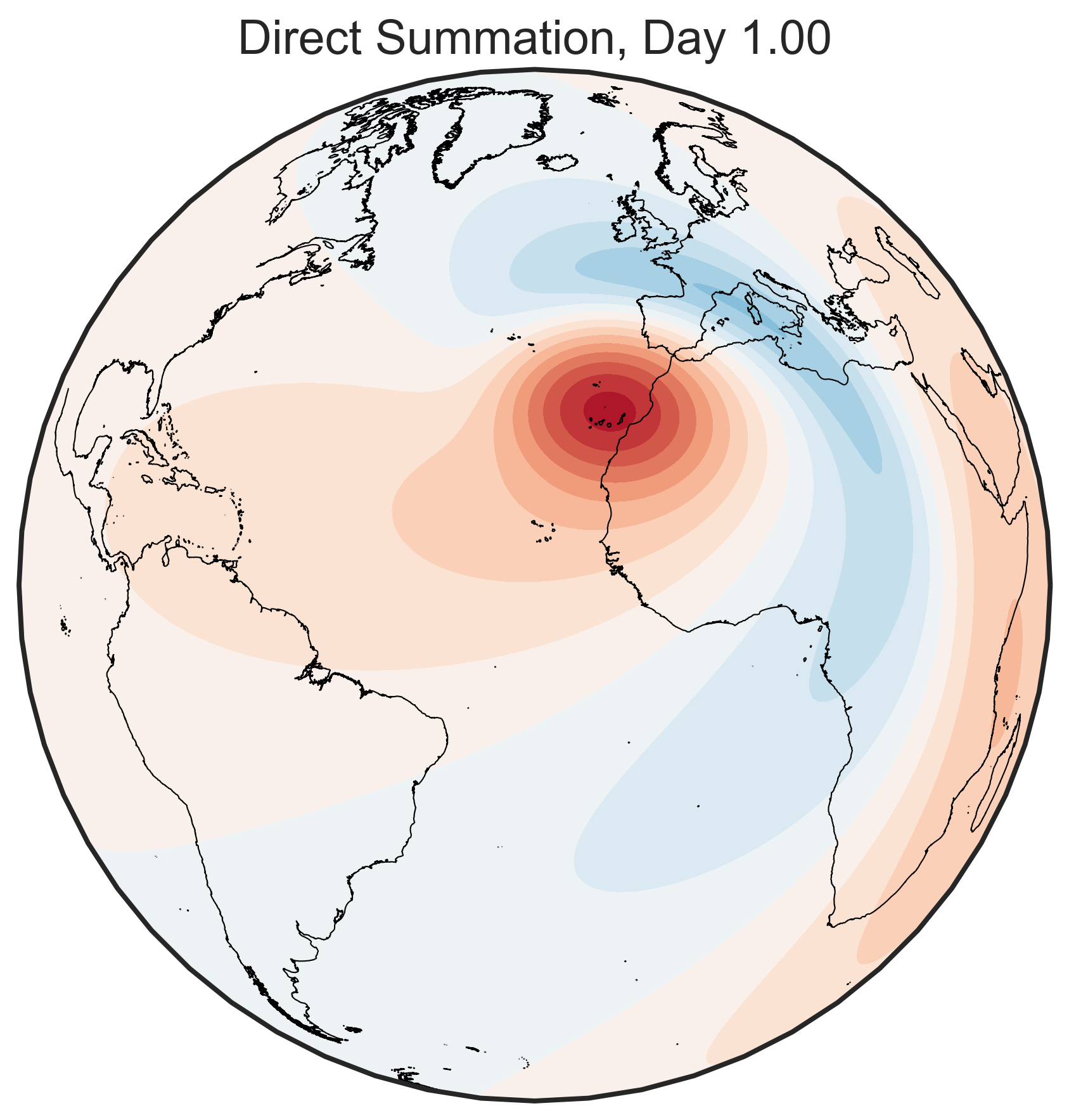}
    \end{subfigure}
    \hspace{0.025\textwidth}
    \begin{subfigure}{0.20\textwidth}
        \includegraphics[width=\linewidth]{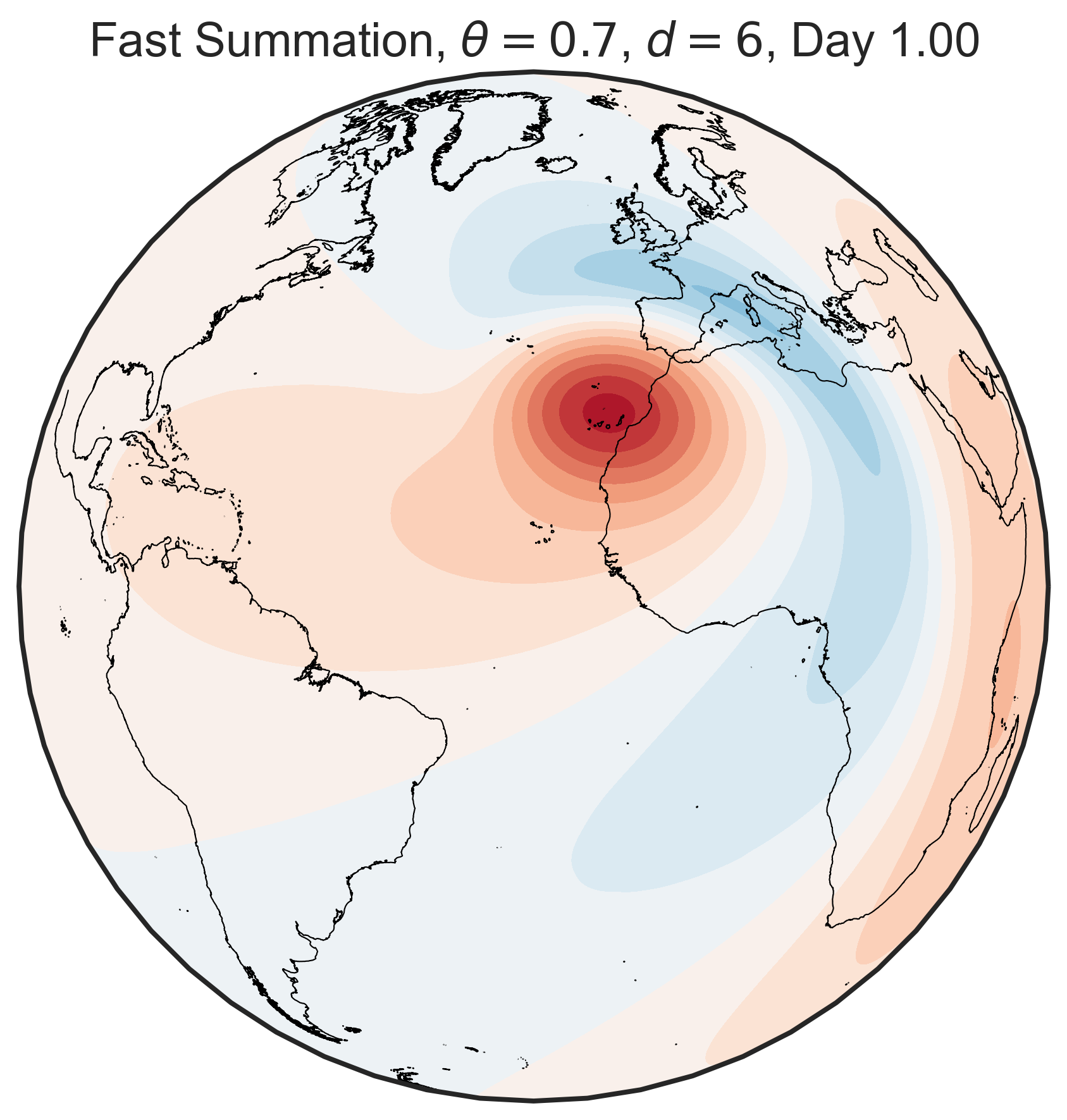}
    \end{subfigure}
    \hspace{0.025\textwidth}
    \begin{subfigure}{0.20\textwidth}
        \includegraphics[width=\linewidth]{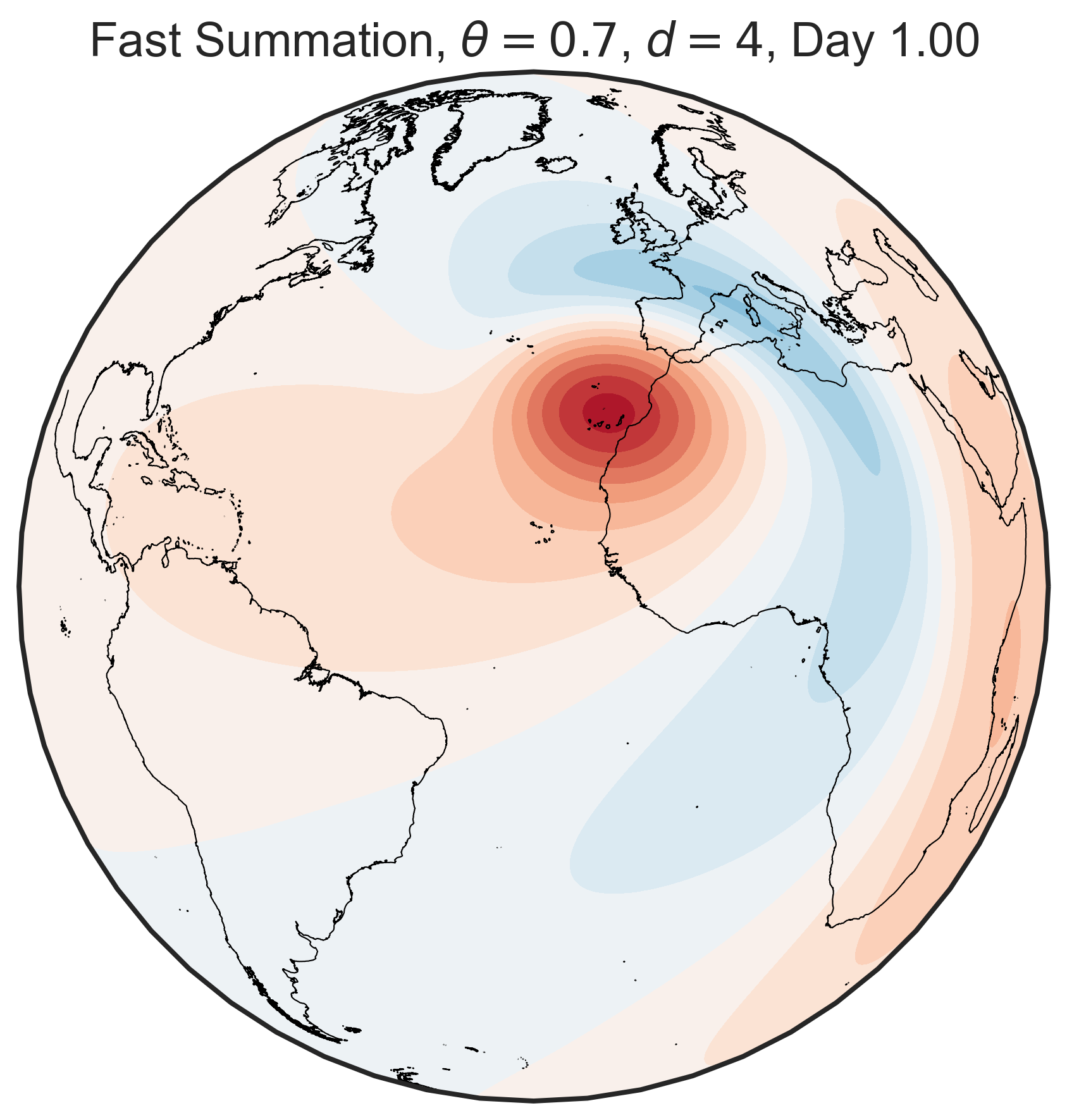}
    \end{subfigure}
    \hspace{0.025\textwidth}
    \begin{subfigure}{0.20\textwidth}
        \includegraphics[width=\linewidth]{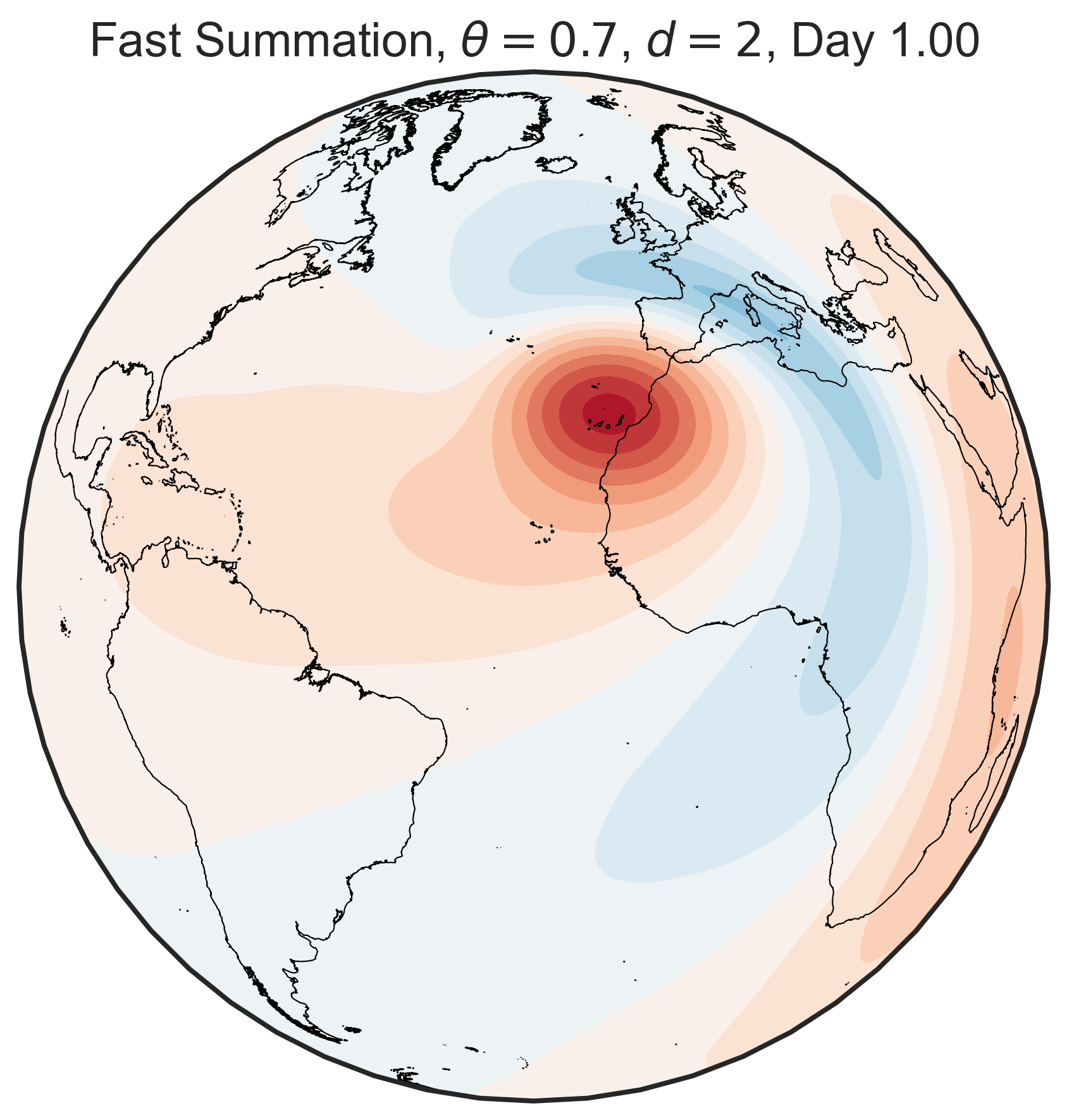}
    \end{subfigure}
    \begin{subfigure}{0.20\textwidth}
        \includegraphics[width=\linewidth]{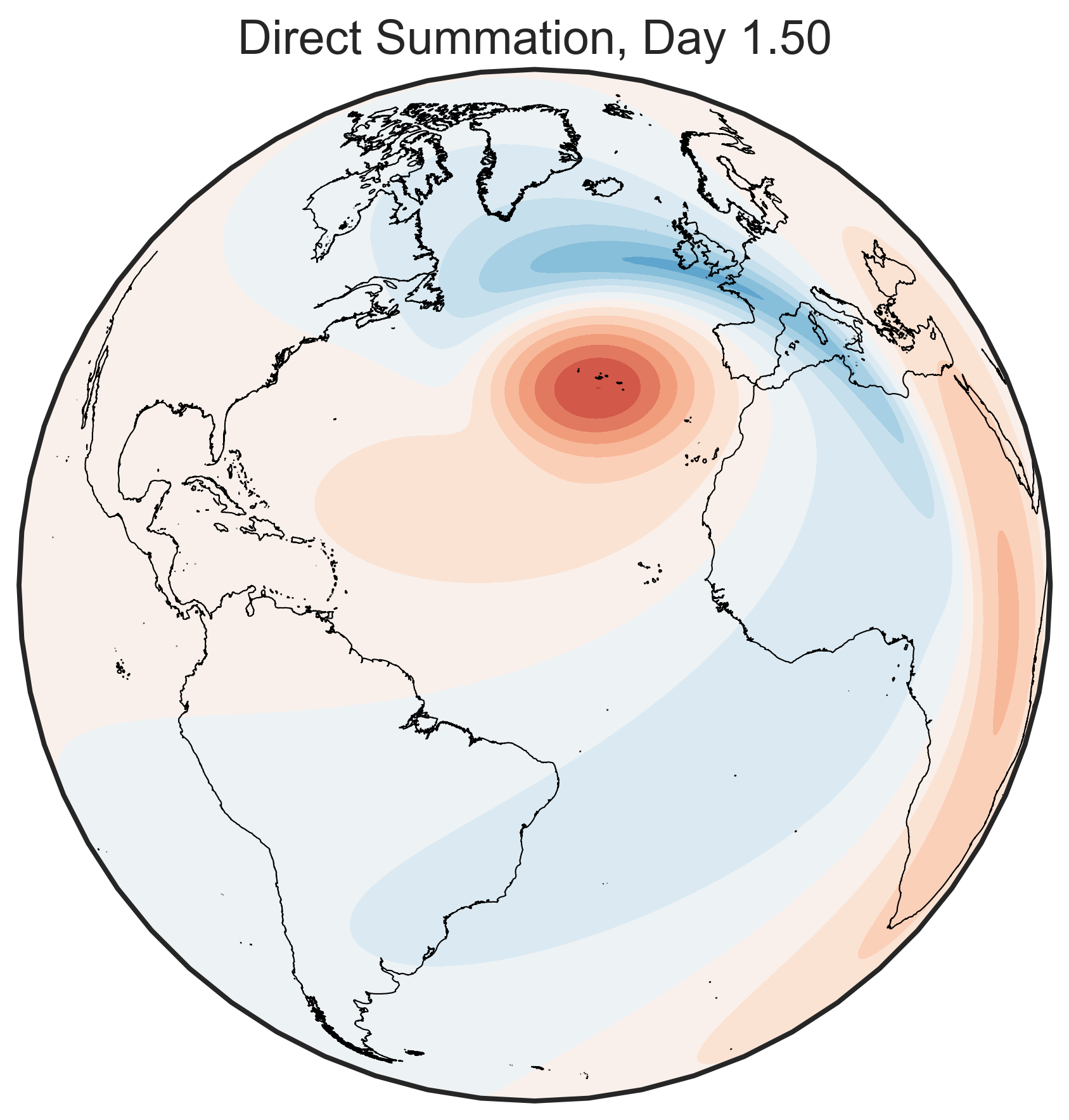}
    \end{subfigure}
    \hspace{0.025\textwidth}
    \begin{subfigure}{0.20\textwidth}
        \includegraphics[width=\linewidth]{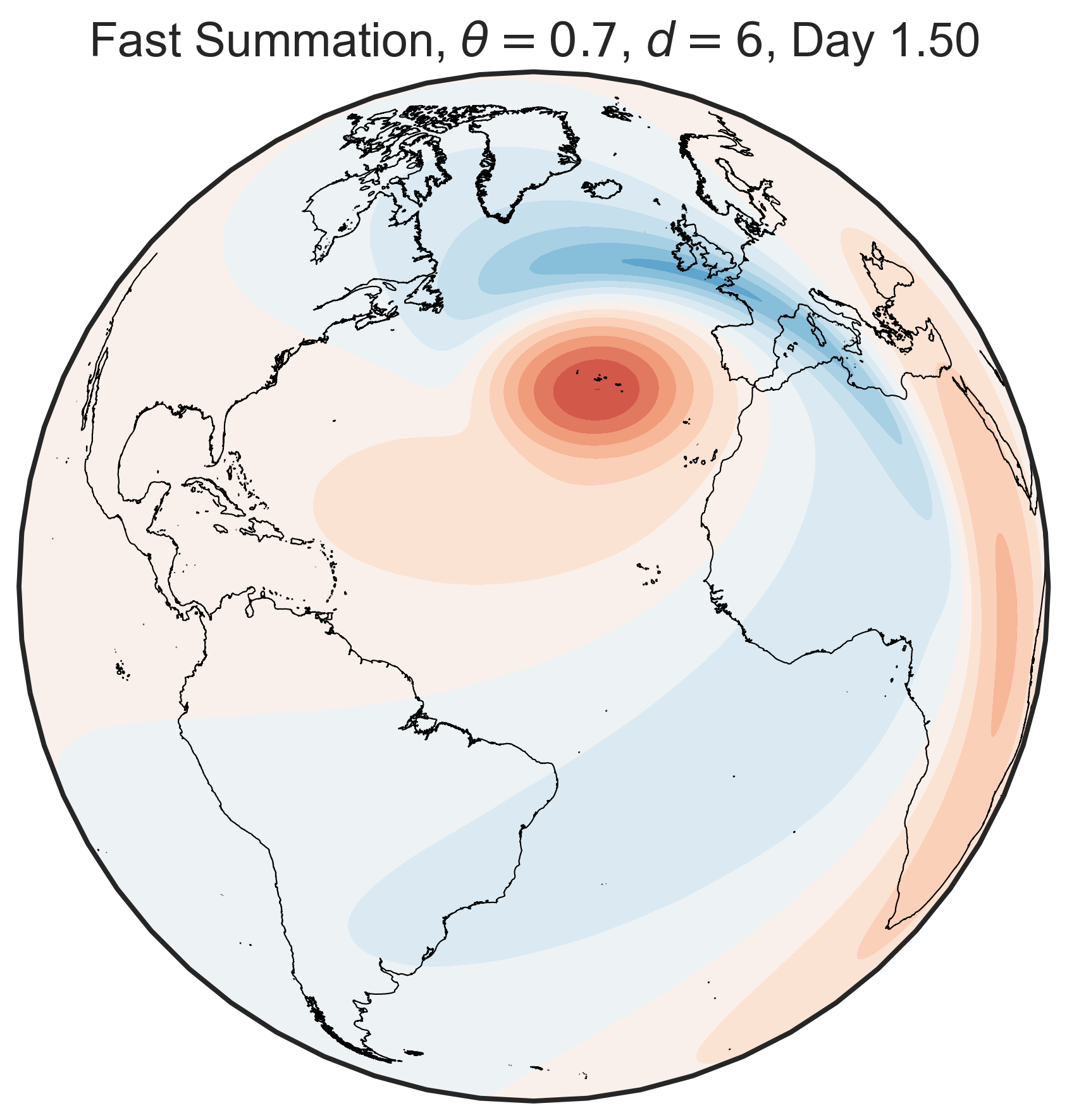}
    \end{subfigure}
    \hspace{0.025\textwidth}
    \begin{subfigure}{0.20\textwidth}
        \includegraphics[width=\linewidth]{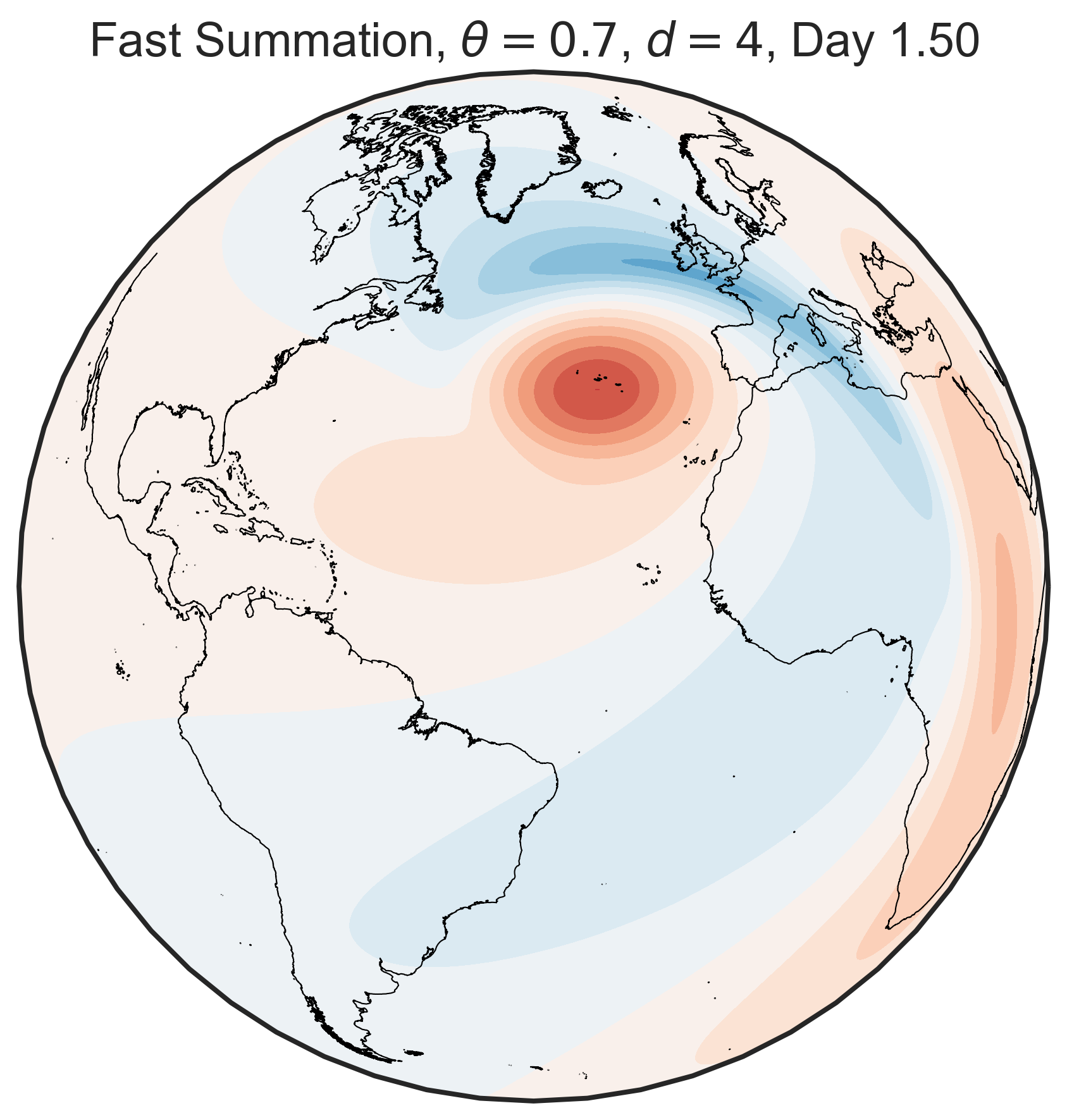}
    \end{subfigure}
    \hspace{0.025\textwidth}
    \begin{subfigure}{0.20\textwidth}
        \includegraphics[width=\linewidth]{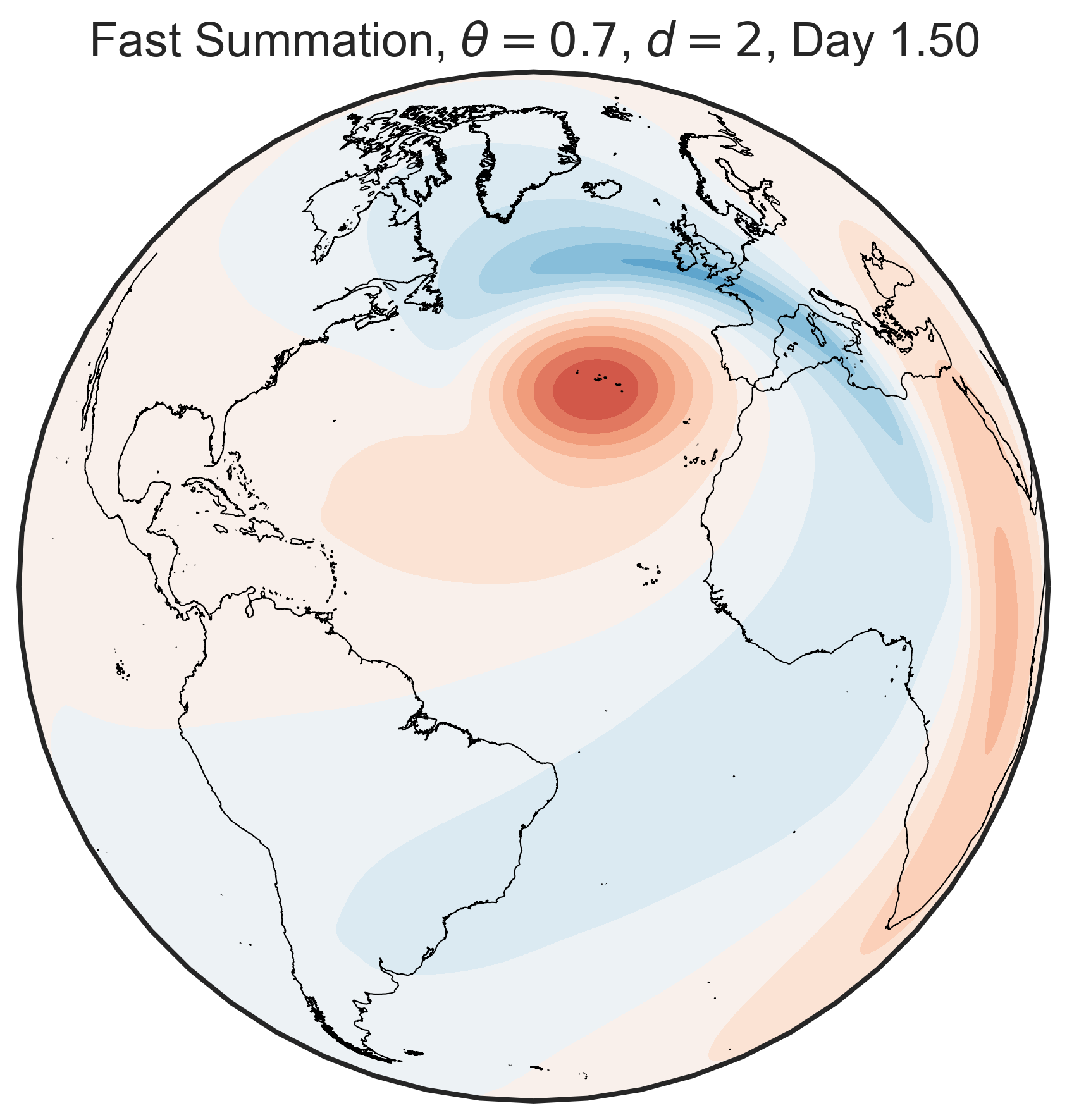}
    \end{subfigure}
    \begin{subfigure}{0.20\textwidth}
        \includegraphics[width=\linewidth]{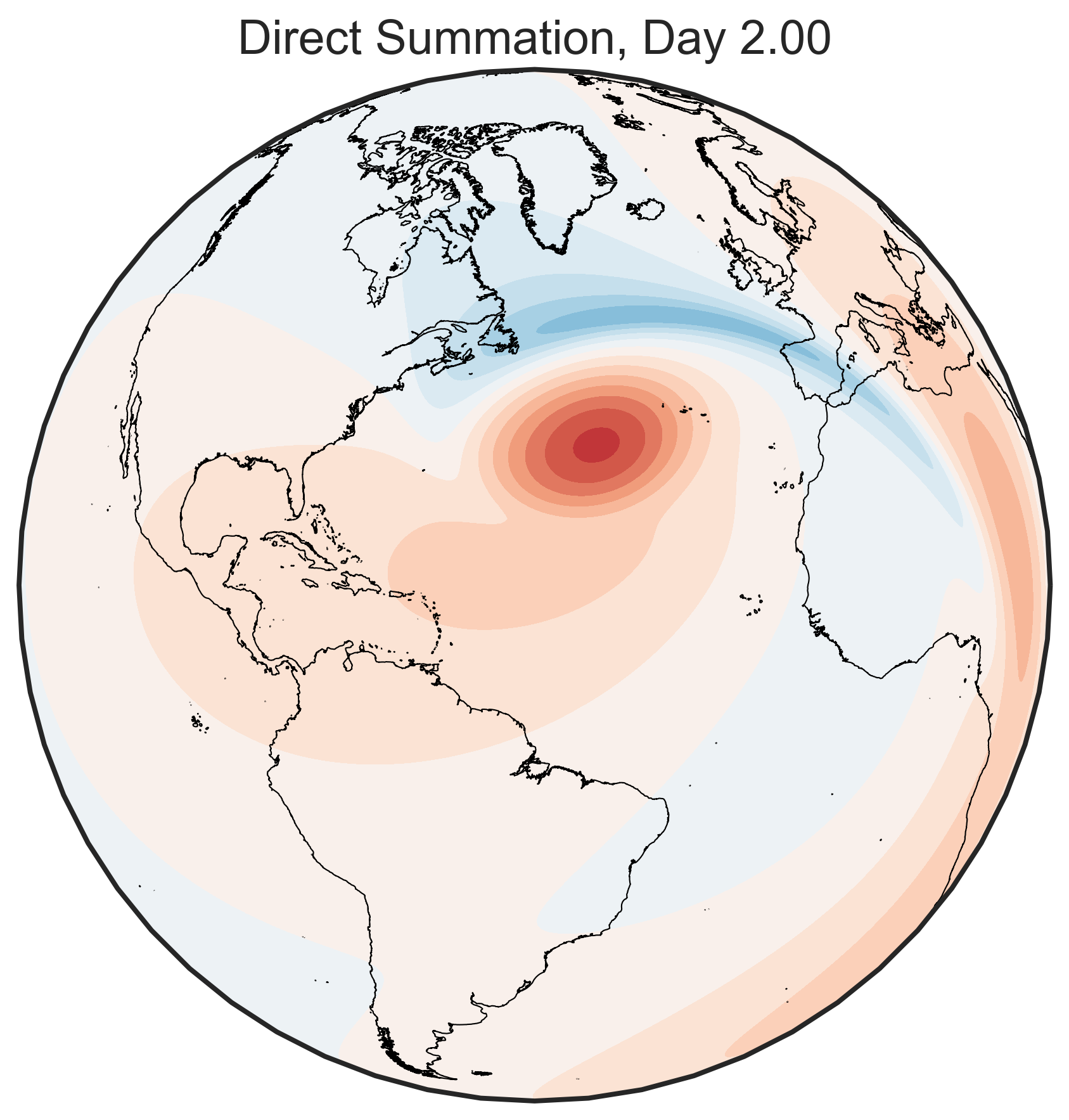}
    \end{subfigure}
    \hspace{0.025\textwidth}
    \begin{subfigure}{0.20\textwidth}
        \includegraphics[width=\linewidth]{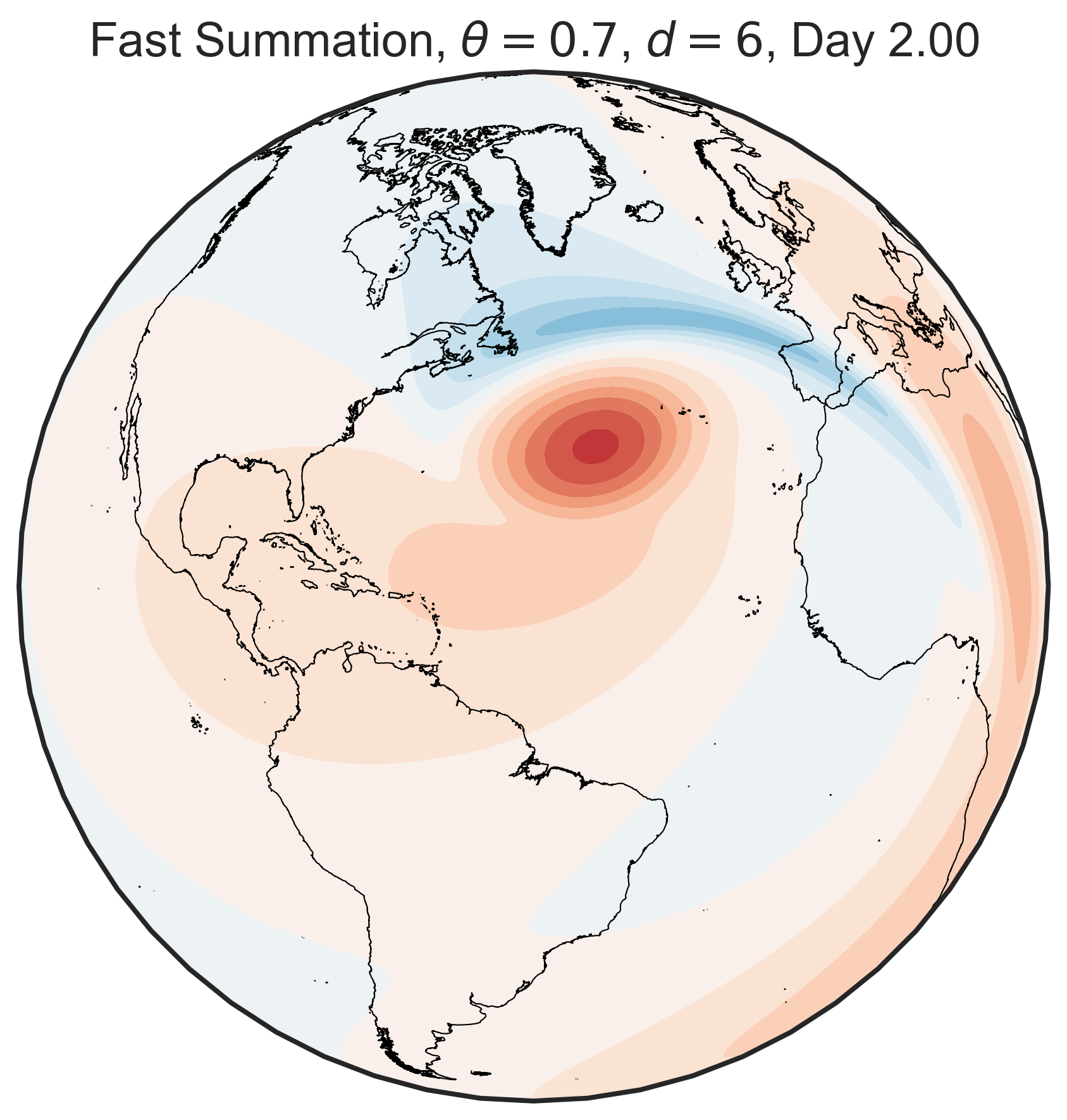}
    \end{subfigure}
    \hspace{0.025\textwidth}
    \begin{subfigure}{0.20\textwidth}
        \includegraphics[width=\linewidth]{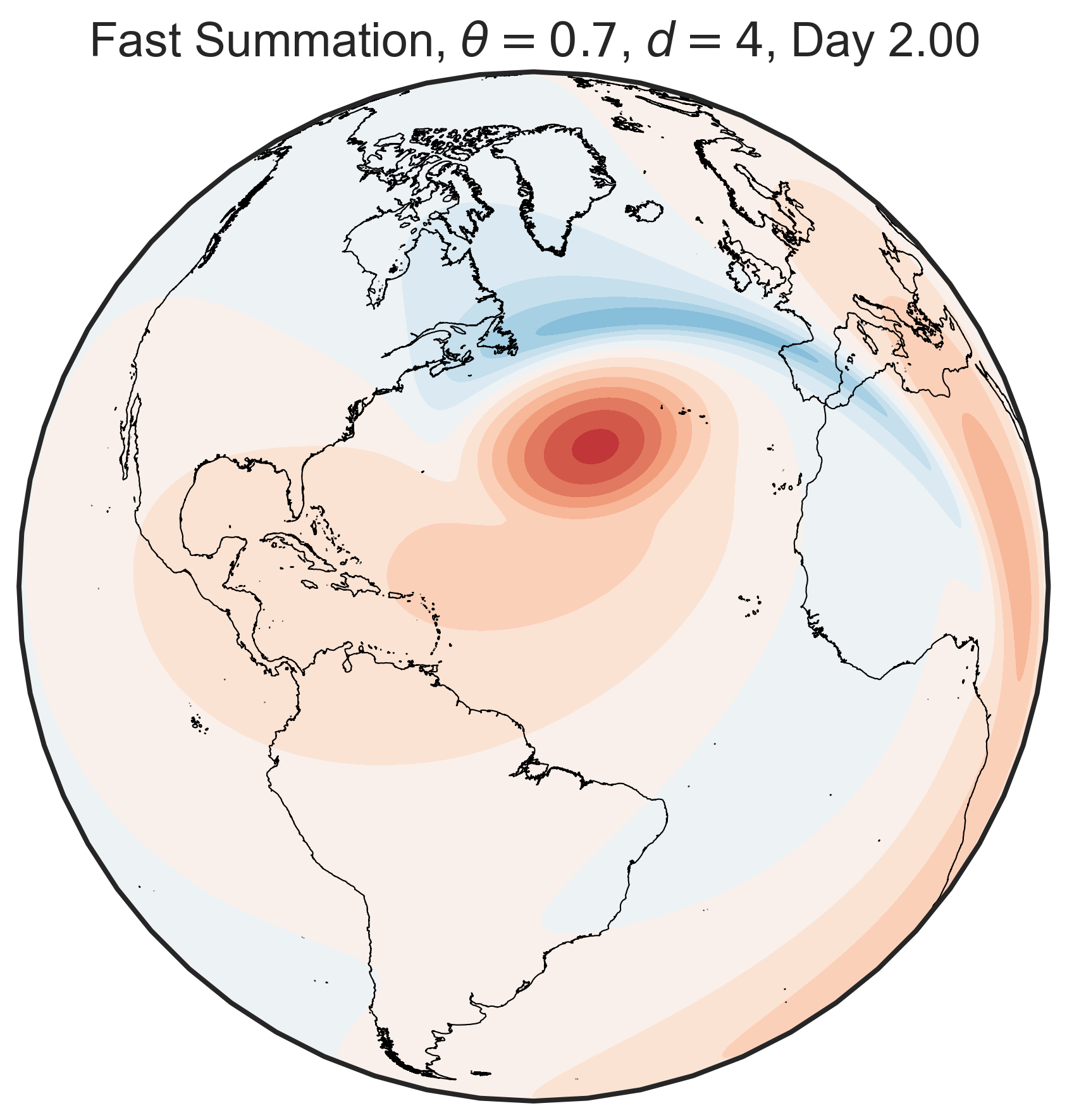}
    \end{subfigure}
    \hspace{0.025\textwidth}
    \begin{subfigure}{0.20\textwidth}
        \includegraphics[width=\linewidth]{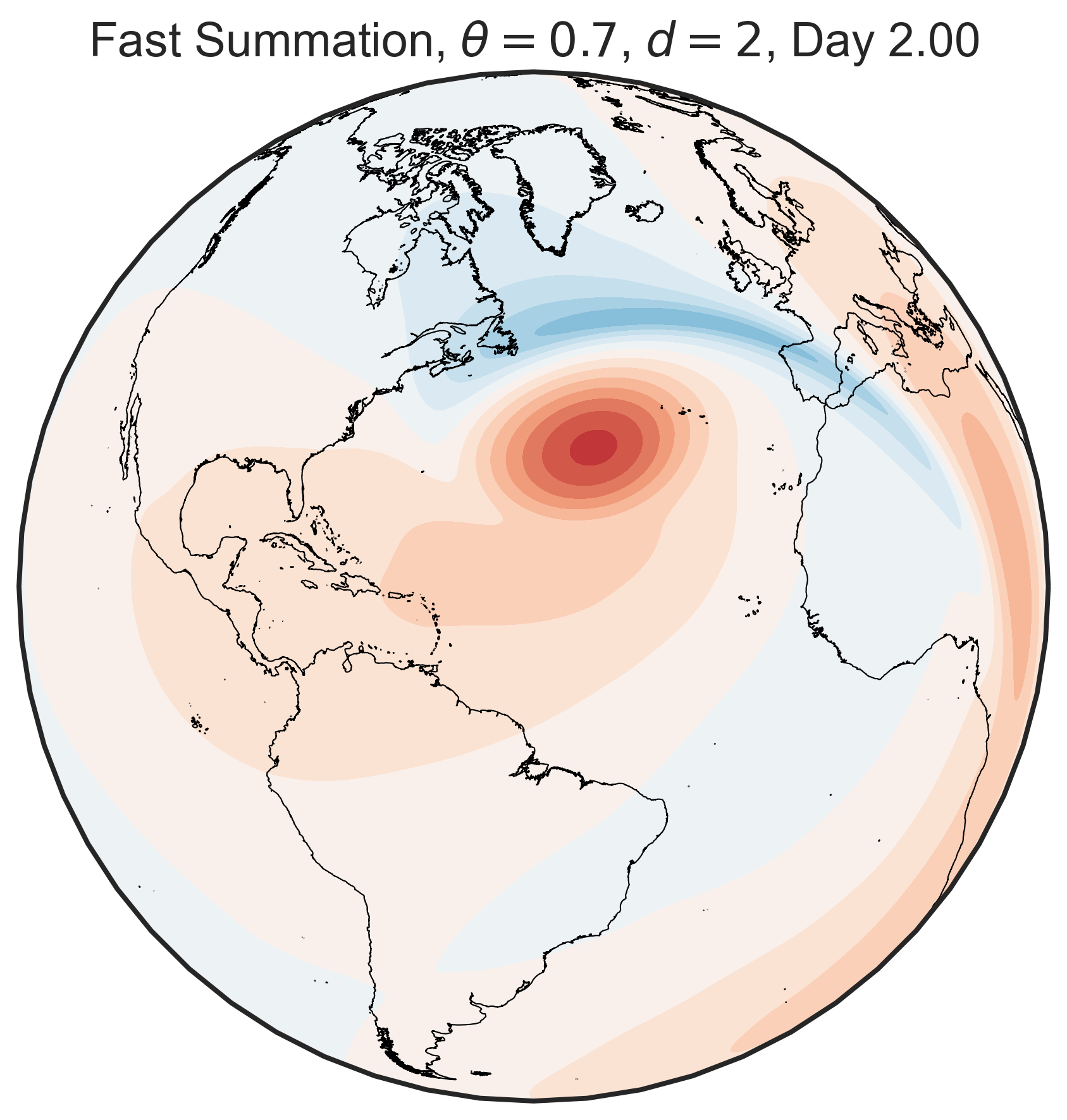}
    \end{subfigure}
    \begin{subfigure}{0.20\textwidth}
        \includegraphics[width=\linewidth]{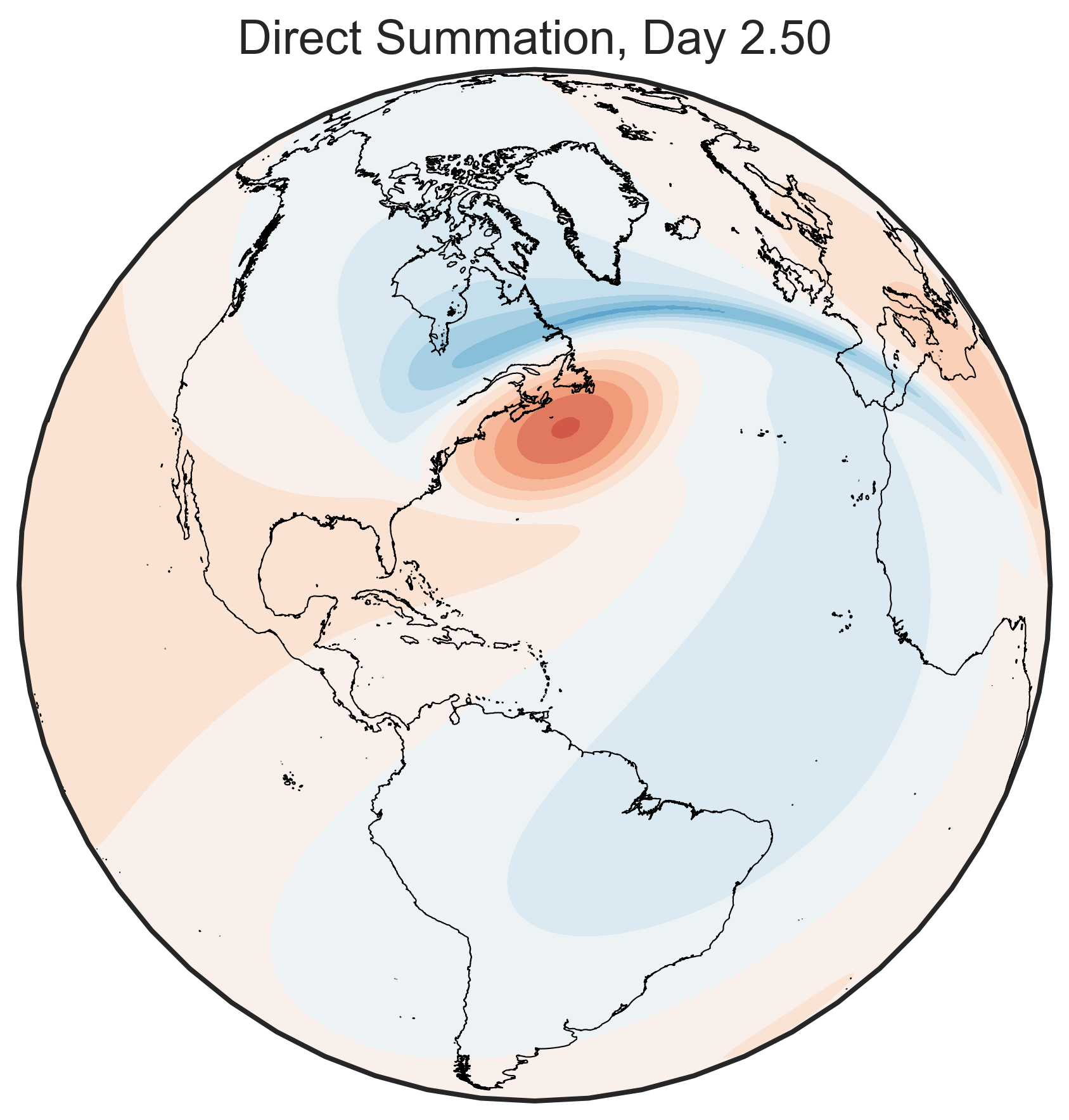}
    \end{subfigure}
    \hspace{0.025\textwidth}
    \begin{subfigure}{0.20\textwidth}
        \includegraphics[width=\linewidth]{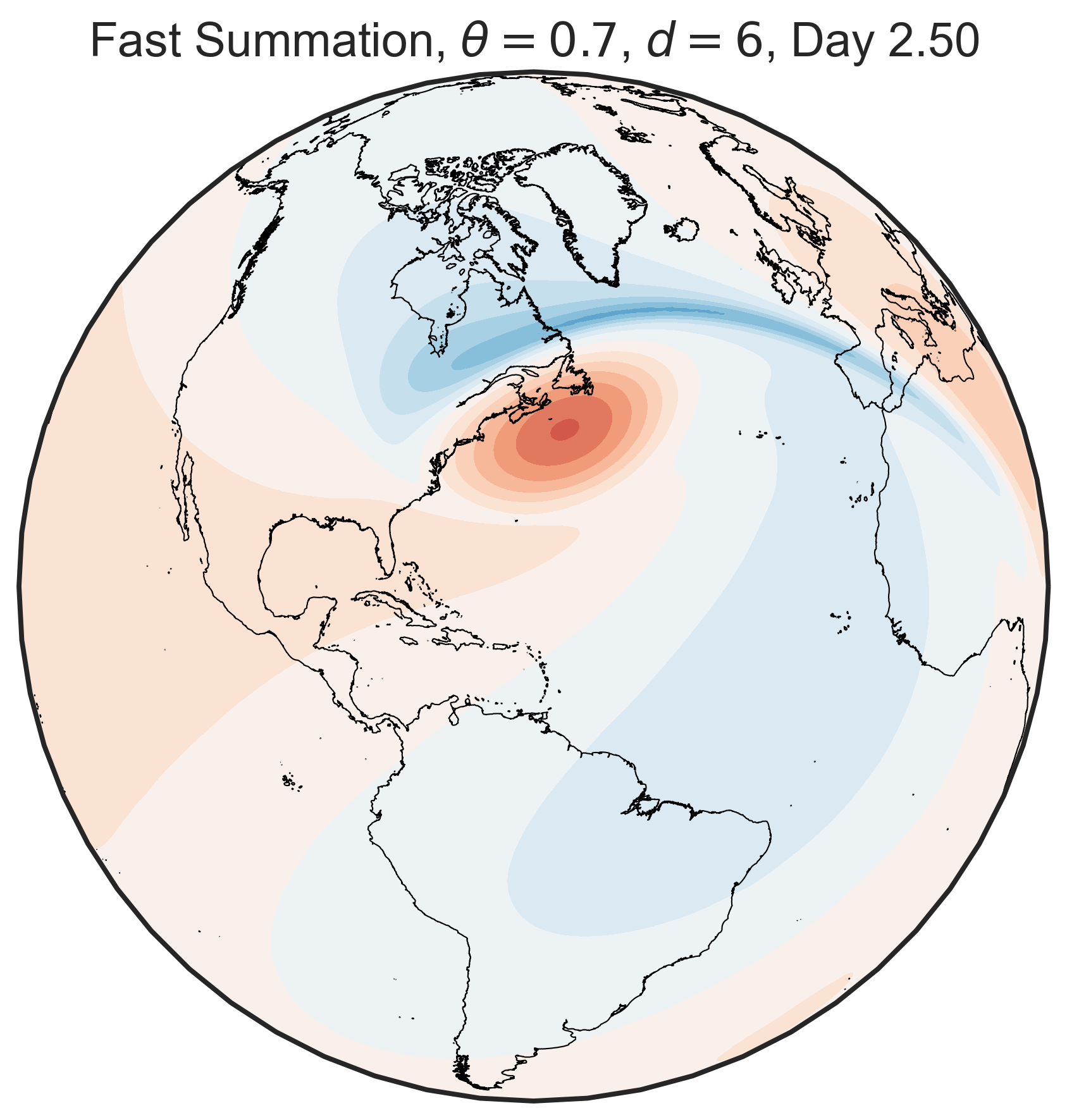}
    \end{subfigure}
    \hspace{0.025\textwidth}
    \begin{subfigure}{0.20\textwidth}
        \includegraphics[width=\linewidth]{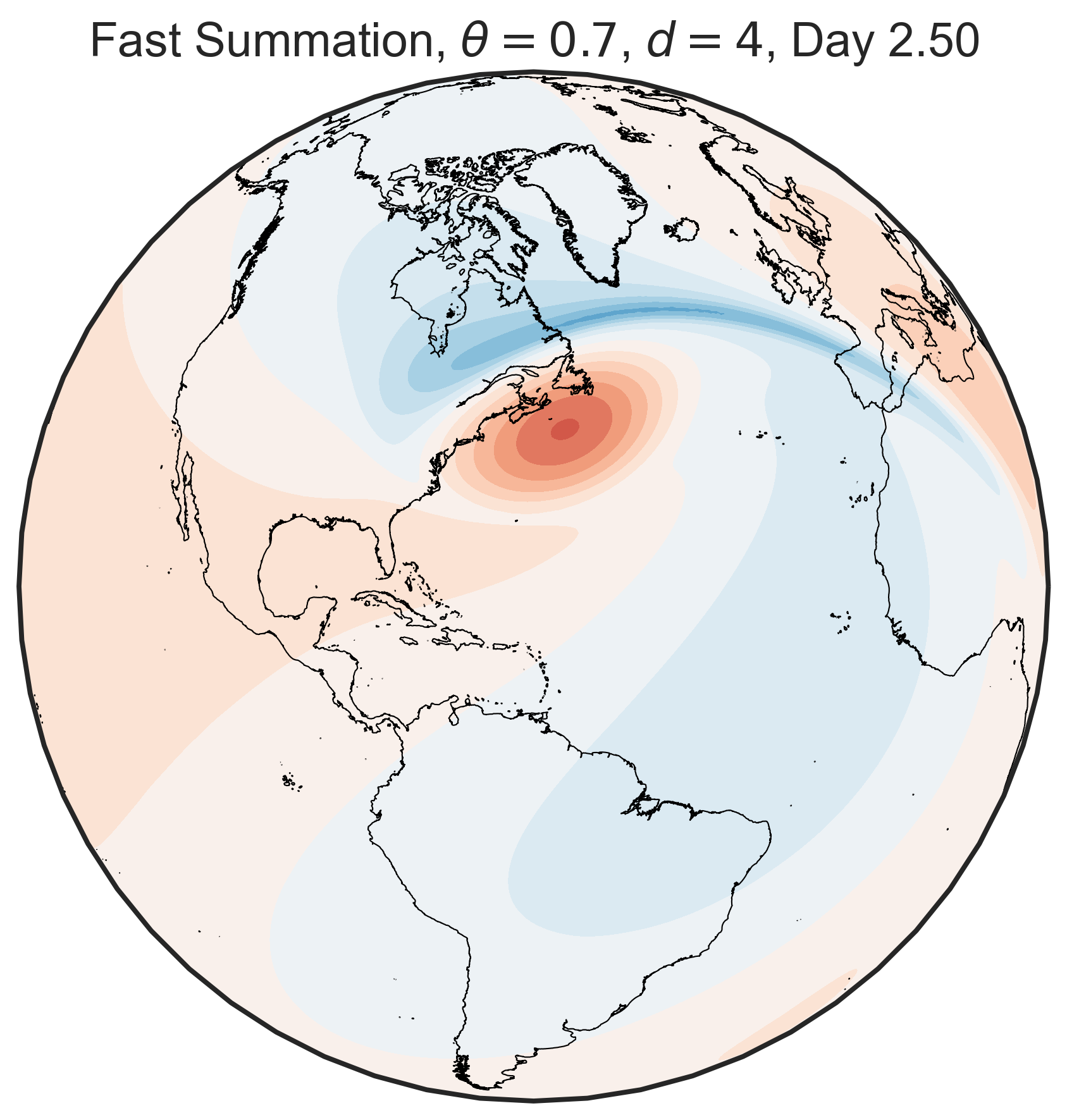}
    \end{subfigure}
    \hspace{0.025\textwidth}
    \begin{subfigure}{0.20\textwidth}
        \includegraphics[width=\linewidth]{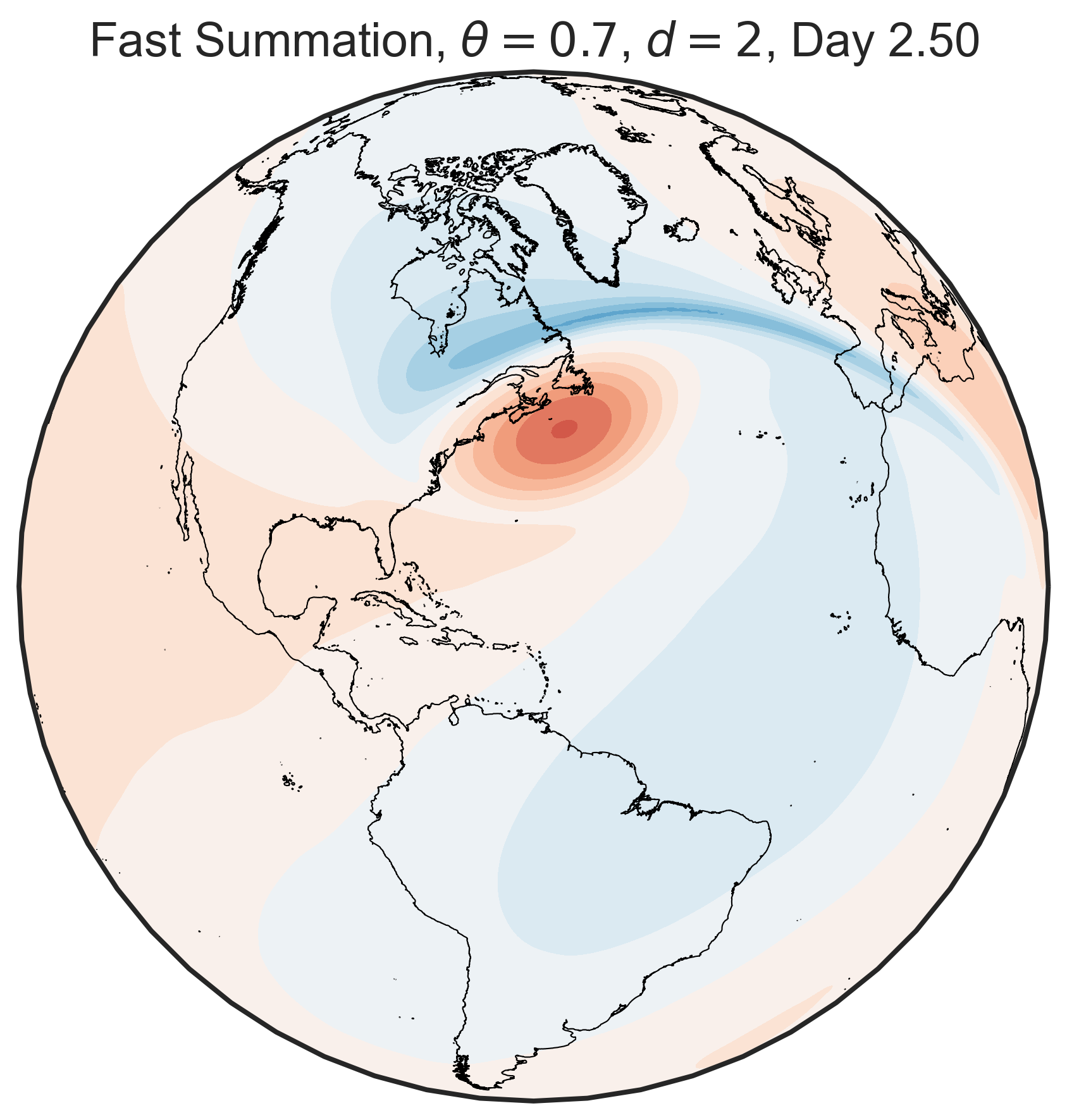}
    \end{subfigure}
    \begin{subfigure}{0.20\textwidth}
        \includegraphics[width=\linewidth]{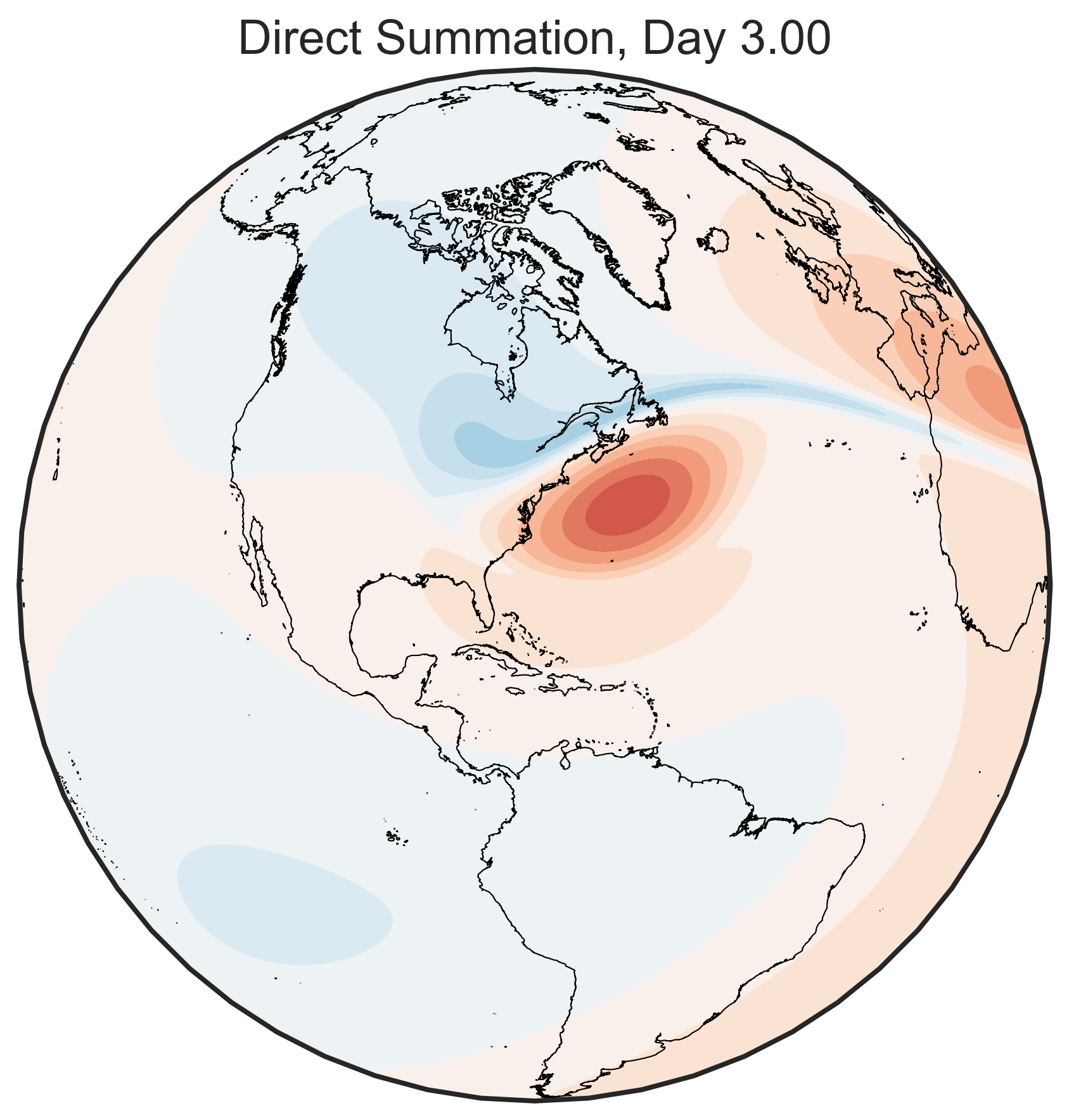}
    \end{subfigure}
    \hspace{0.025\textwidth}
    \begin{subfigure}{0.20\textwidth}
        \includegraphics[width=\linewidth]{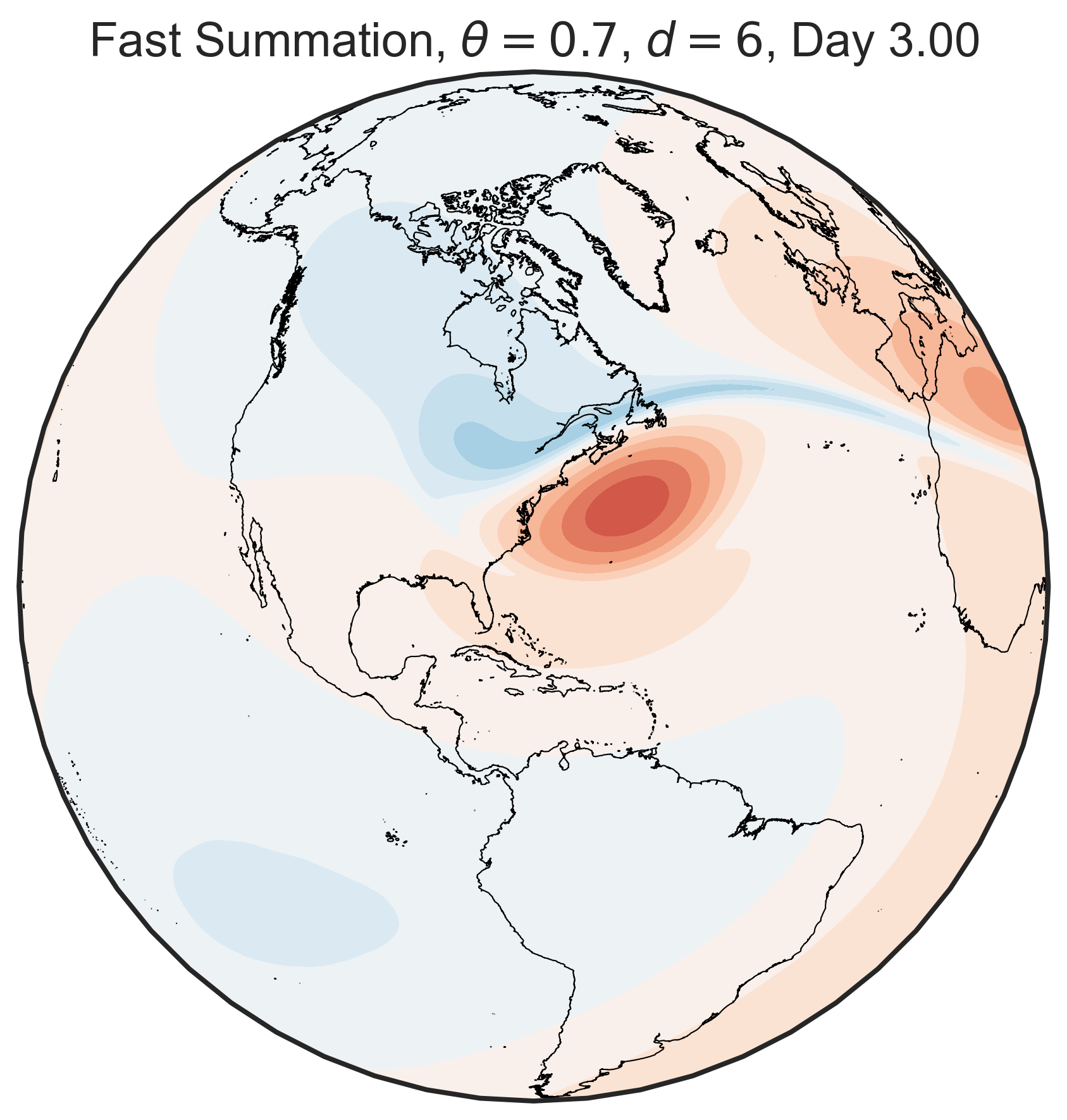}
    \end{subfigure}
    \hspace{0.025\textwidth}
    \begin{subfigure}{0.20\textwidth}
        \includegraphics[width=\linewidth]{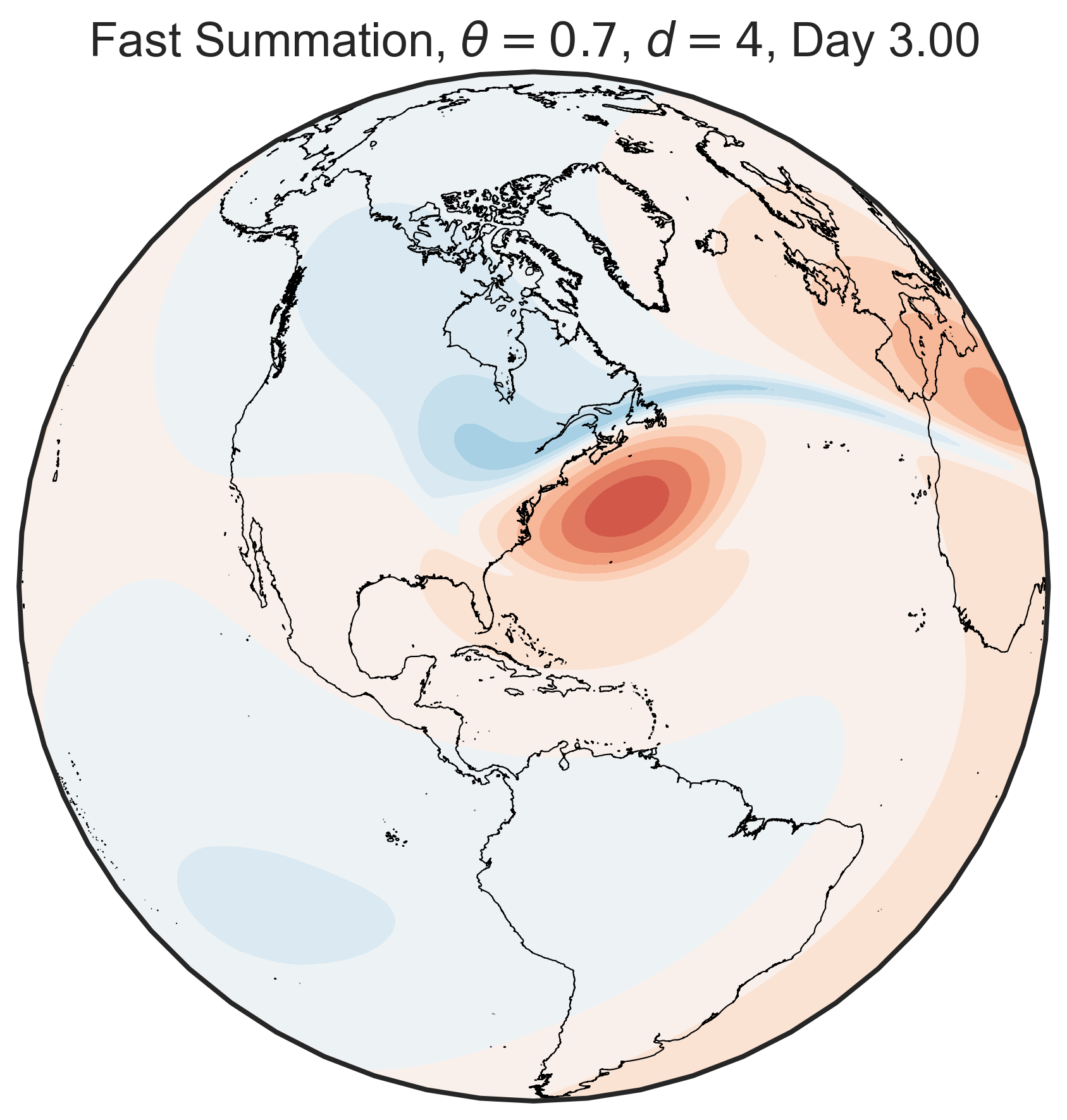}
    \end{subfigure}
    \hspace{0.025\textwidth}
    \begin{subfigure}{0.20\textwidth}
        \includegraphics[width=\linewidth]{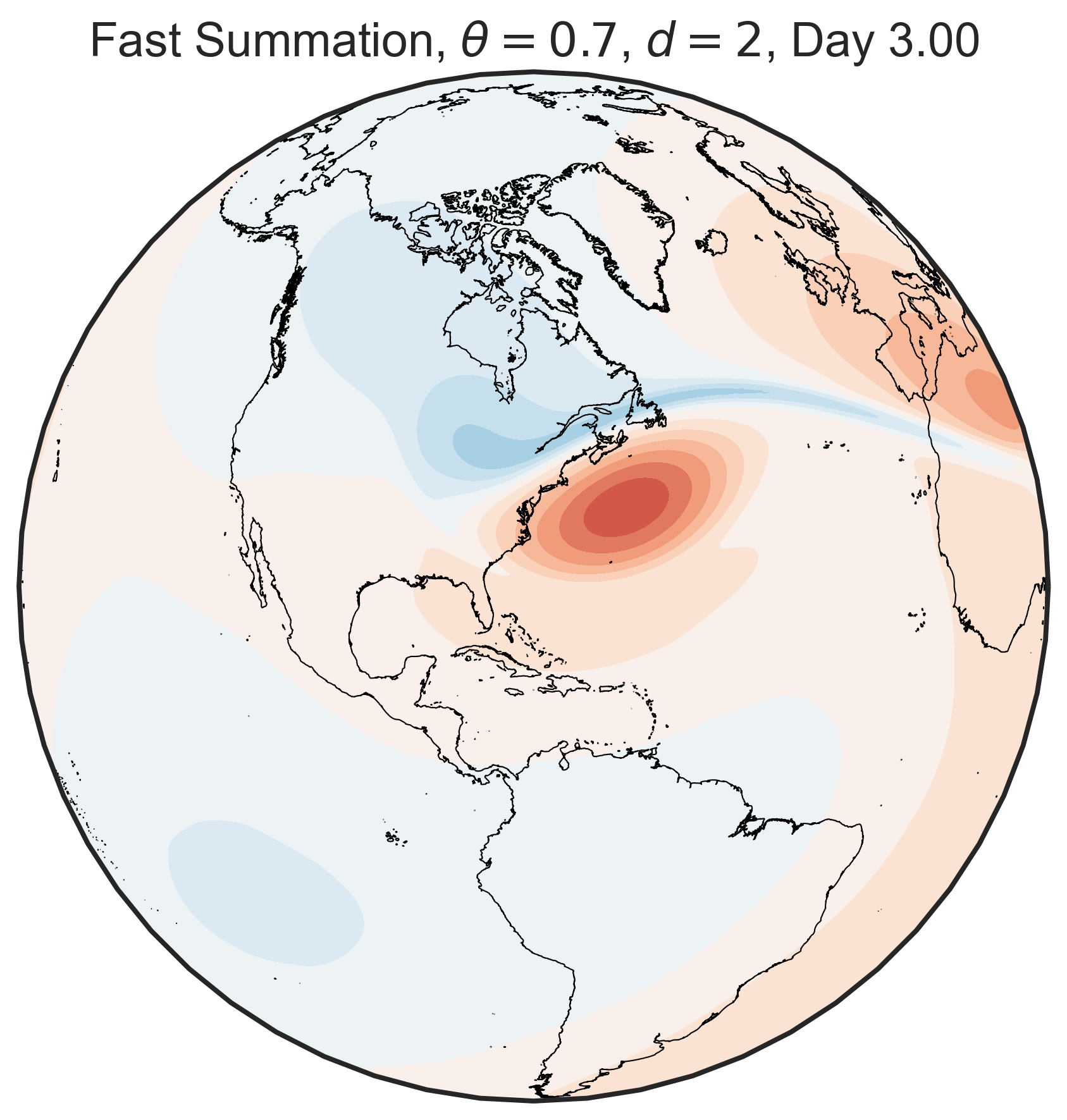}
    \end{subfigure}
    \begin{subfigure}{\textwidth}
        \centering
        \includegraphics[scale=0.25,trim = {0 1.75cm 0 0.31cm}]{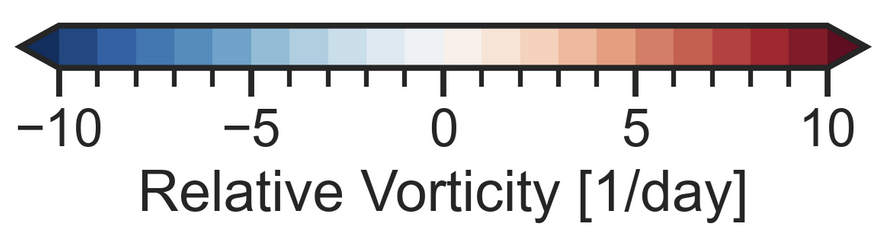}
    \end{subfigure}
    \caption{The results of the Gaussian vortex test case. The left column is computed with direct summation. The remaining three columns are computed with fast summation, with interpolation degree $d=6$ in the second column, $d=4$ in the third column, and $d=2$ in the fourth column. The top row is at $t=0.5$ days, and each row is another $0.5$ days later. We use a time step of $0.01$ days. }
    \label{fig:gvresults}
\end{figure}

We present the results in Fig.~\ref{fig:gvresults}. The Gaussian vortex behaves as expected, moving to the northwest while also developing a counter vortex and a long vortex filament. We observe good agreement in the results between the computation with direct summation and the computations with fast summation when done with $d=4$ and $d=6$. However, when $d=2$, we observe some minute differences. Thus, we can see that the fast summation does not introduce any error that is visible to the eye after 3 days. In these runs, we start with $10242$ particles and we choose AMR parameters $\varepsilon_1=0.0025$ and $\varepsilon_2=0.2$. These runs end with around $36000$ particles. 

\subsection{Polar Vortex Collapse}

\begin{figure}
    \centering
    \begin{subfigure}{0.45\textwidth}
        \includegraphics[width=\linewidth]{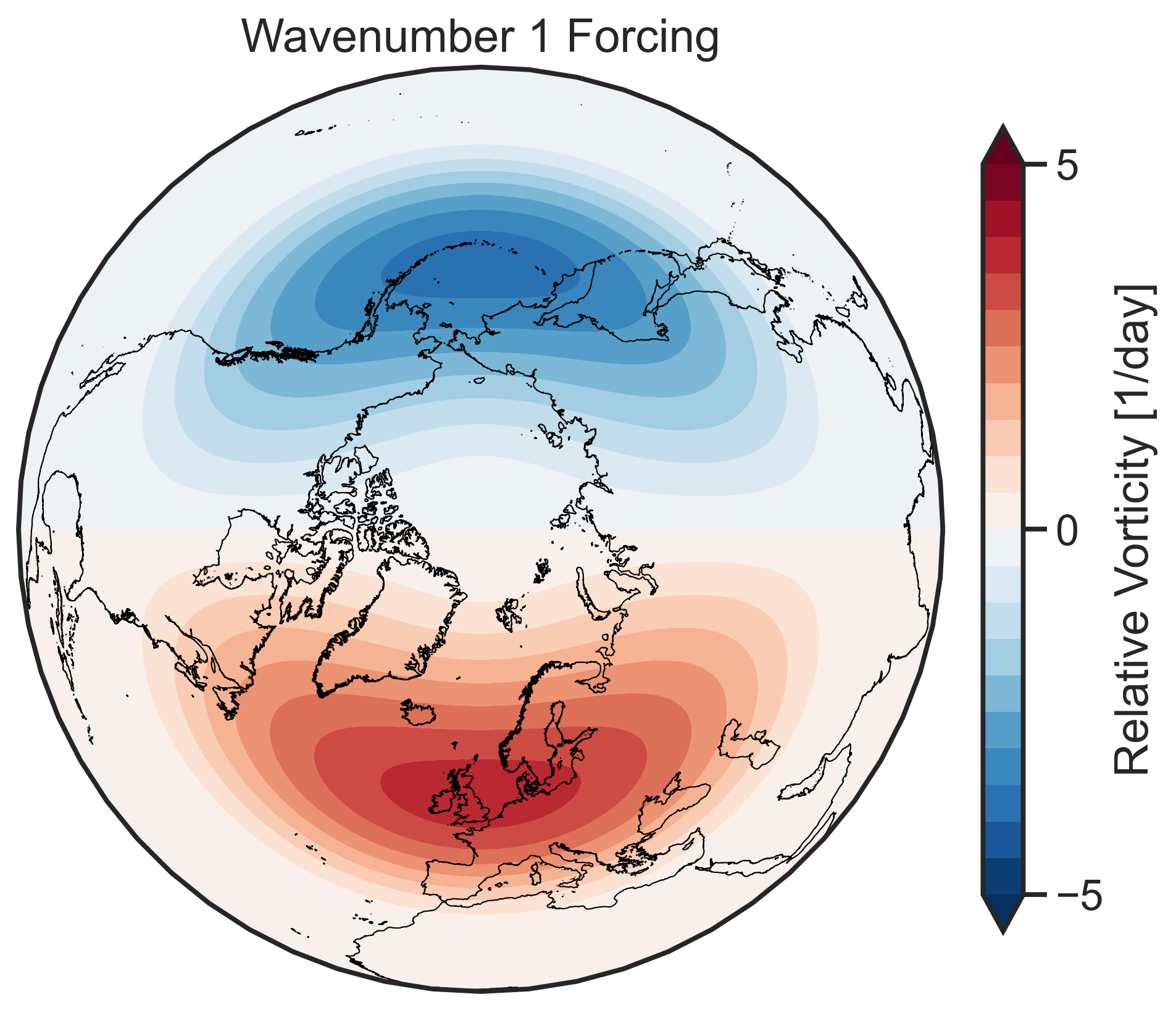}
    \end{subfigure}
    \begin{subfigure}{0.45\textwidth}
        \includegraphics[width=\linewidth]{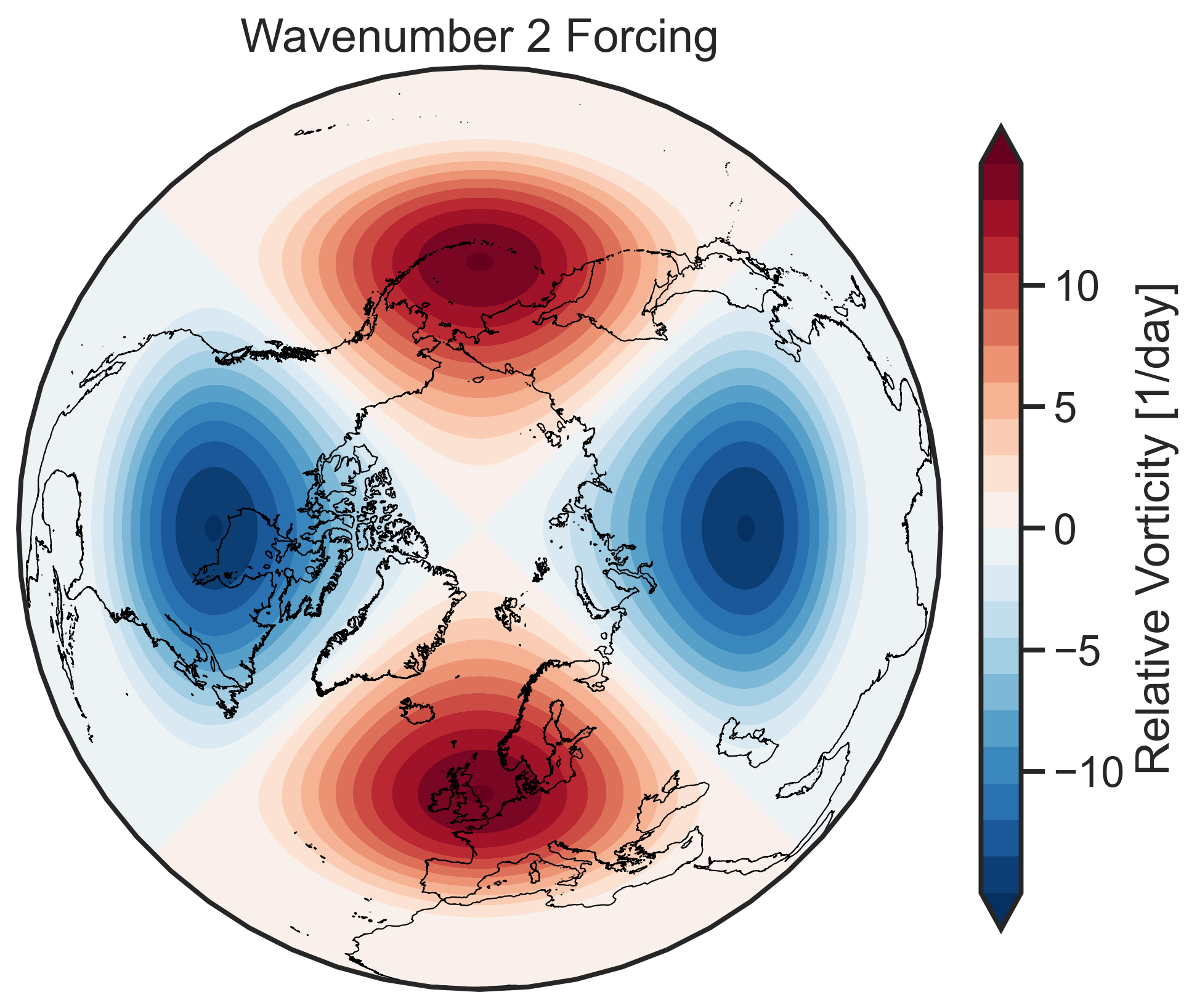}
    \end{subfigure}
    \caption{On the left, a visualization of the wavenumber 1 forcing, and on the right, a visualization of the wavenumber 2 forcing. }
    \label{fig:sswforcing}
\end{figure}

Sudden Stratospheric Warmings (SSW) occur when the polar vortex in the stratosphere is displaced or split, leading to dramatic rises in temperatures in the stratosphere with significant tropospheric impacts~\cite{baldwin2021sudden}. The BVE cannot capture the wave forcing that drives such collapses of the polar vortex, but we can mimic the causes of a collapse with a vorticity forcing, following the work of~\cite{juckes1987high}. We start with an initial vorticity distribution given as 
\begin{equation}
    \zeta_0(\theta,\phi)=\pi e^{-2\beta^2(1-\cos(\theta_0-\theta))}\qty(2\beta^2\cos\theta\sin(\theta_0-\theta)+\sin\theta)
\end{equation} where we take parameters $\beta=1.5$ and $\theta_0=\frac{15}{32}\pi$. We plot this in Fig.~\ref{fig:pv0}. Then, for the forcing we first define 
\begin{equation}
    A(t)=\begin{cases}\frac{1-\cos(\frac{t}{T_p}\pi)}{2} & t < T_p \\ 1 & T_p \leq t < T_f-T_p \\ \frac{1-\cos(\pi+\frac{t-T_f+T_p}{T_p}\pi)}{2} & T_f-T_p \leq t < T_f \\ 0 & t\geq T_f
    \end{cases}
\end{equation} where $0$ to $T_p$ is the ramp up time of the forcing, $T_p$ to $T_f-T_p$ is the period of maximum forcing, and $T_f-T_p$ to $T_f$ is the period of time where the forcing decreases. We take $T_p=4$ and $T_f=15$ days. Next, we have 
\begin{equation}
    B(\theta)=\begin{cases}
        \frac{\tan^2\theta_1}{\tan^2\theta}e^{1-\frac{\tan^2\theta_1}{\tan^2\theta}}&0<\theta\leq\frac{\pi}{2}\\ 0&-\frac{\pi}{2}\leq\theta\leq0
    \end{cases}
\end{equation} where $\theta_1=\frac{\pi}{3}$ is the latitude of maximum forcing. Then, with $A(t)$ and $B(\theta)$ we define the forcing as 
\begin{equation}
    F_R(\theta,\phi,t)=\frac35\Omega k^2A(t)B(\theta)\cos(k\phi)
\end{equation} where $k$ is the zonal wavenumber of the forcing. The choice to scale the amplitude of the forcing as the wavenumber squared is informed by the similar scaling of the Rossby Haurwitz wave in the wavenumber, as seen in Eq.~\ref{eq:rhwave}. With the forcing, the Biot-Savart integral then becomes 
\begin{equation}
    \mathbf{u}(\mathbf{x},t)=-\frac{1}{4\pi}\int_S\frac{\mathbf{x}\cross\mathbf{y}}{1-\mathbf{x}\cdot\mathbf{y}}(\zeta(\mathbf{y})-F_R(\mathbf{y},t))dS(\mathbf{y})
\end{equation}
The forcing is plotted for $k=1$ and $k=2$ in Fig.~\ref{fig:sswforcing}. The Rossby number of the forcing is $0.3k^2$. We choose these two wavenumber forcings because these have been tied to different types of polar vortex collapses~\cite{baldwin2021sudden}. In particular, wavenumber 1 forcing has been tied to displacements of the polar vortex~\cite{castanheira2010dynamical}, while wavenumber 2 forcing has been tied to splits of the polar vortex~\cite{martius2009blocking}. 

\begin{figure}
    \centering
    \begin{subfigure}{0.20\textwidth}
        \includegraphics[width=\linewidth]{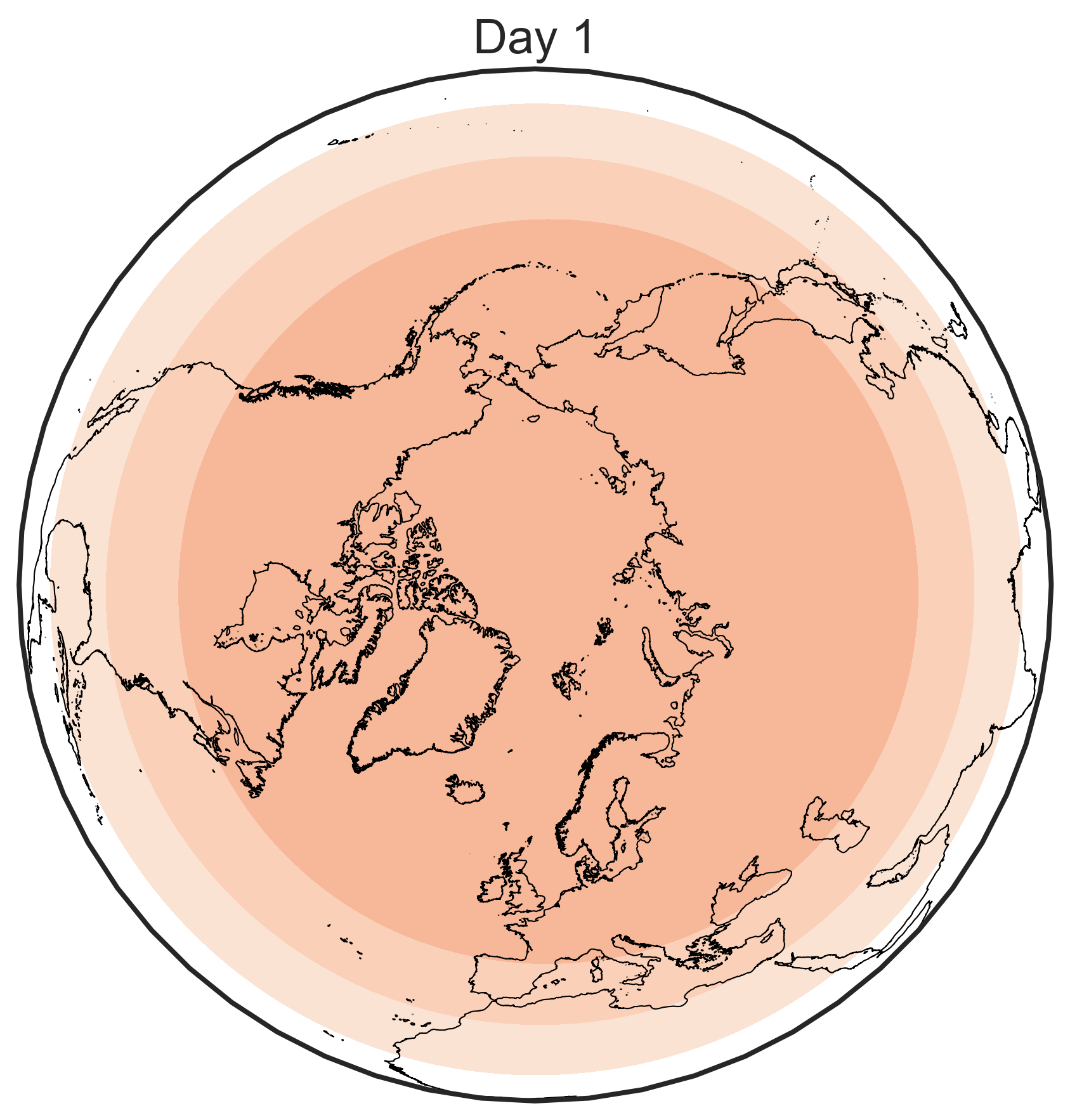}
    \end{subfigure}
    \hspace{0.025\textwidth}
    \begin{subfigure}{0.20\textwidth}
        \includegraphics[width=\linewidth]{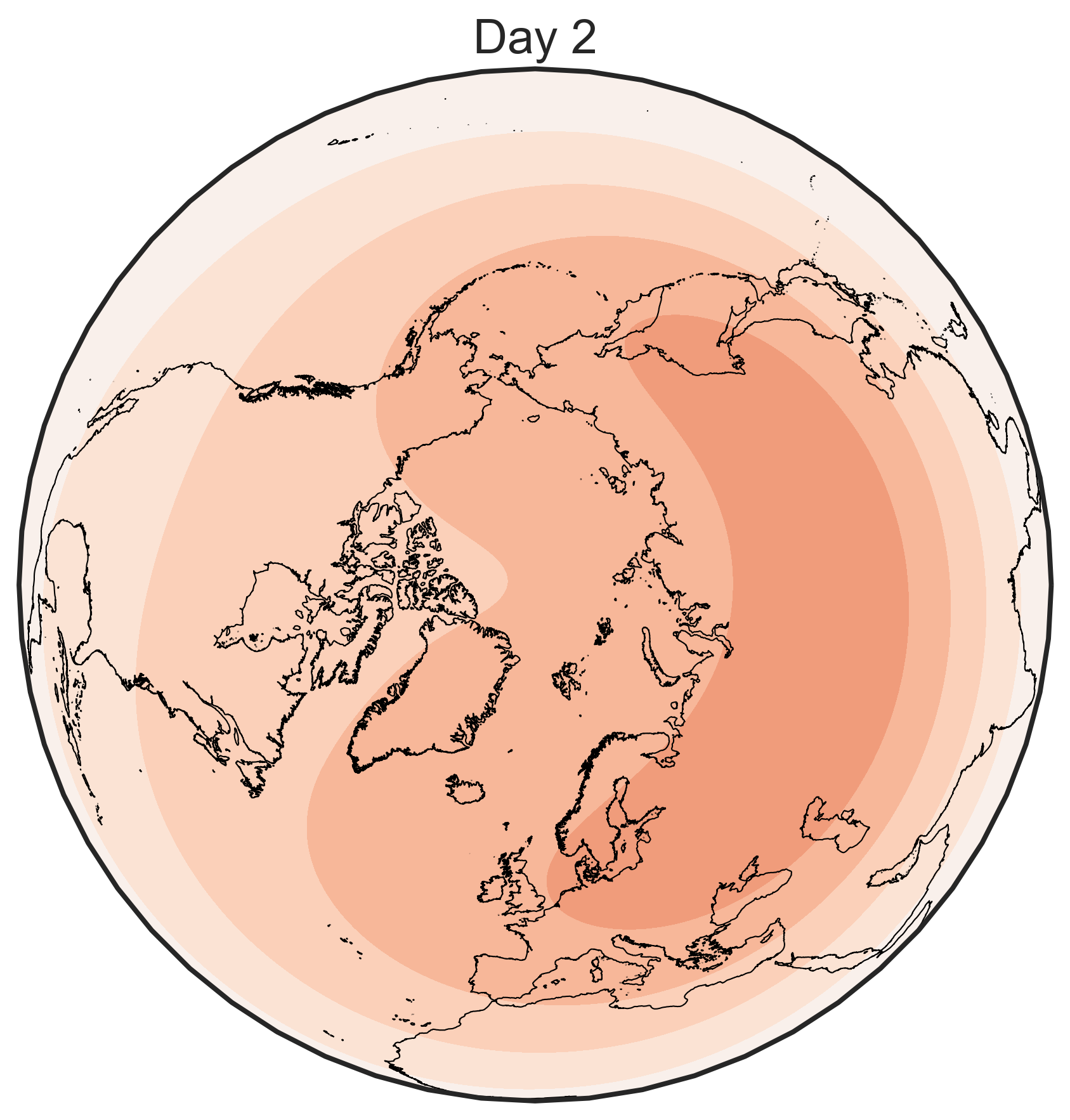}
    \end{subfigure}
    \hspace{0.025\textwidth}
    \begin{subfigure}{0.20\textwidth}
        \includegraphics[width=\linewidth]{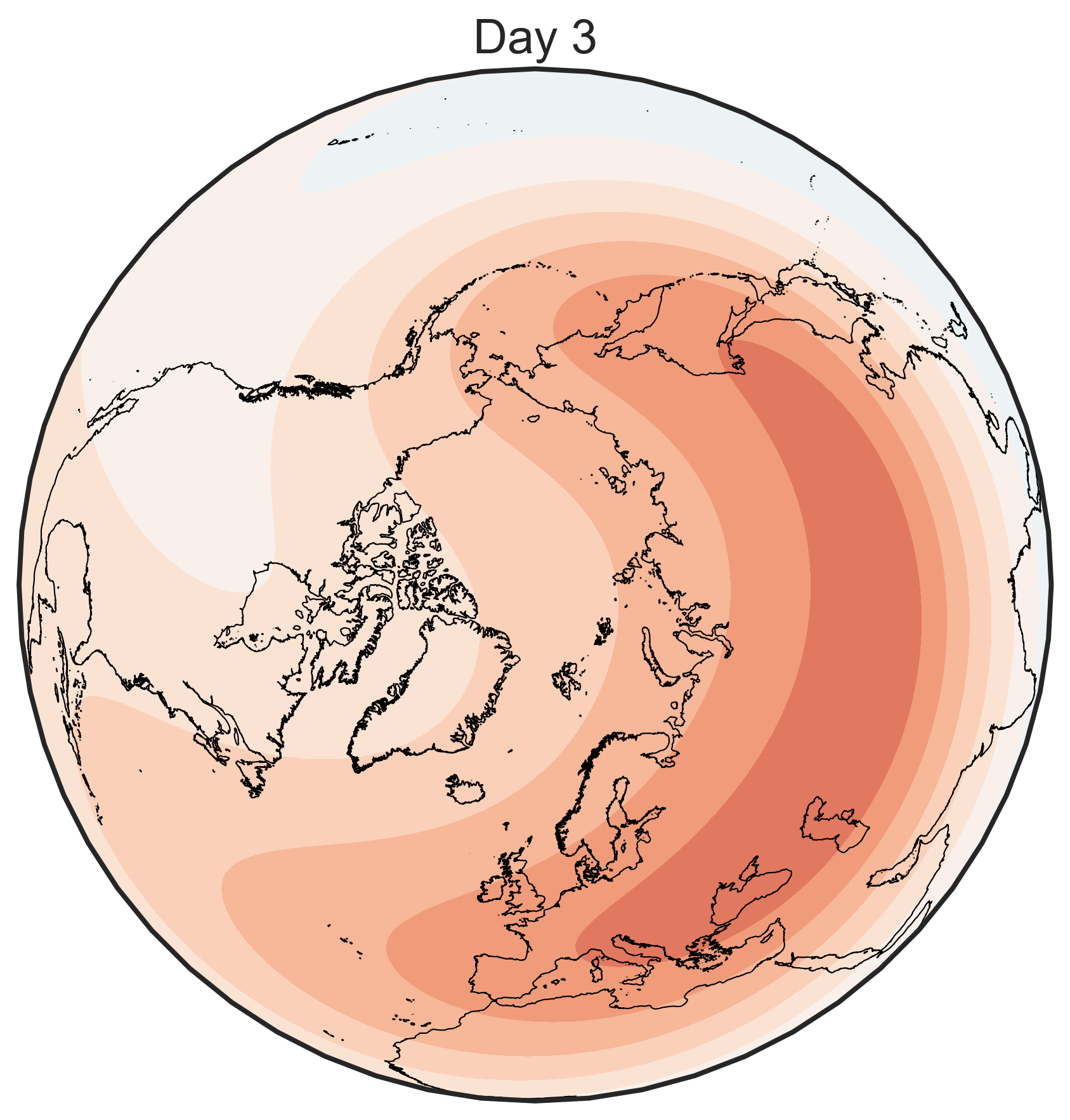}
    \end{subfigure}
    \hspace{0.025\textwidth}
    \begin{subfigure}{0.20\textwidth}
        \includegraphics[width=\linewidth]{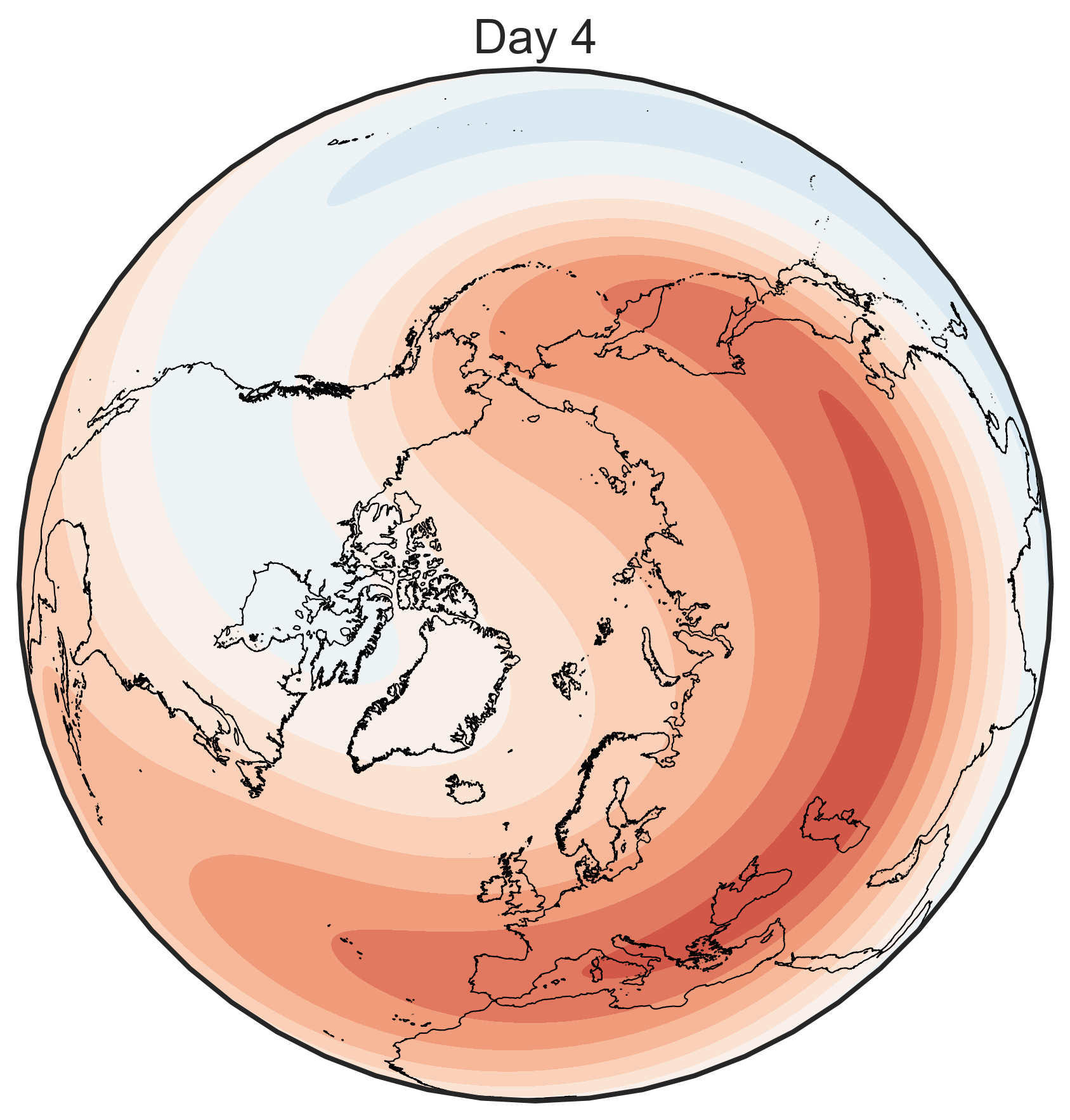}
    \end{subfigure}
    \begin{subfigure}{0.20\textwidth}
        \includegraphics[width=\linewidth]{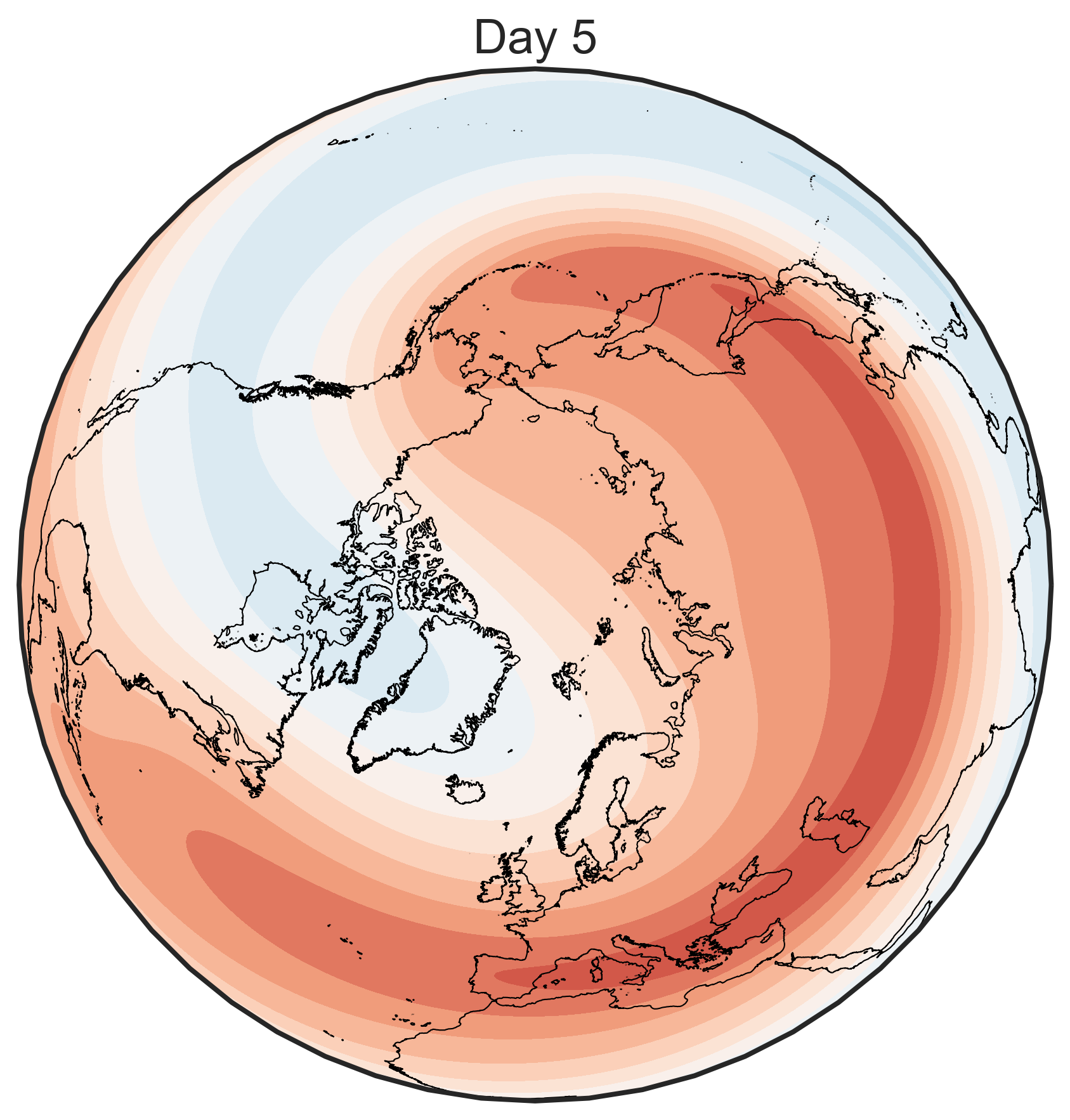}
    \end{subfigure}
    \hspace{0.025\textwidth}
    \begin{subfigure}{0.20\textwidth}
        \includegraphics[width=\linewidth]{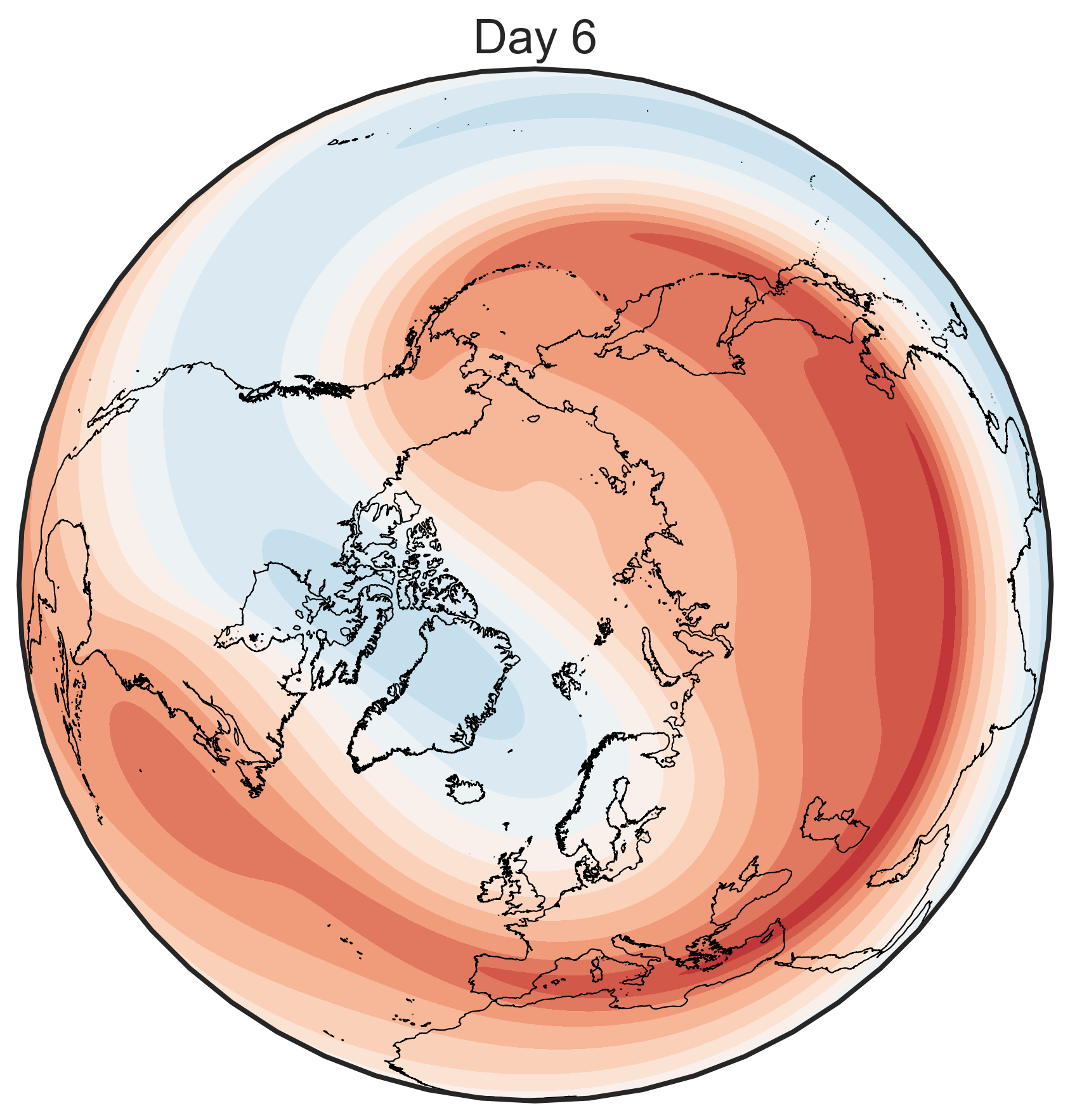}
    \end{subfigure}
    \hspace{0.025\textwidth}
    \begin{subfigure}{0.20\textwidth}
        \includegraphics[width=\linewidth]{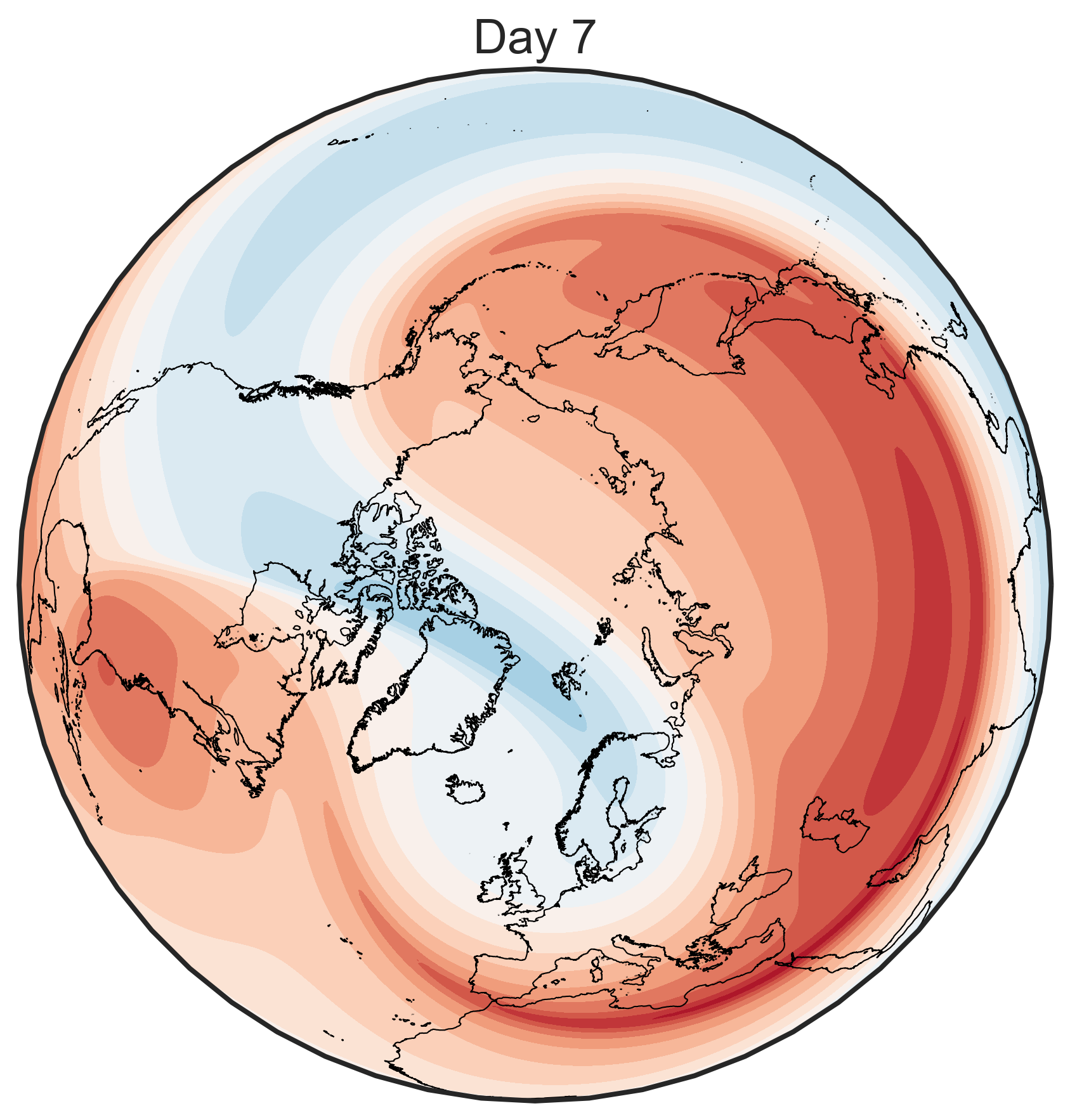}
    \end{subfigure}
    \hspace{0.025\textwidth}
    \begin{subfigure}{0.20\textwidth}
        \includegraphics[width=\linewidth]{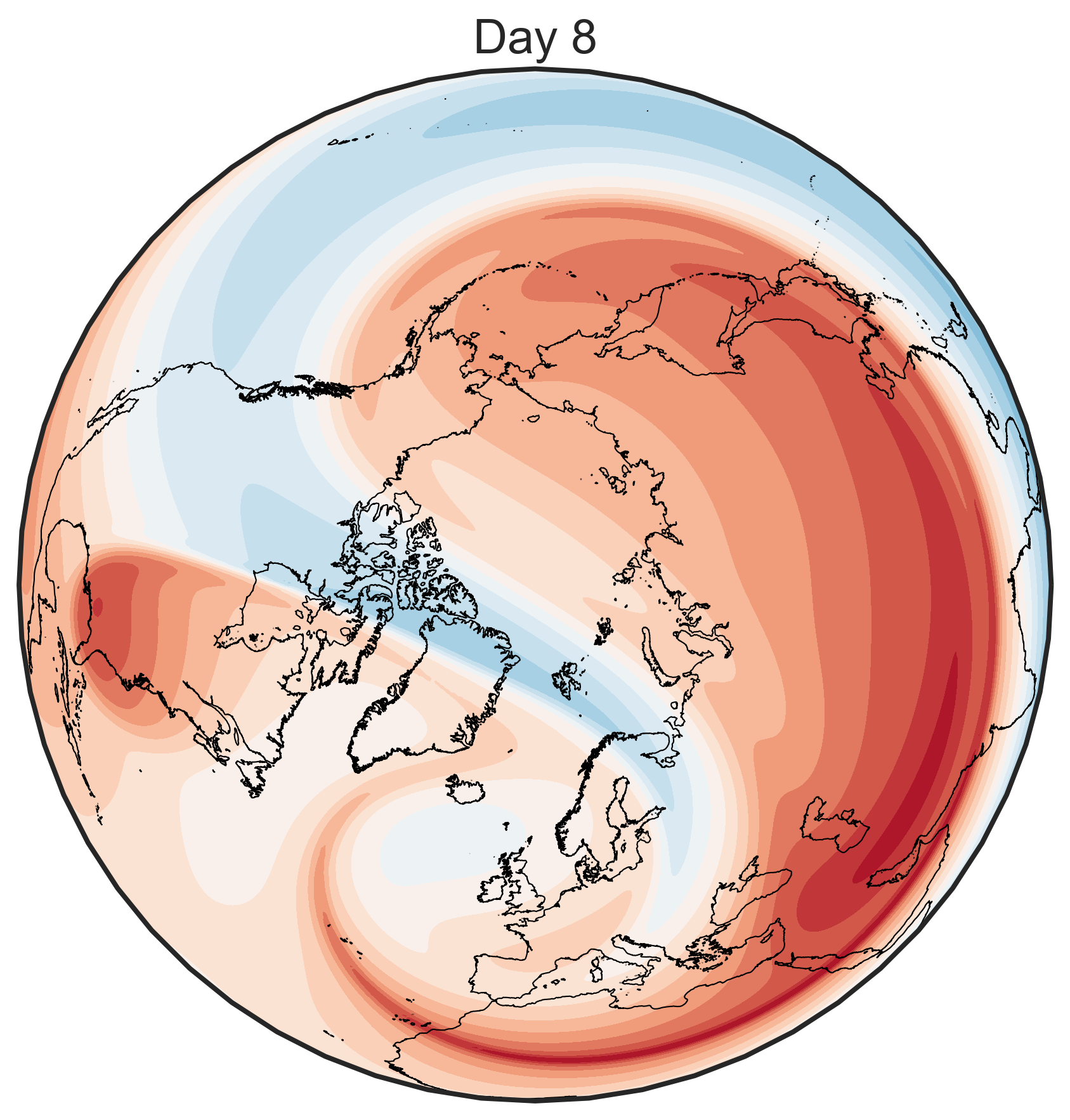}
    \end{subfigure}
    \begin{subfigure}{0.20\textwidth}
        \includegraphics[width=\linewidth]{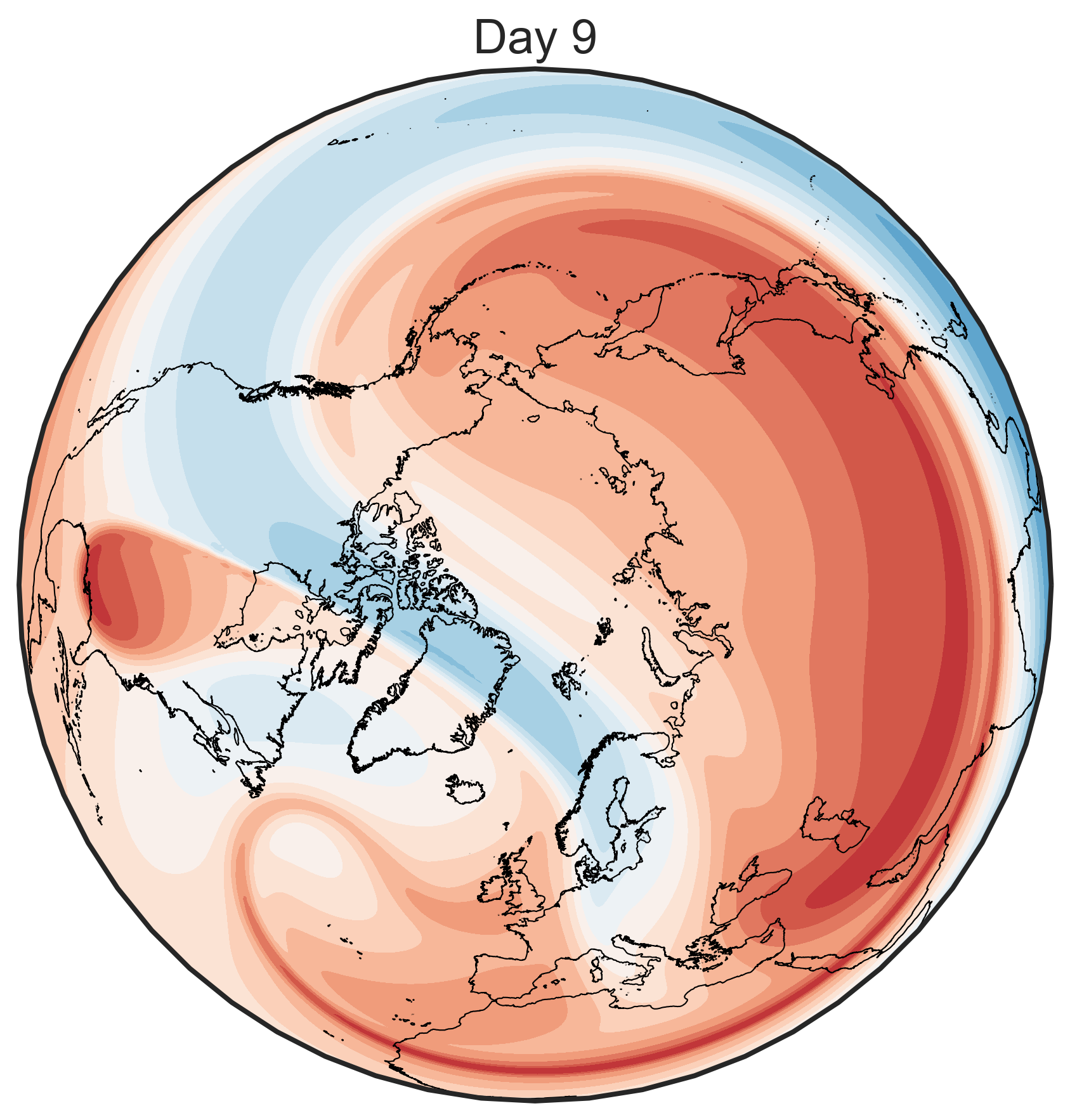}
    \end{subfigure}
    \hspace{0.025\textwidth}
    \begin{subfigure}{0.20\textwidth}
        \includegraphics[width=\linewidth]{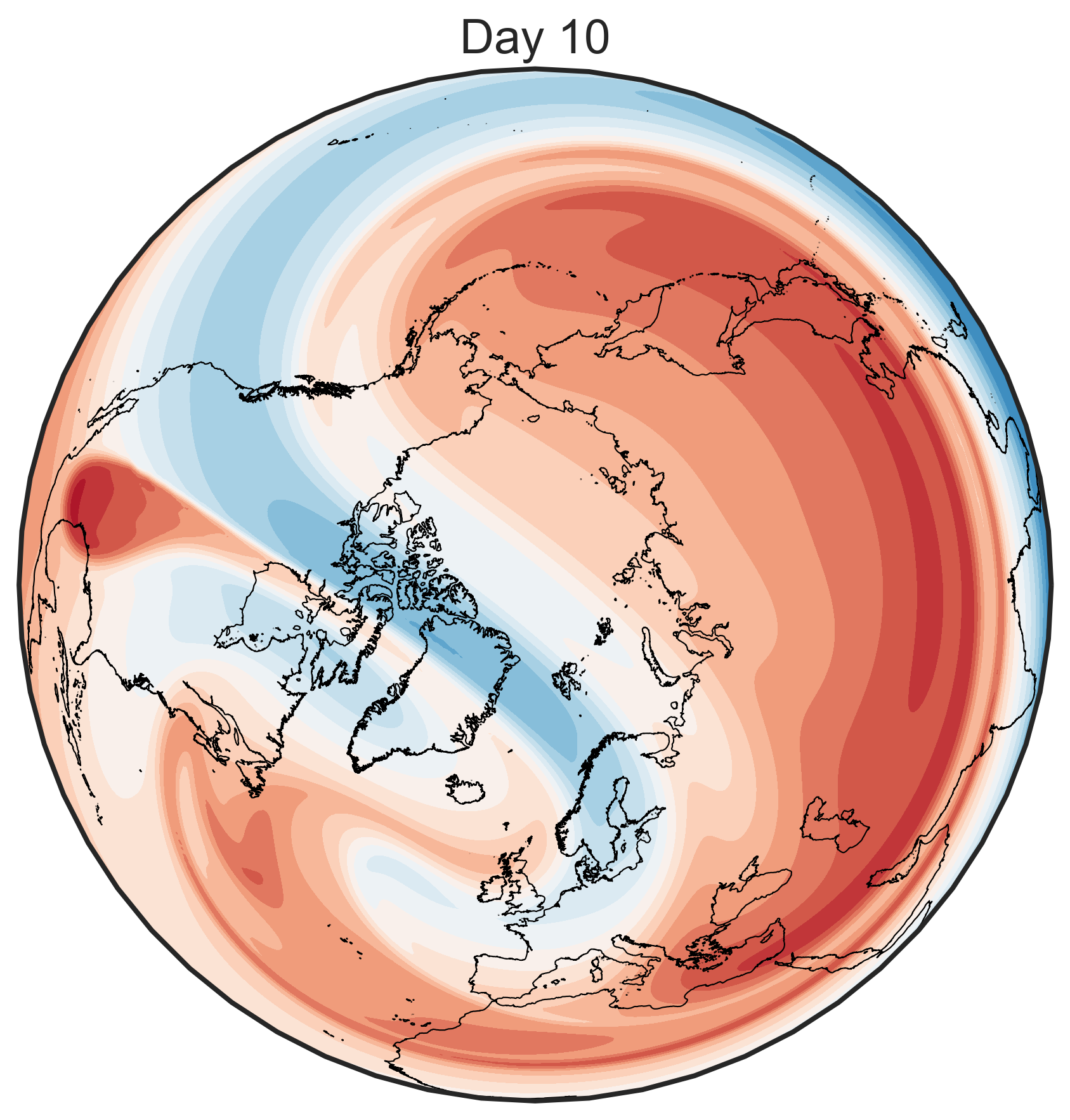}
    \end{subfigure}
    \hspace{0.025\textwidth}
    \begin{subfigure}{0.20\textwidth}
        \includegraphics[width=\linewidth]{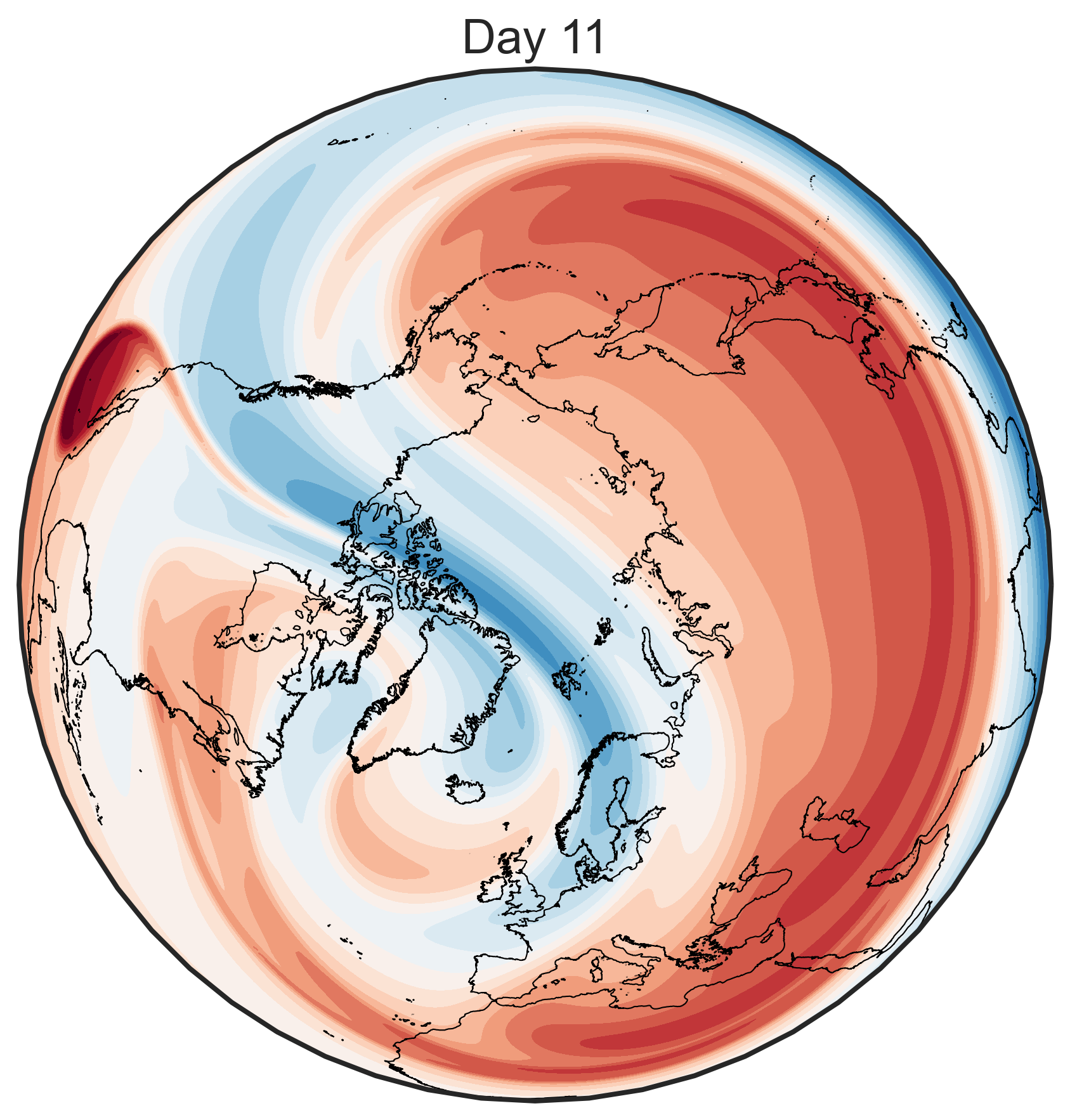}
    \end{subfigure}
    \hspace{0.025\textwidth}
    \begin{subfigure}{0.20\textwidth}
        \includegraphics[width=\linewidth]{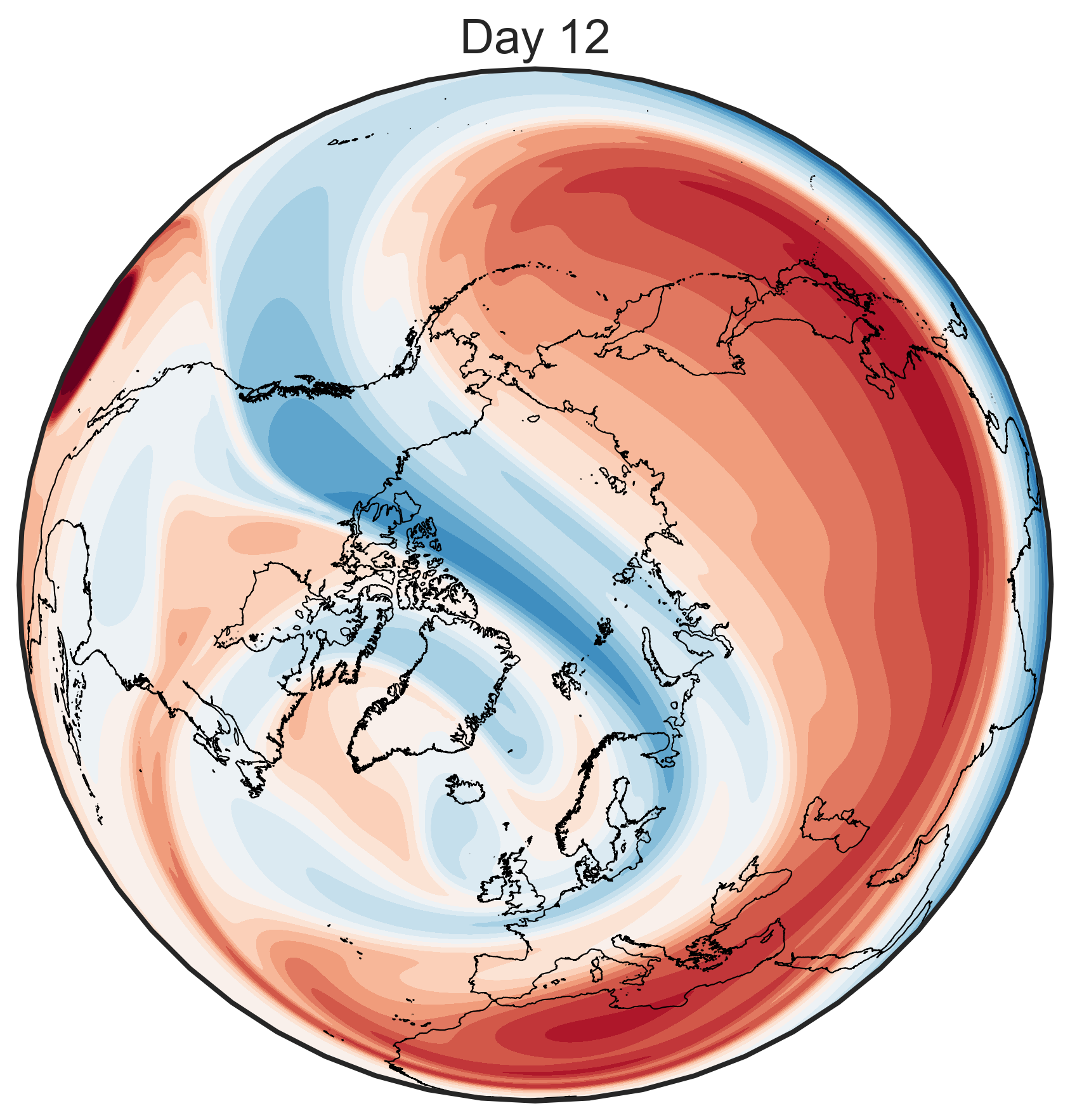}
    \end{subfigure}
    \begin{subfigure}{0.20\textwidth}
        \includegraphics[width=\linewidth]{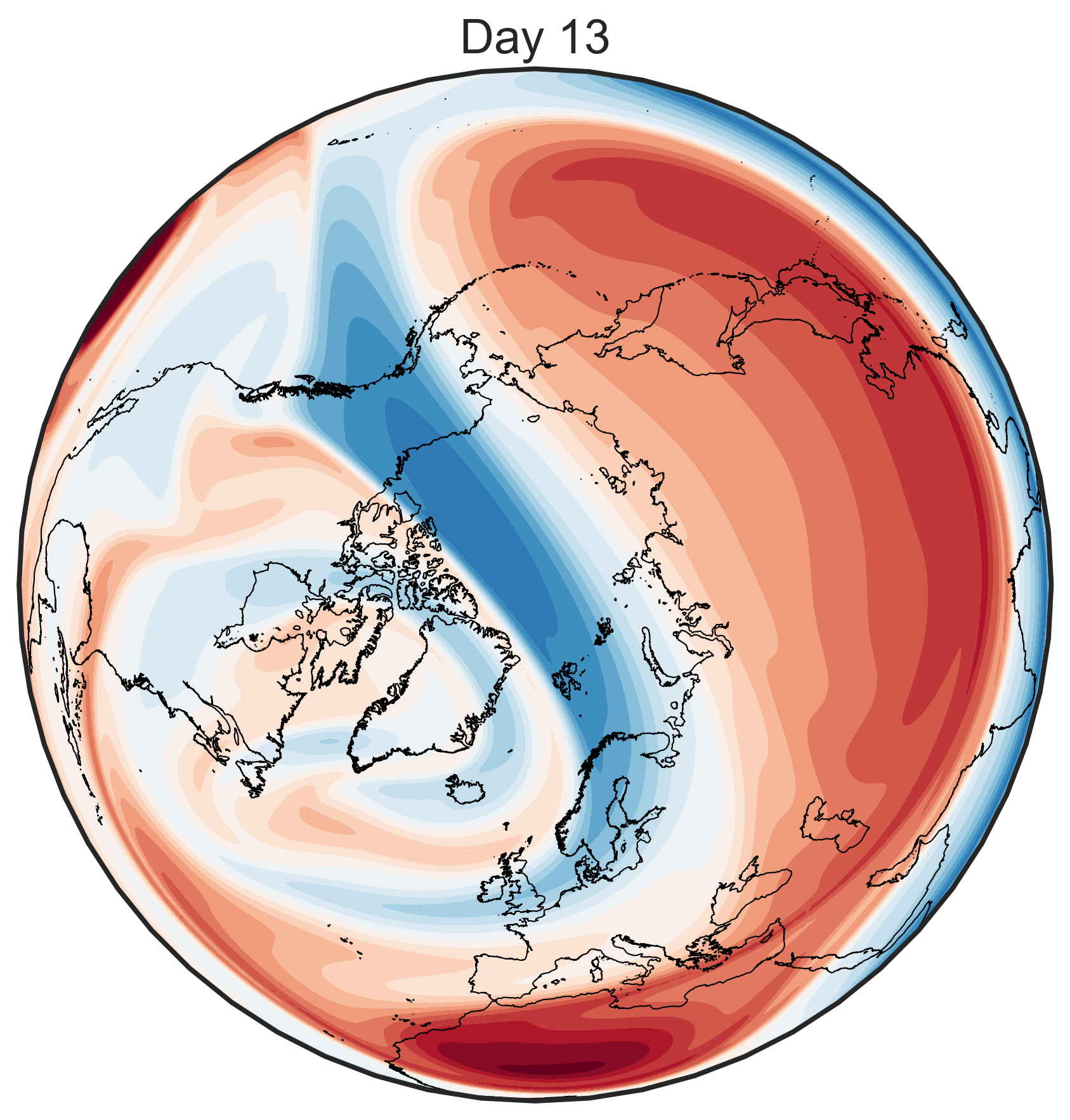}
    \end{subfigure}
    \hspace{0.025\textwidth}
    \begin{subfigure}{0.20\textwidth}
        \includegraphics[width=\linewidth]{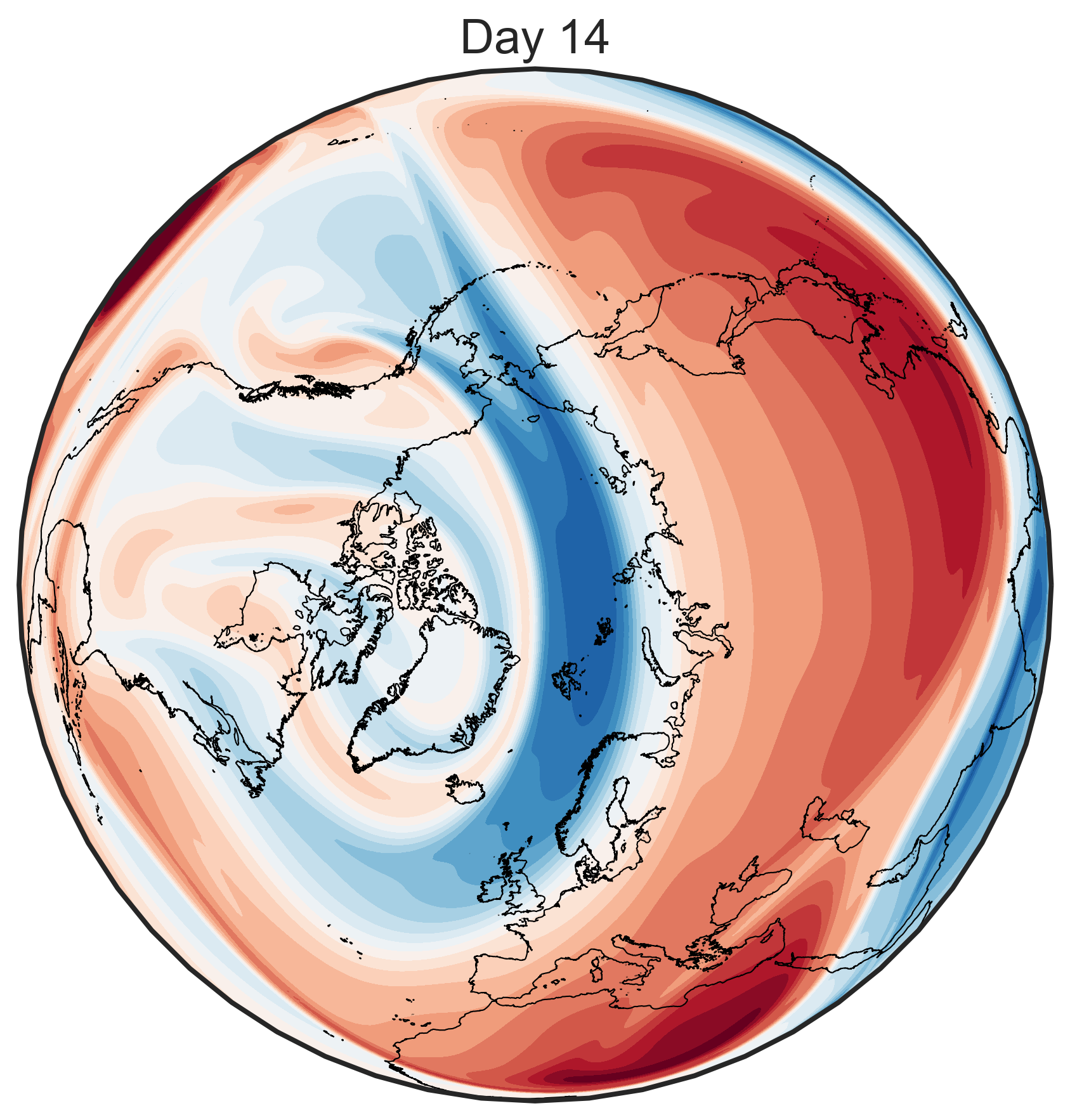}
    \end{subfigure}
    \hspace{0.025\textwidth}
    \begin{subfigure}{0.20\textwidth}
        \includegraphics[width=\linewidth]{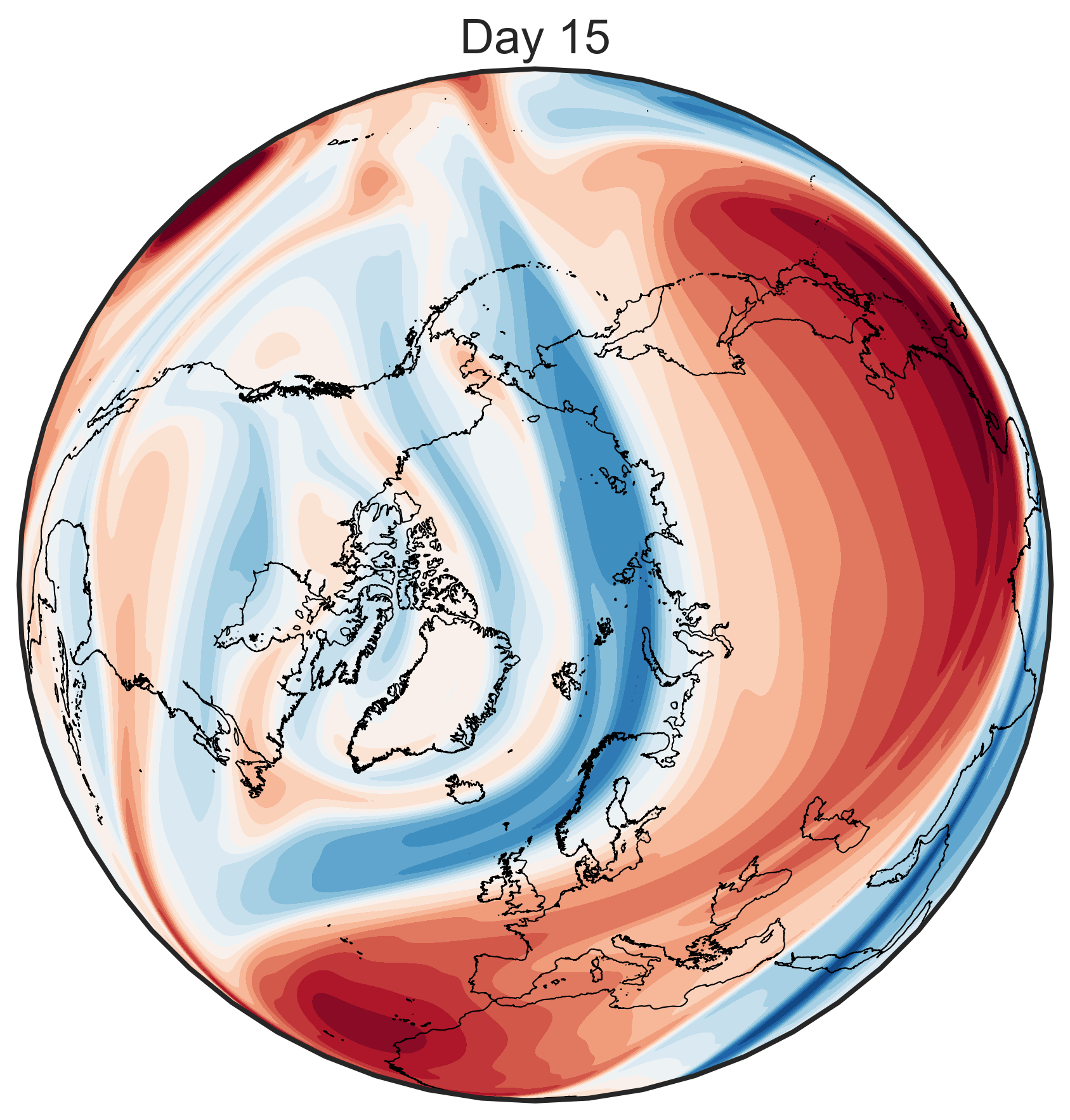}
    \end{subfigure}
    \hspace{0.025\textwidth}
    \begin{subfigure}{0.20\textwidth}
        \includegraphics[width=\linewidth]{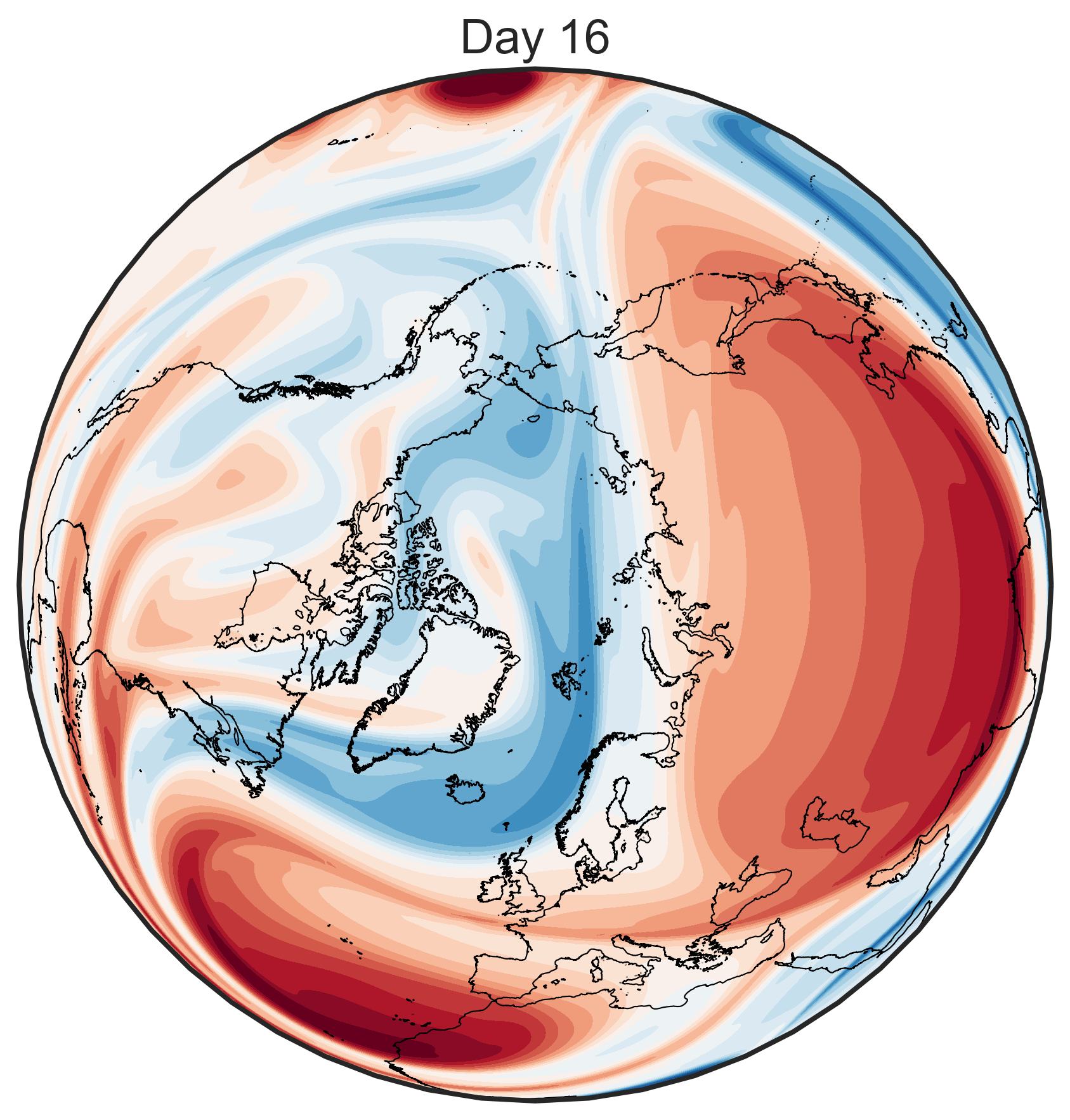}
    \end{subfigure}
    \begin{subfigure}{0.20\textwidth}
        \includegraphics[width=\linewidth]{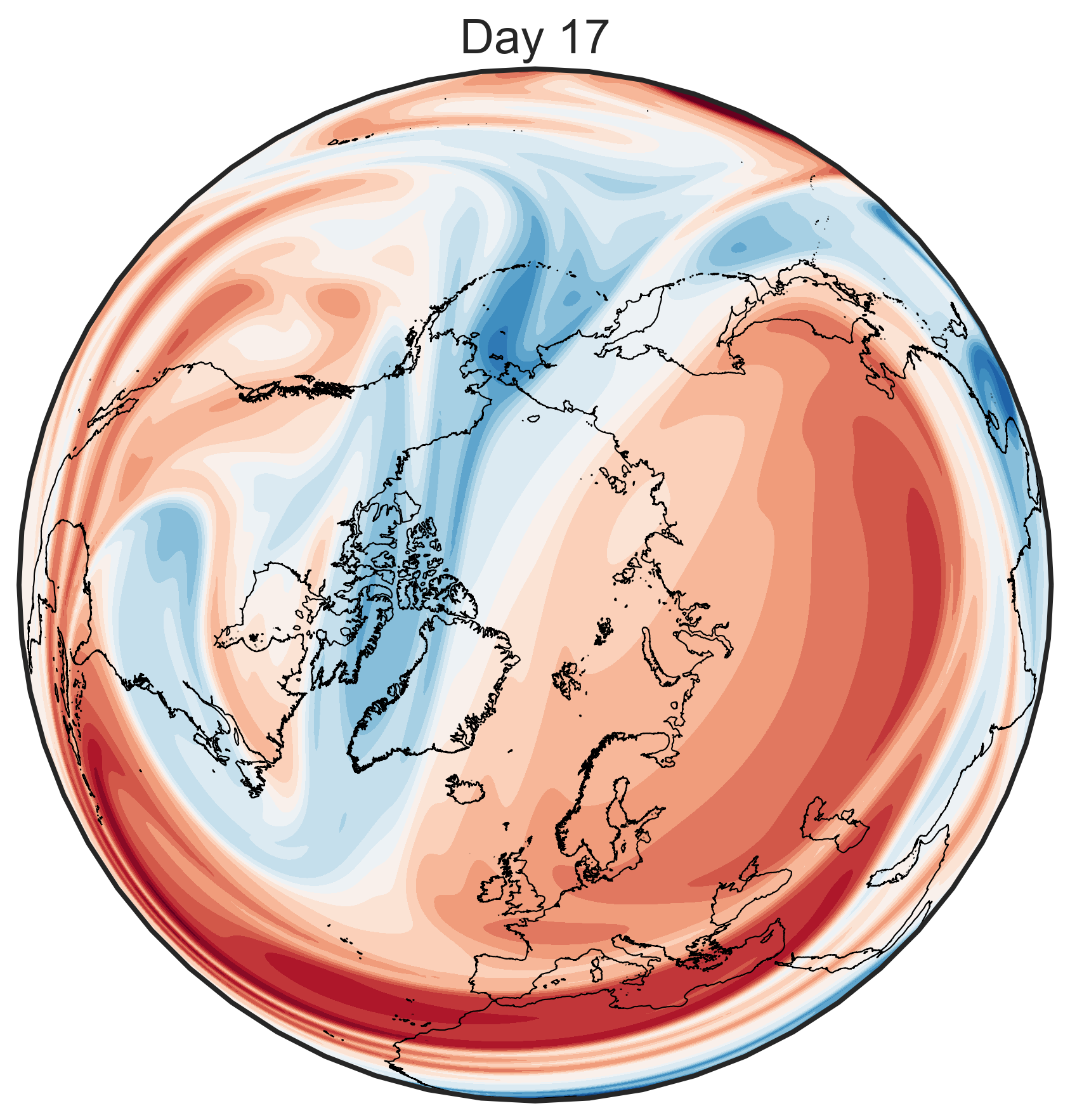}
    \end{subfigure}
    \hspace{0.025\textwidth}
    \begin{subfigure}{0.20\textwidth}
        \includegraphics[width=\linewidth]{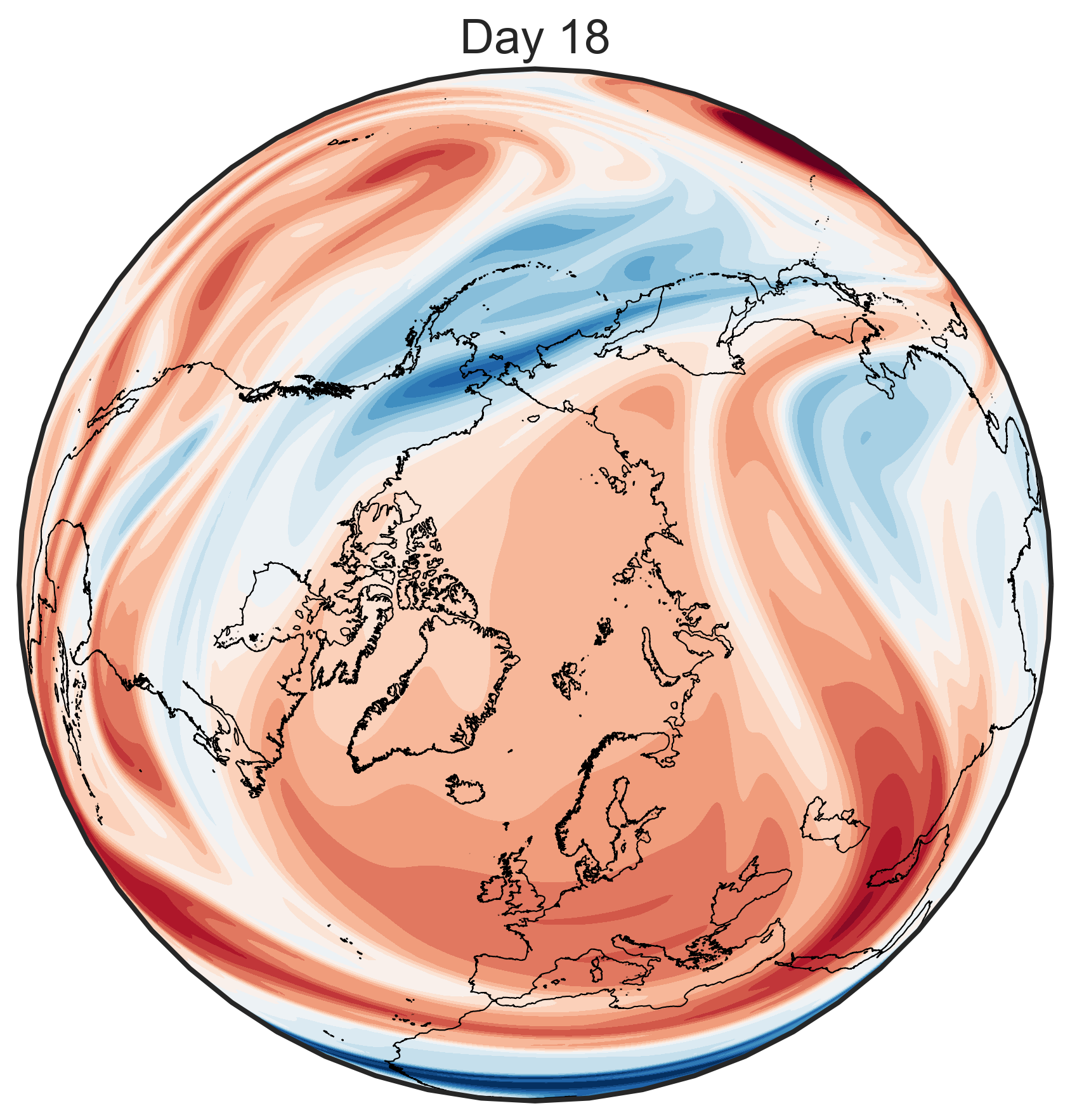}
    \end{subfigure}
    \hspace{0.025\textwidth}
    \begin{subfigure}{0.20\textwidth}
        \includegraphics[width=\linewidth]{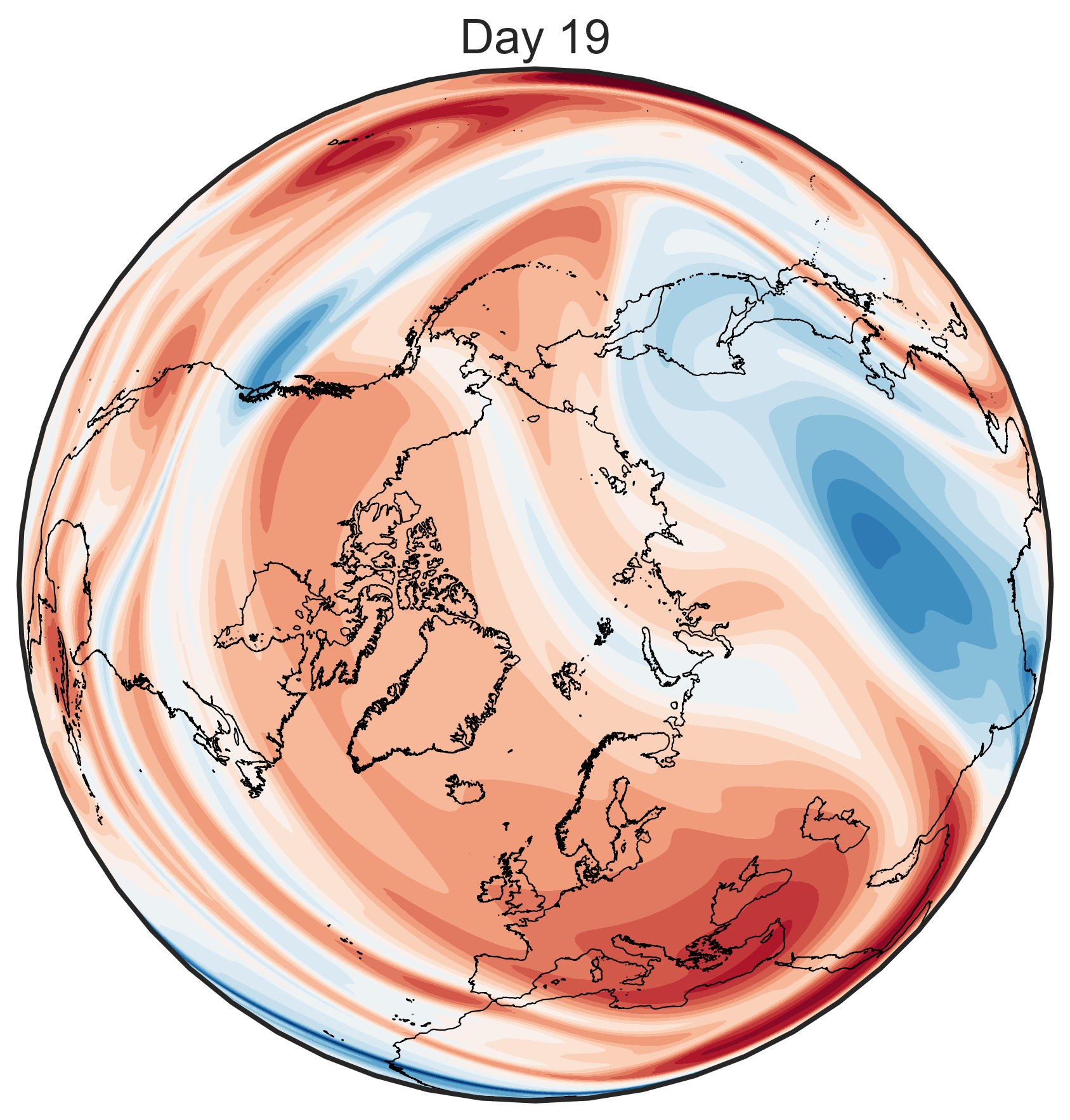}
    \end{subfigure}
    \hspace{0.025\textwidth}
    \begin{subfigure}{0.20\textwidth}
        \includegraphics[width=\linewidth]{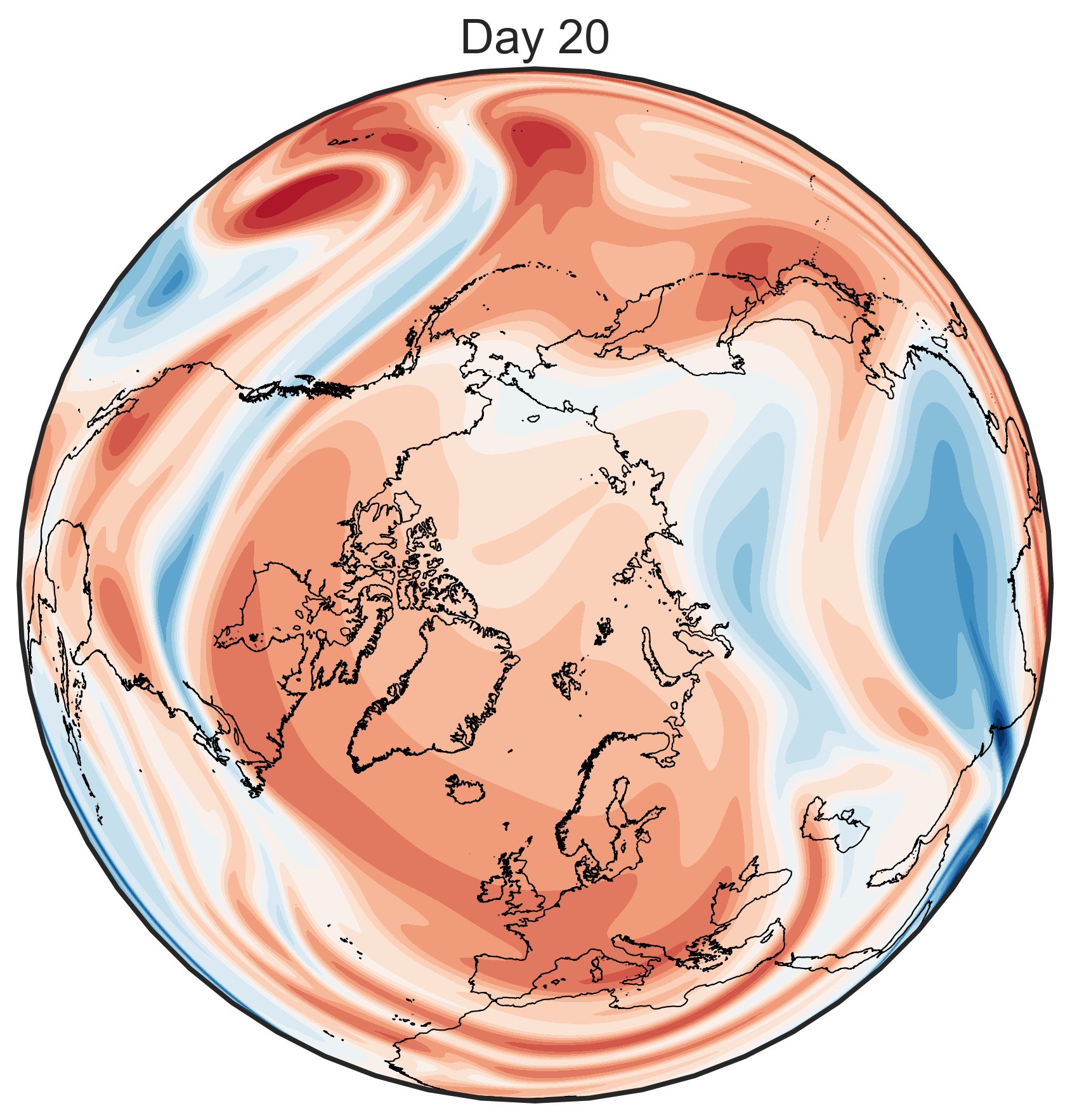}
    \end{subfigure}
    \begin{subfigure}{0.20\textwidth}
        \includegraphics[width=\linewidth]{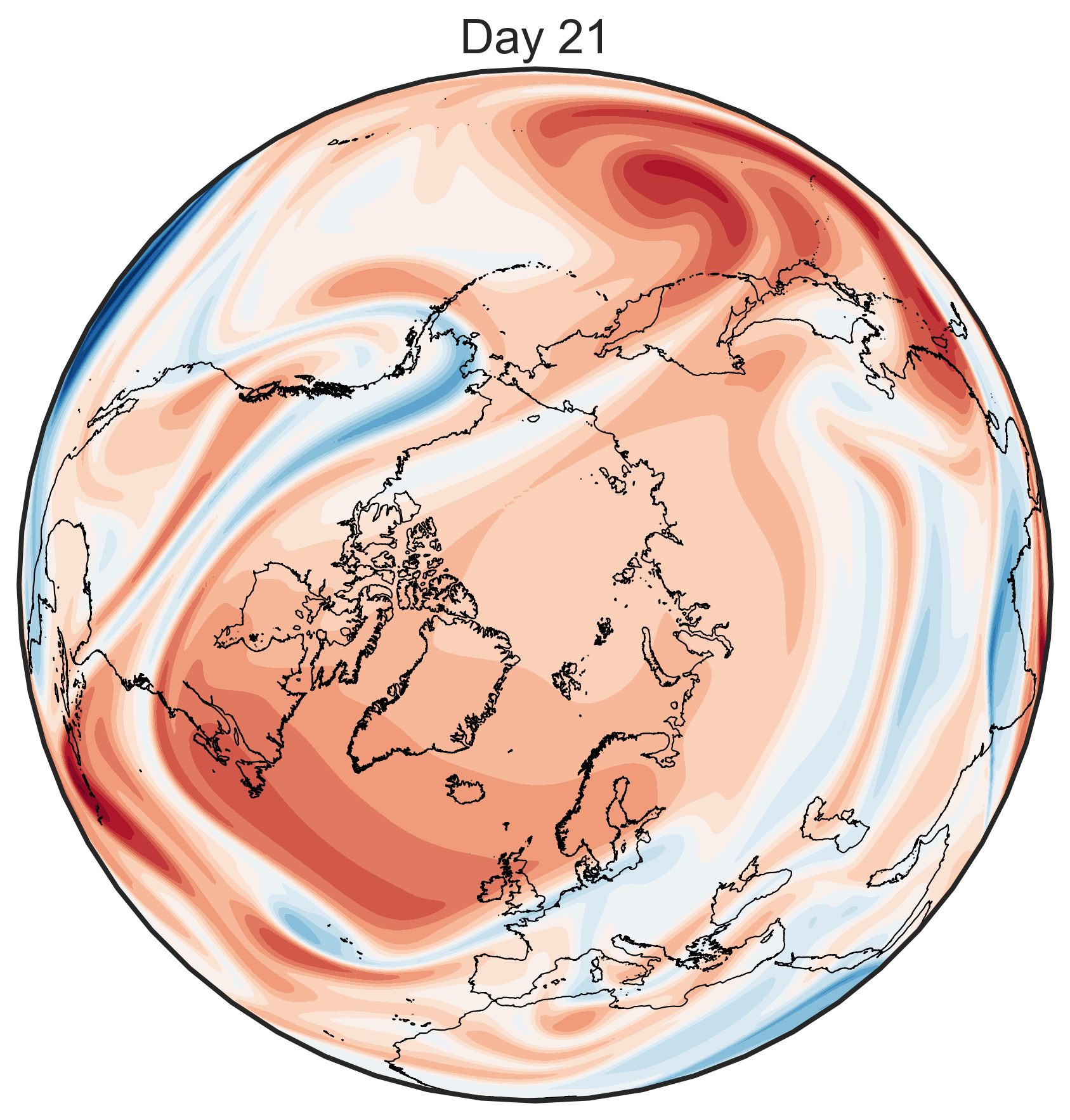}
    \end{subfigure}
    \hspace{0.025\textwidth}
    \begin{subfigure}{0.20\textwidth}
        \includegraphics[width=\linewidth]{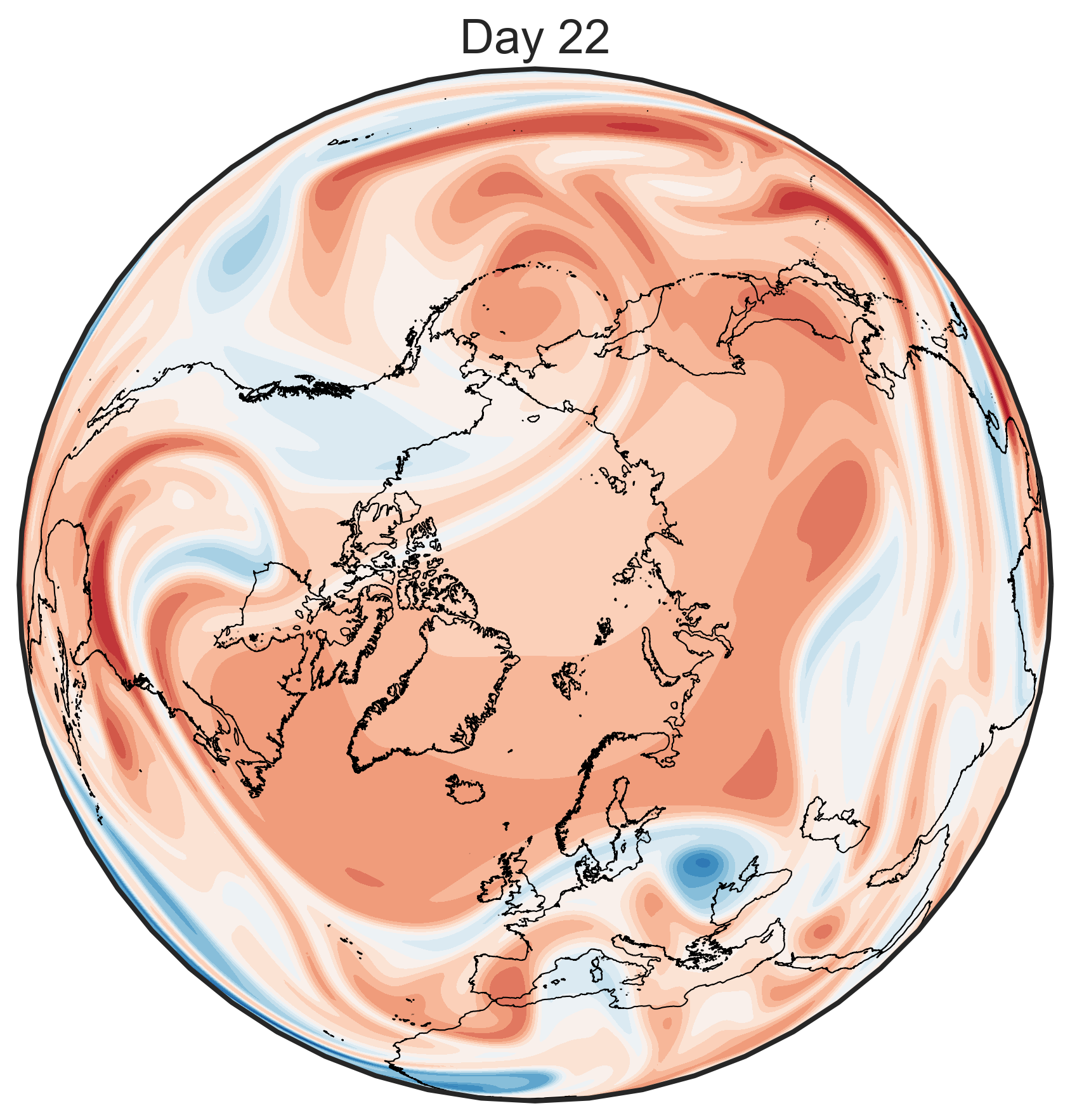}
    \end{subfigure}
    \hspace{0.025\textwidth}
    \begin{subfigure}{0.20\textwidth}
        \includegraphics[width=\linewidth]{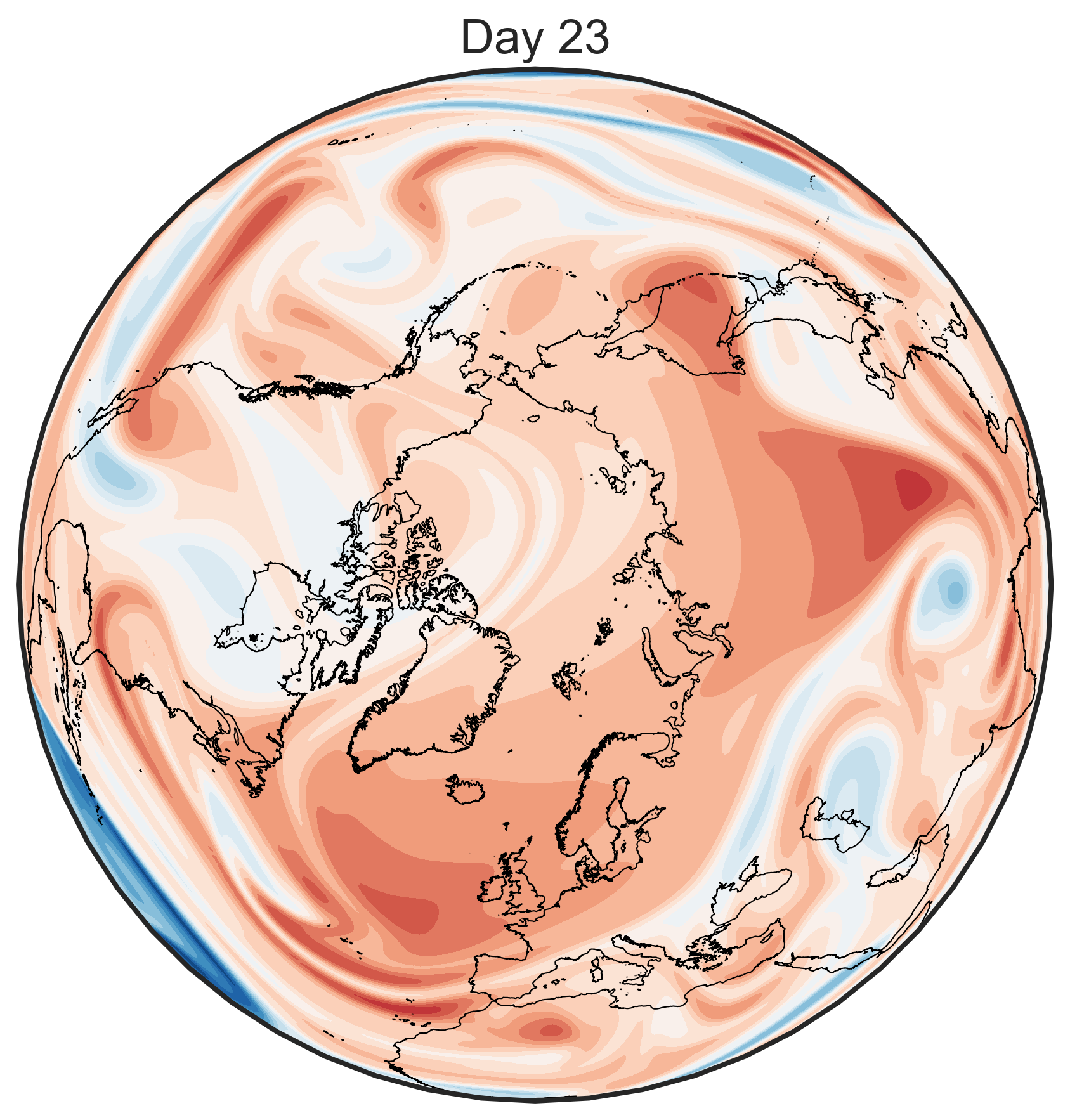}
    \end{subfigure}
    \hspace{0.025\textwidth}
    \begin{subfigure}{0.20\textwidth}
        \includegraphics[width=\linewidth]{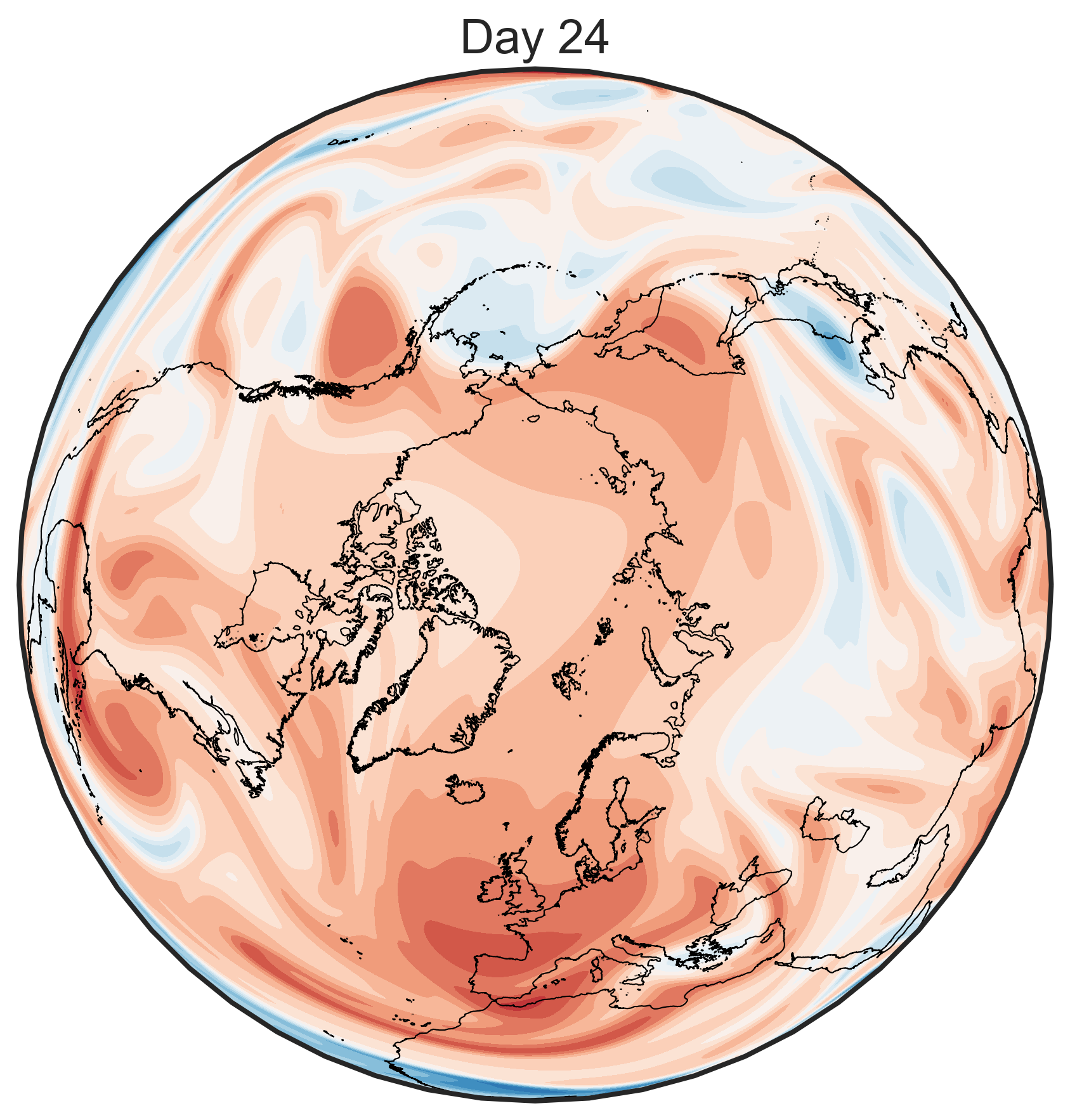}
    \end{subfigure}
    \begin{subfigure}{\textwidth}
        \centering
        \includegraphics[scale=0.25,trim = {0 1.75cm 0 0.31cm}]{cbar.png}
    \end{subfigure}
    \caption{The results of the polar vortex collapse with a wavenumber 1 forcing. The first plot is at day 1, and each plot increments by 1 day. The forcing shuts off at day 15.}
    \label{fig:ss1results}
\end{figure}

\begin{figure}
    \centering
    \begin{subfigure}{0.20\textwidth}
        \includegraphics[width=\linewidth]{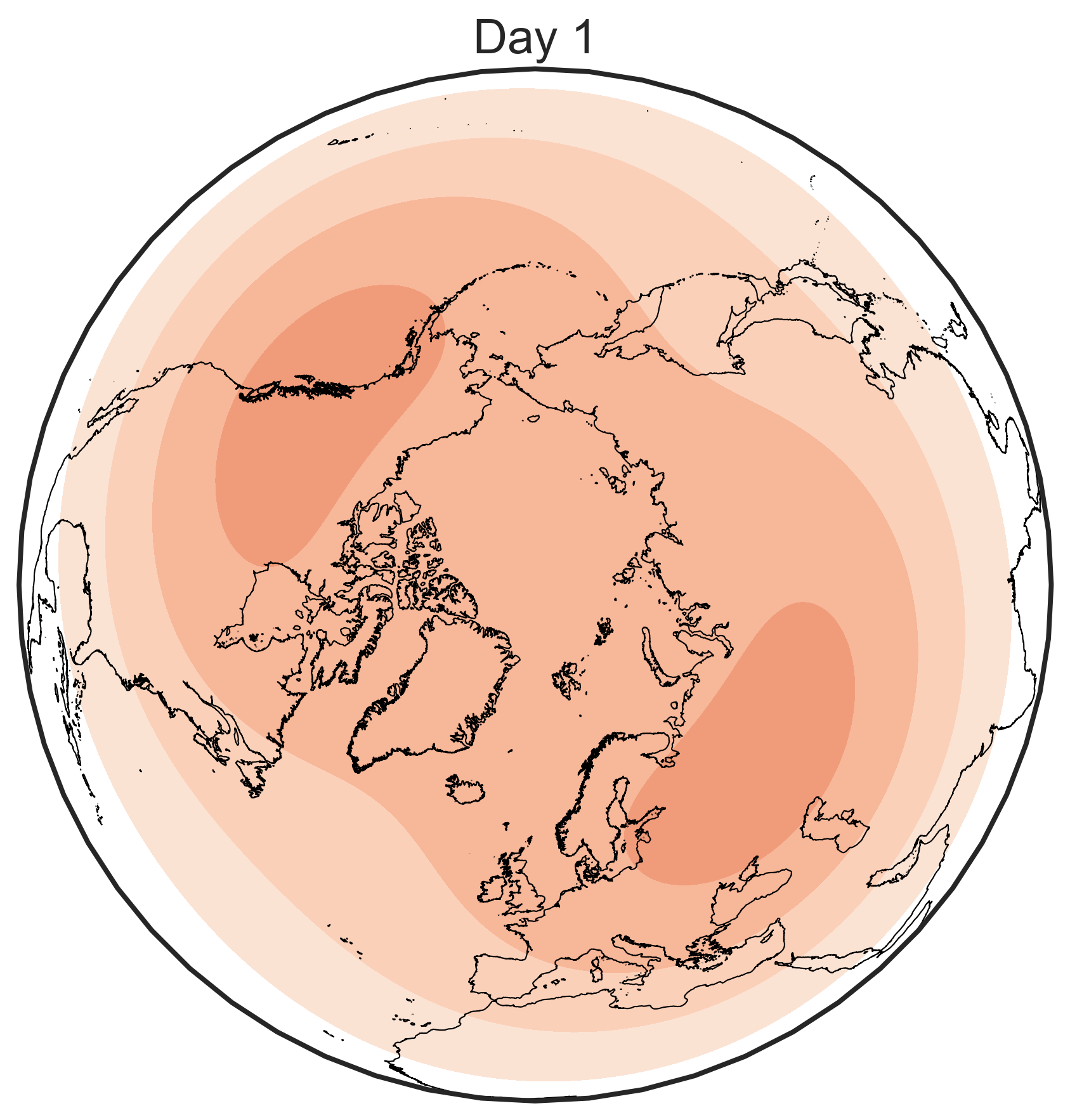}
    \end{subfigure}
    \hspace{0.025\textwidth}
    \begin{subfigure}{0.20\textwidth}
        \includegraphics[width=\linewidth]{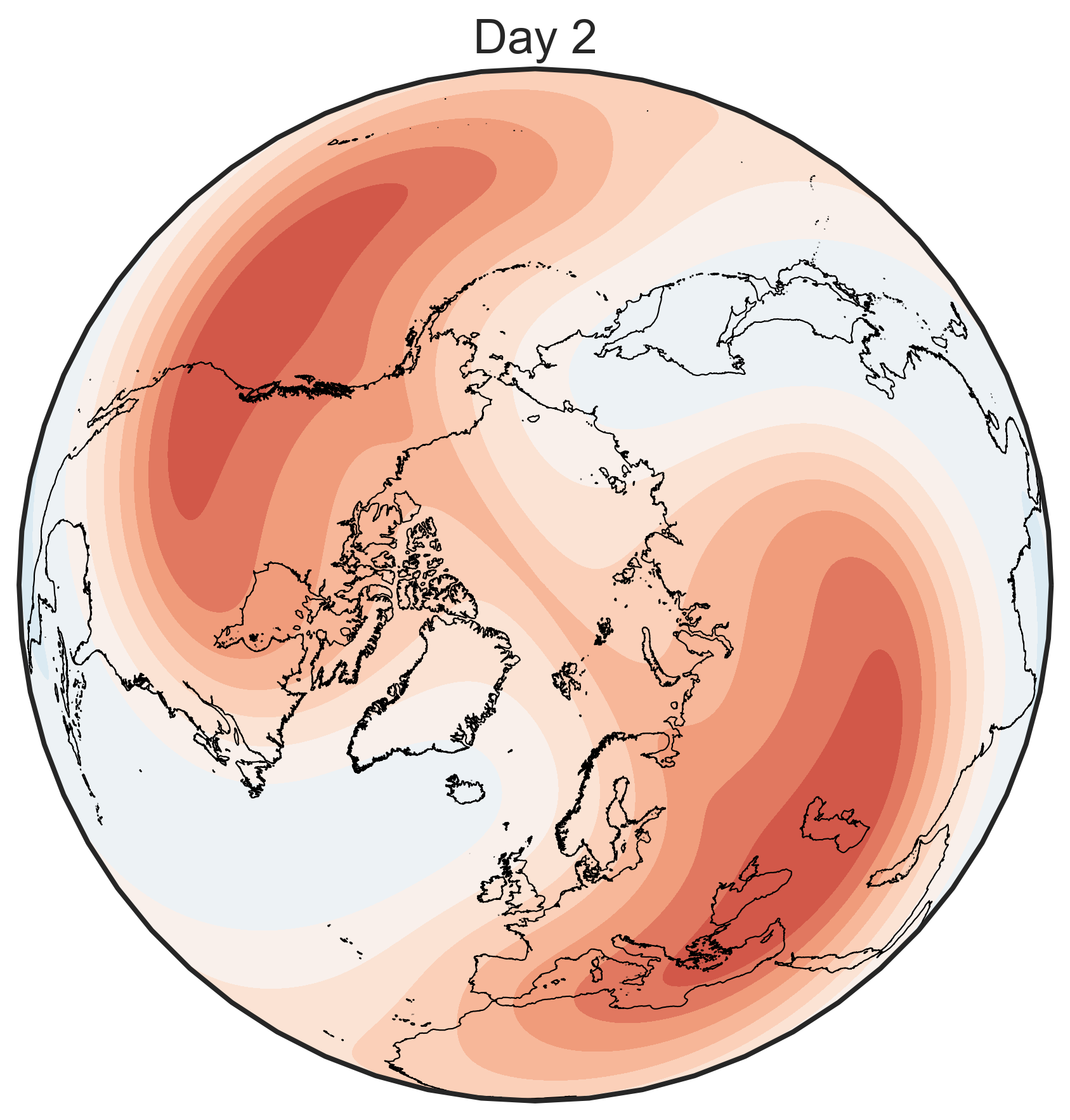}
    \end{subfigure}
    \hspace{0.025\textwidth}
    \begin{subfigure}{0.20\textwidth}
        \includegraphics[width=\linewidth]{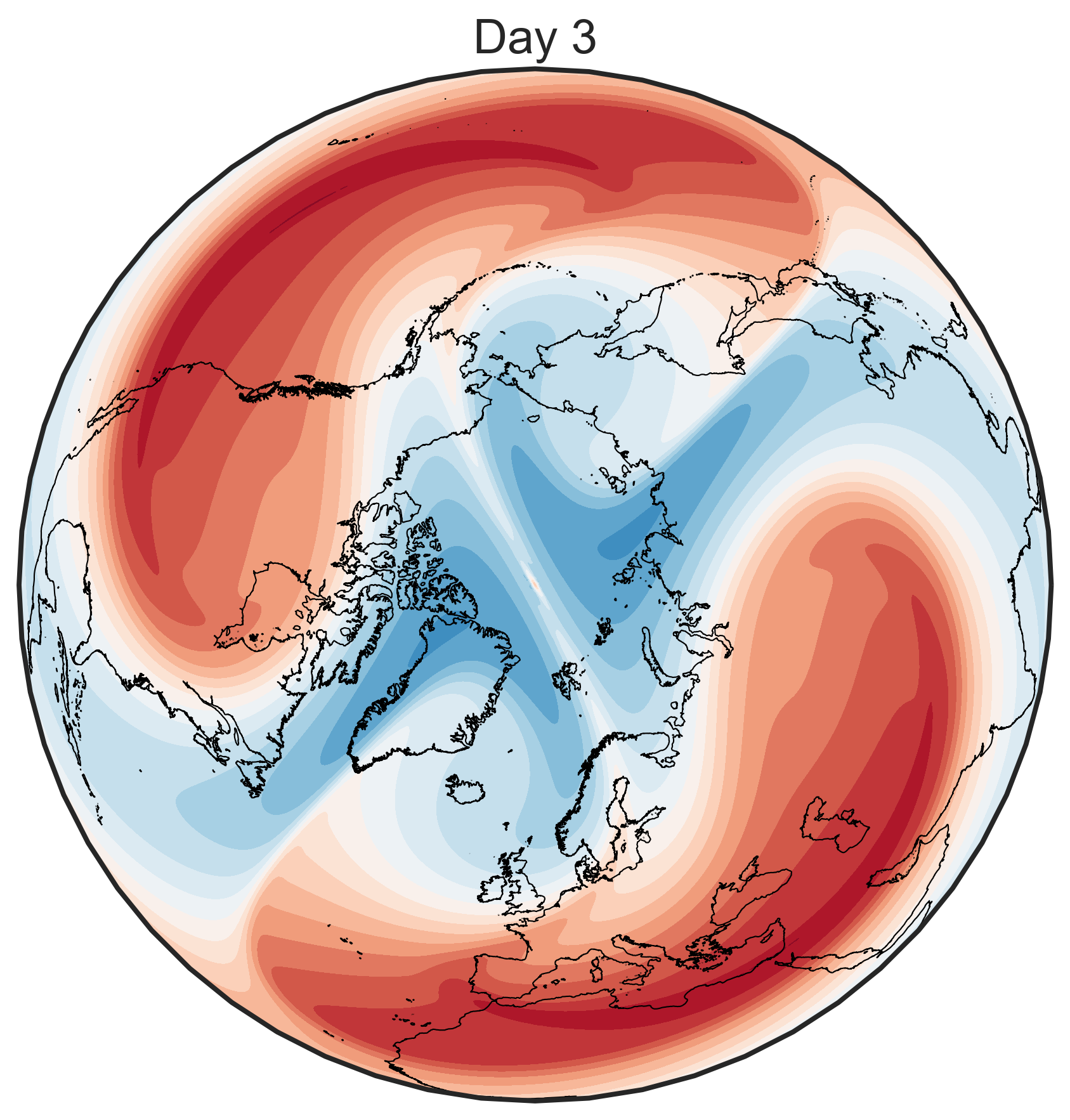}
    \end{subfigure}
    \hspace{0.025\textwidth}
    \begin{subfigure}{0.20\textwidth}
        \includegraphics[width=\linewidth]{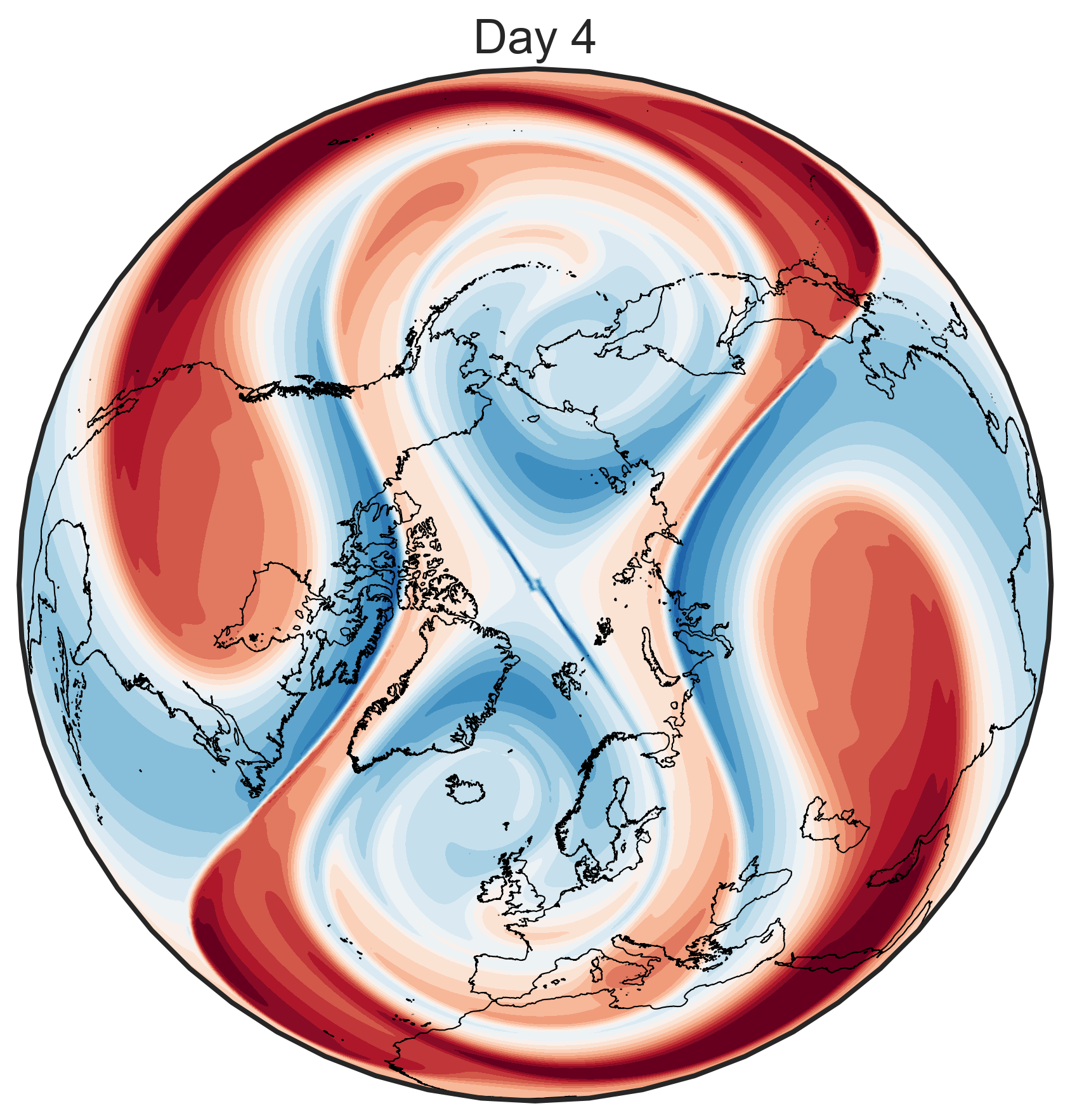}
    \end{subfigure}
    \begin{subfigure}{0.20\textwidth}
        \includegraphics[width=\linewidth]{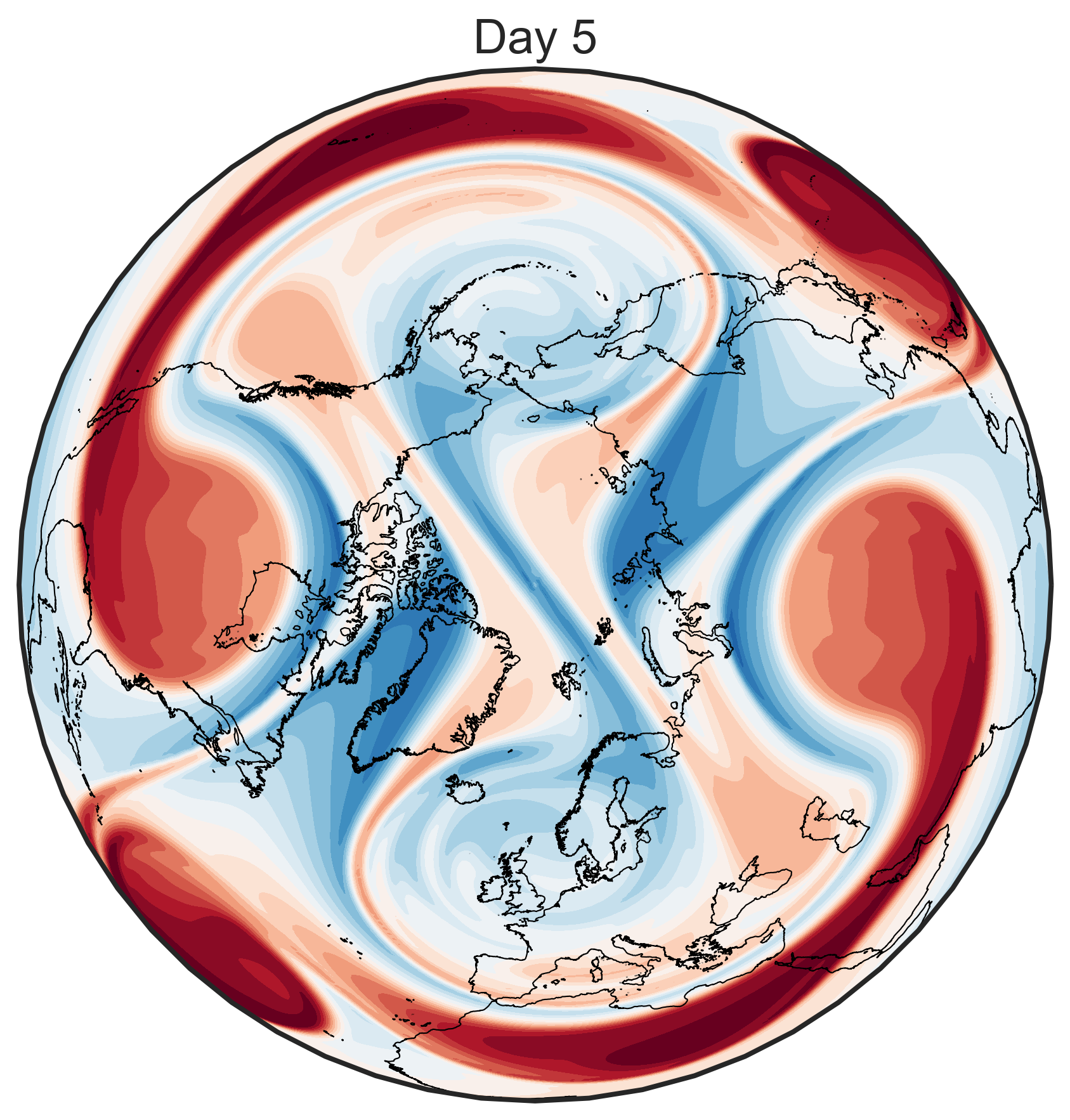}
    \end{subfigure}
    \hspace{0.025\textwidth}
    \begin{subfigure}{0.20\textwidth}
        \includegraphics[width=\linewidth]{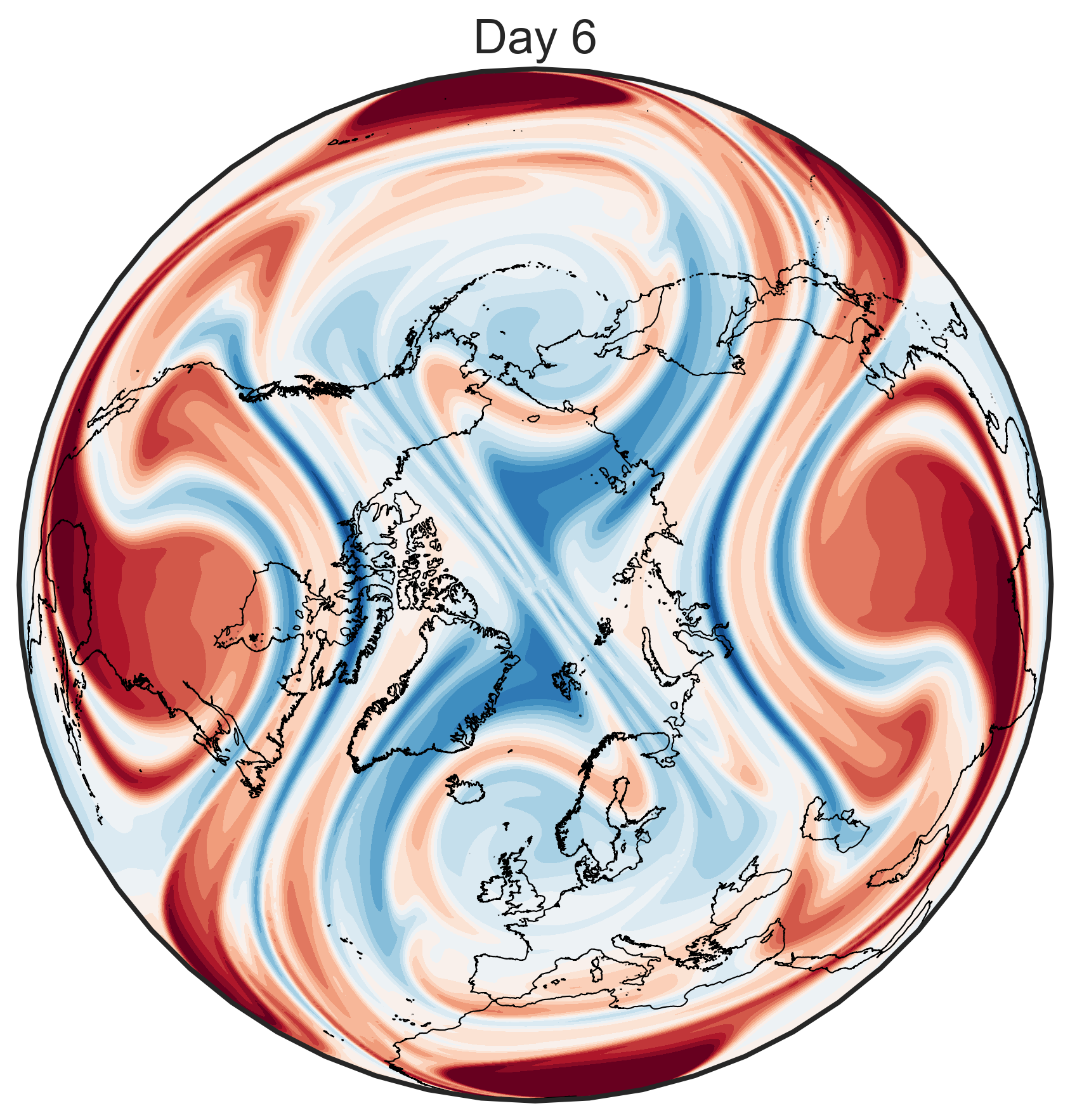}
    \end{subfigure}
    \hspace{0.025\textwidth}
    \begin{subfigure}{0.20\textwidth}
        \includegraphics[width=\linewidth]{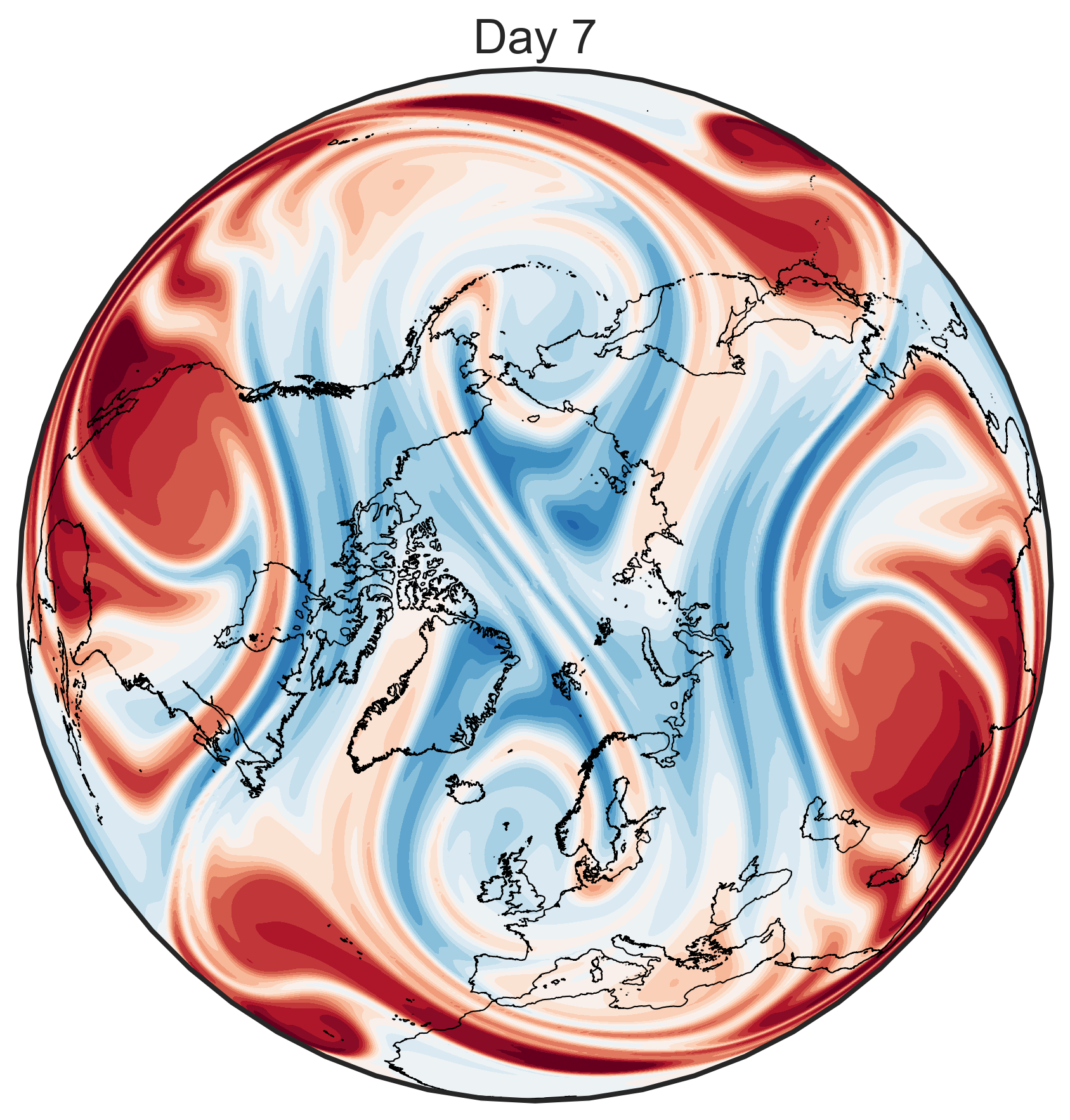}
    \end{subfigure}
    \hspace{0.025\textwidth}
    \begin{subfigure}{0.20\textwidth}
        \includegraphics[width=\linewidth]{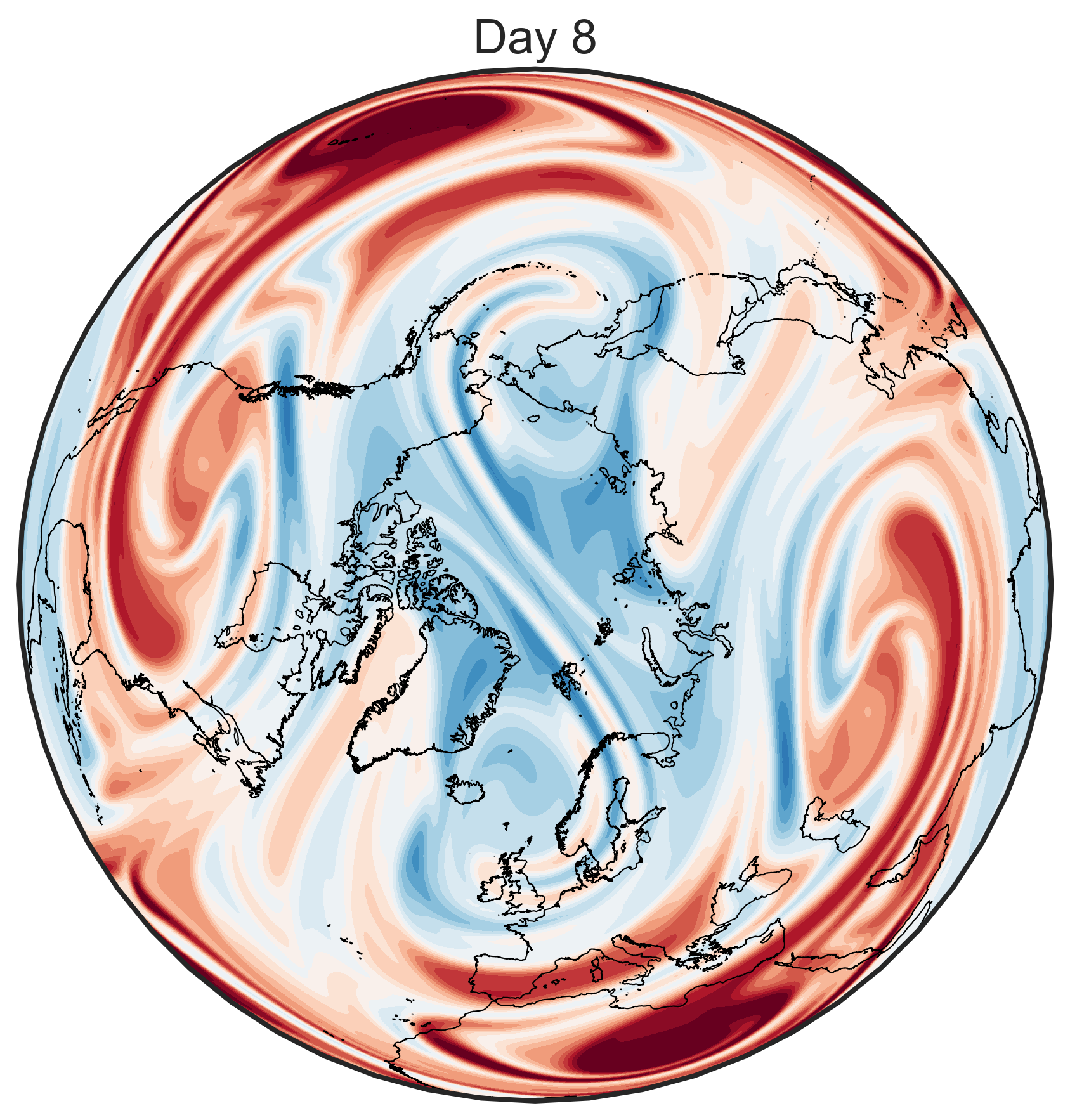}
    \end{subfigure}
    \begin{subfigure}{0.20\textwidth}
        \includegraphics[width=\linewidth]{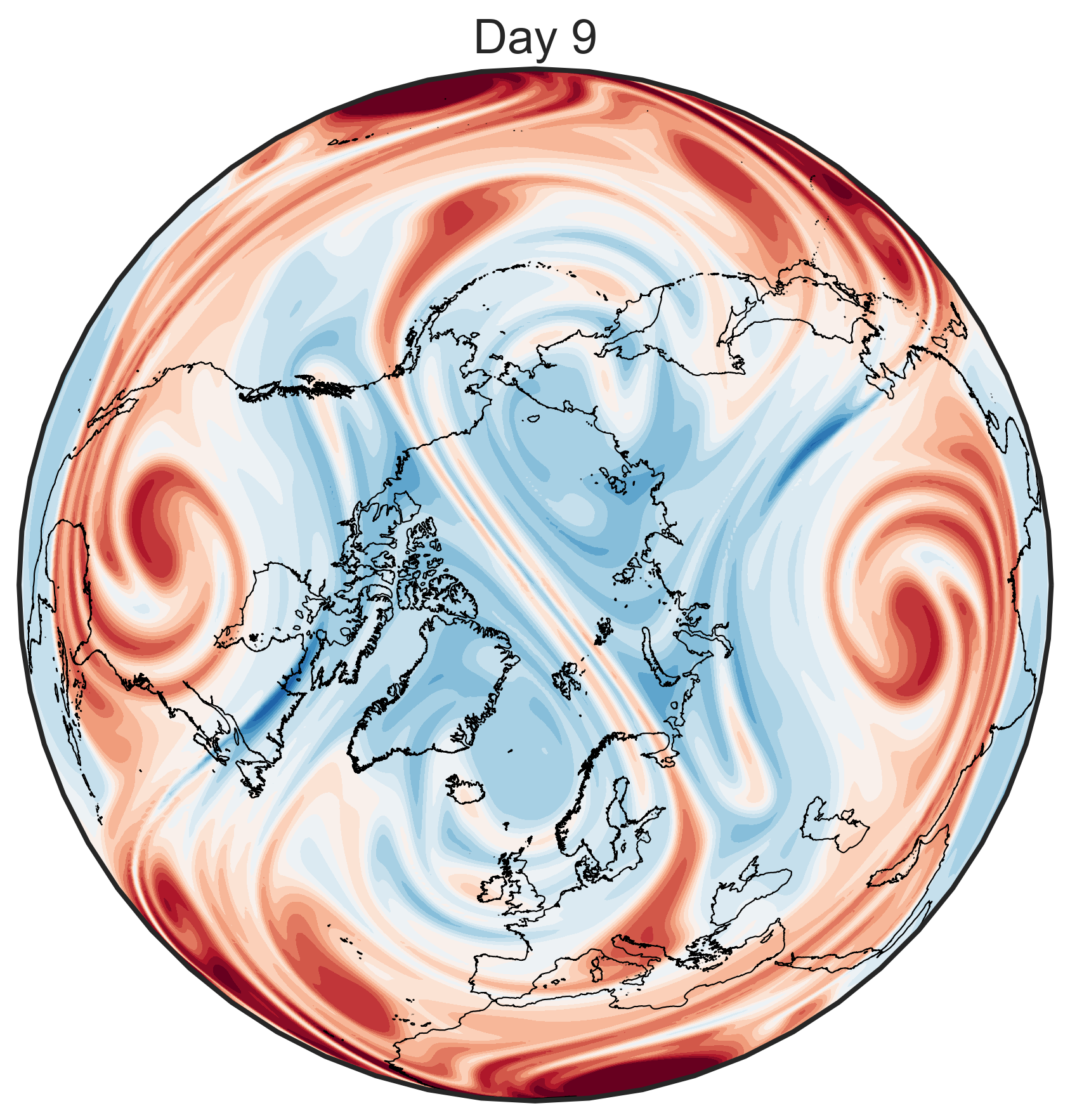}
    \end{subfigure}
    \hspace{0.025\textwidth}
    \begin{subfigure}{0.20\textwidth}
        \includegraphics[width=\linewidth]{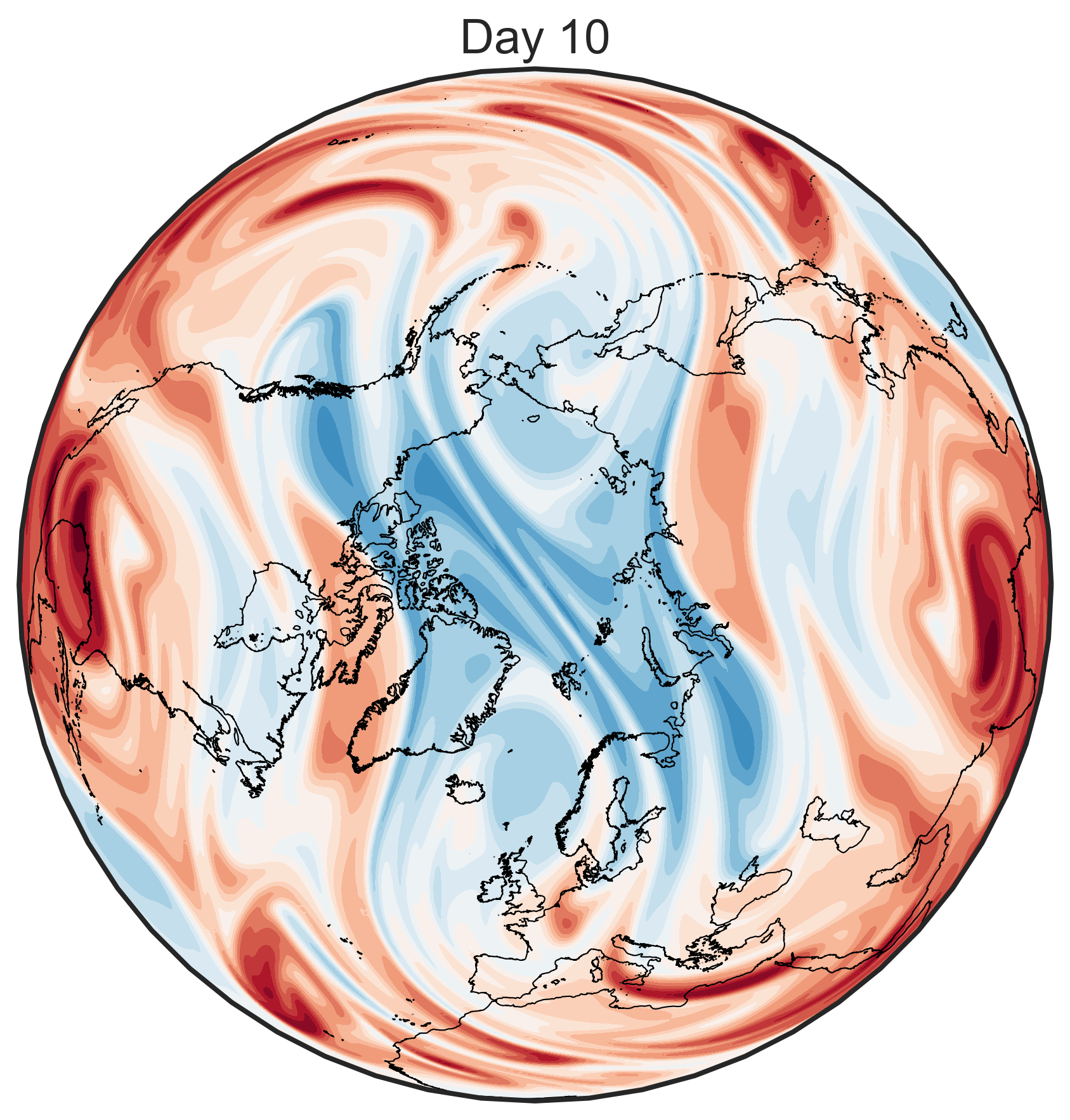}
    \end{subfigure}
    \hspace{0.025\textwidth}
    \begin{subfigure}{0.20\textwidth}
        \includegraphics[width=\linewidth]{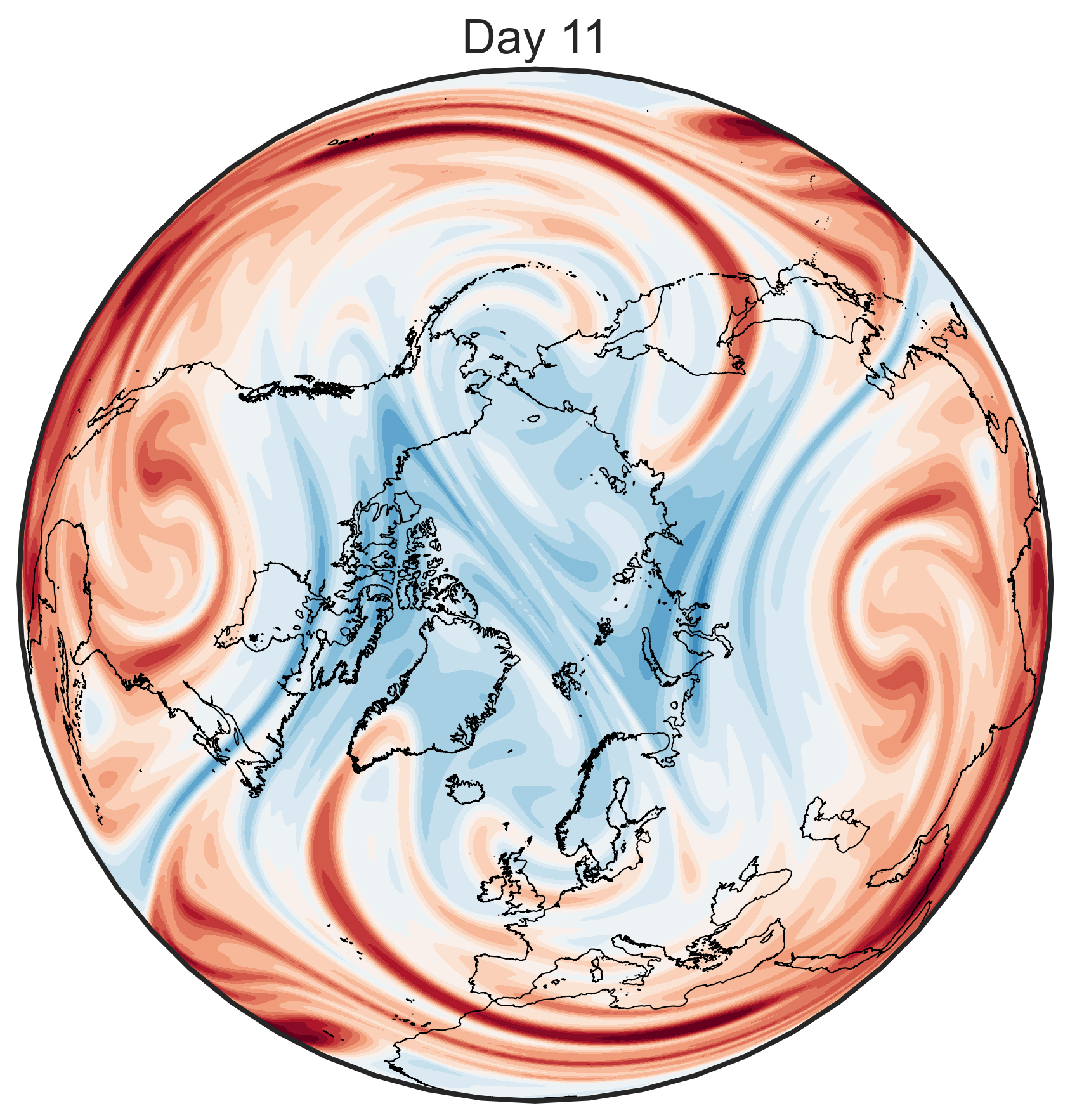}
    \end{subfigure}
    \hspace{0.025\textwidth}
    \begin{subfigure}{0.20\textwidth}
        \includegraphics[width=\linewidth]{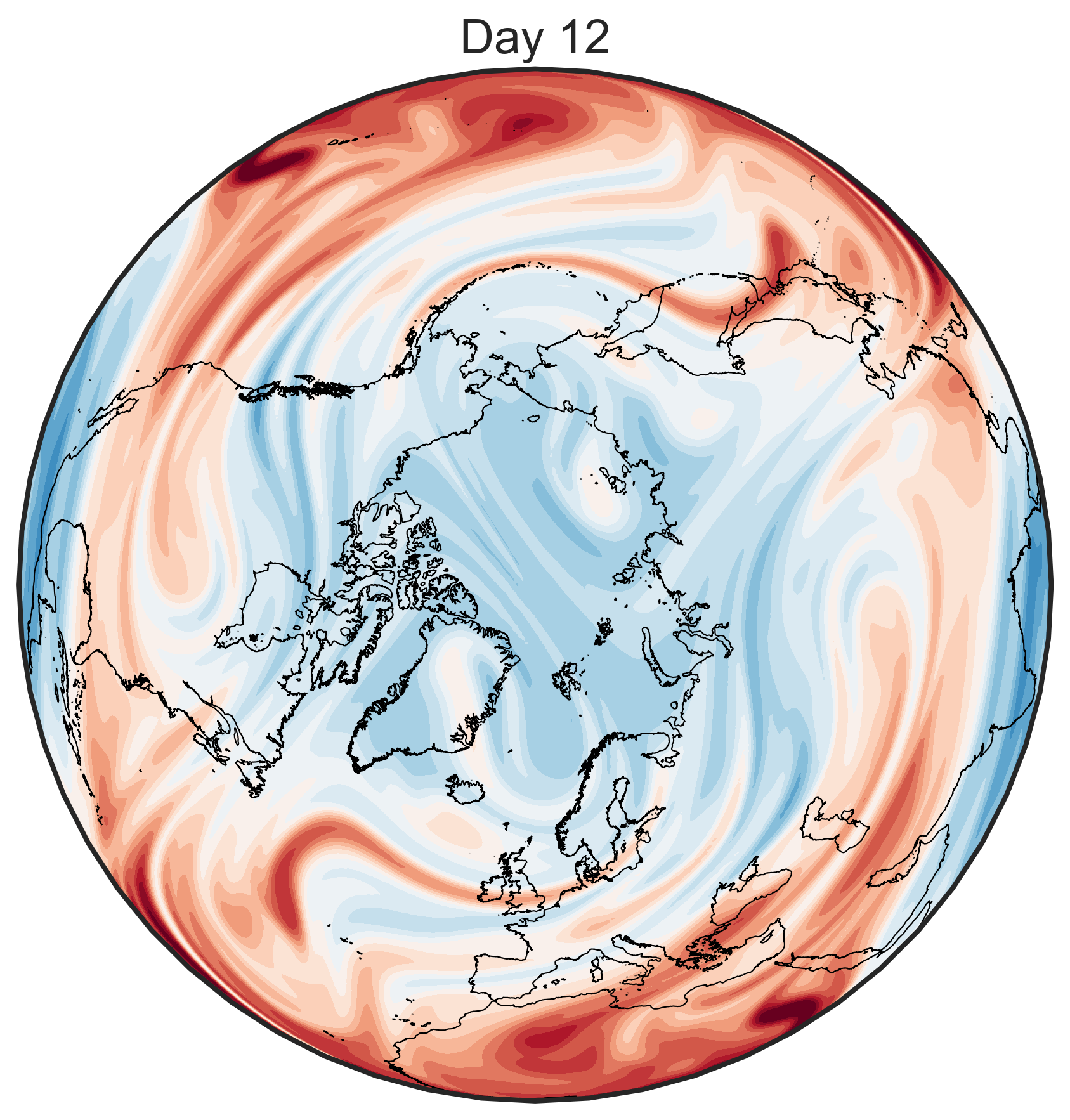}
    \end{subfigure}
    \begin{subfigure}{0.20\textwidth}
        \includegraphics[width=\linewidth]{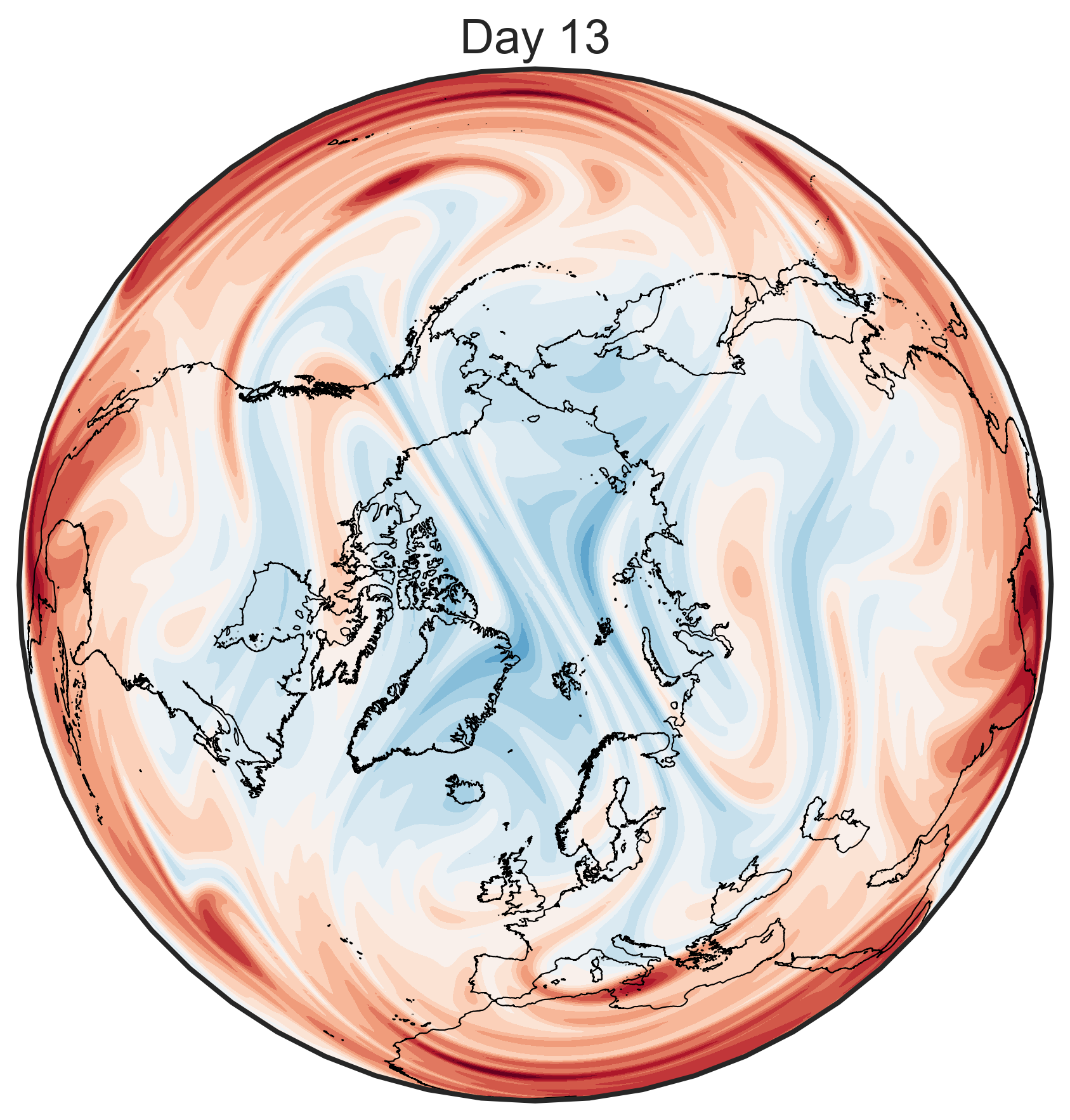}
    \end{subfigure}
    \hspace{0.025\textwidth}
    \begin{subfigure}{0.20\textwidth}
        \includegraphics[width=\linewidth]{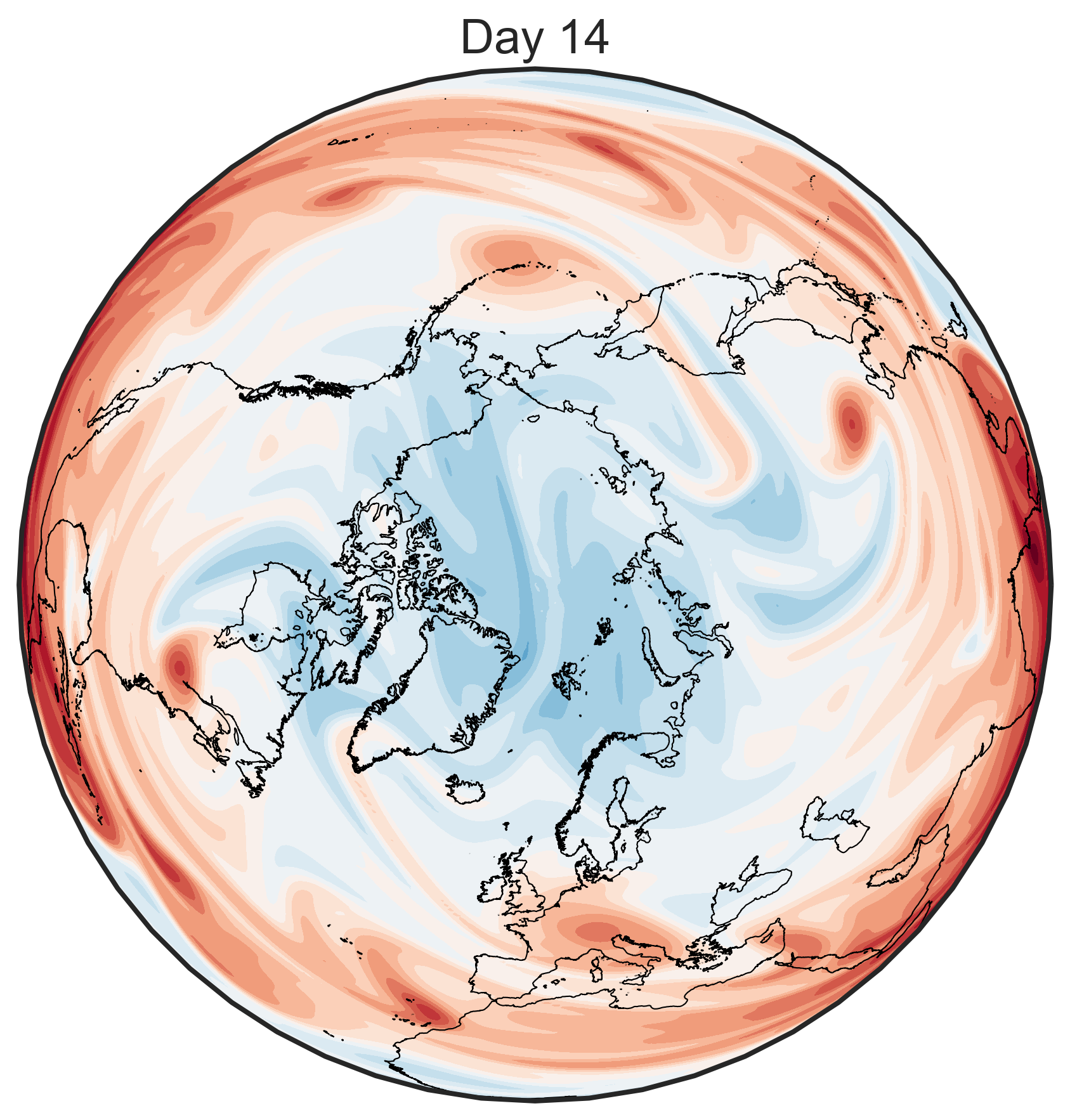}
    \end{subfigure}
    \hspace{0.025\textwidth}
    \begin{subfigure}{0.20\textwidth}
        \includegraphics[width=\linewidth]{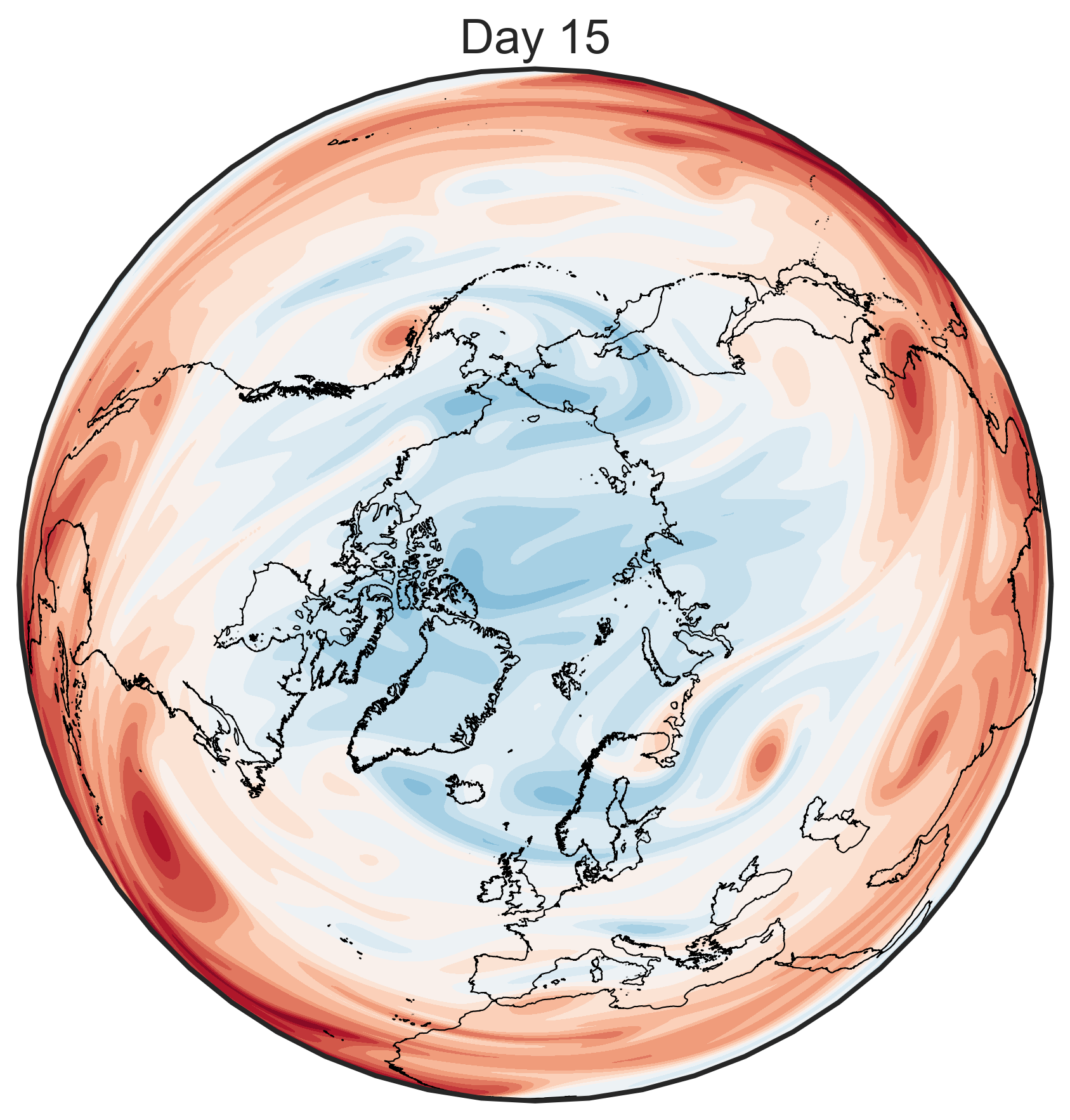}
    \end{subfigure}
    \hspace{0.025\textwidth}
    \begin{subfigure}{0.20\textwidth}
        \includegraphics[width=\linewidth]{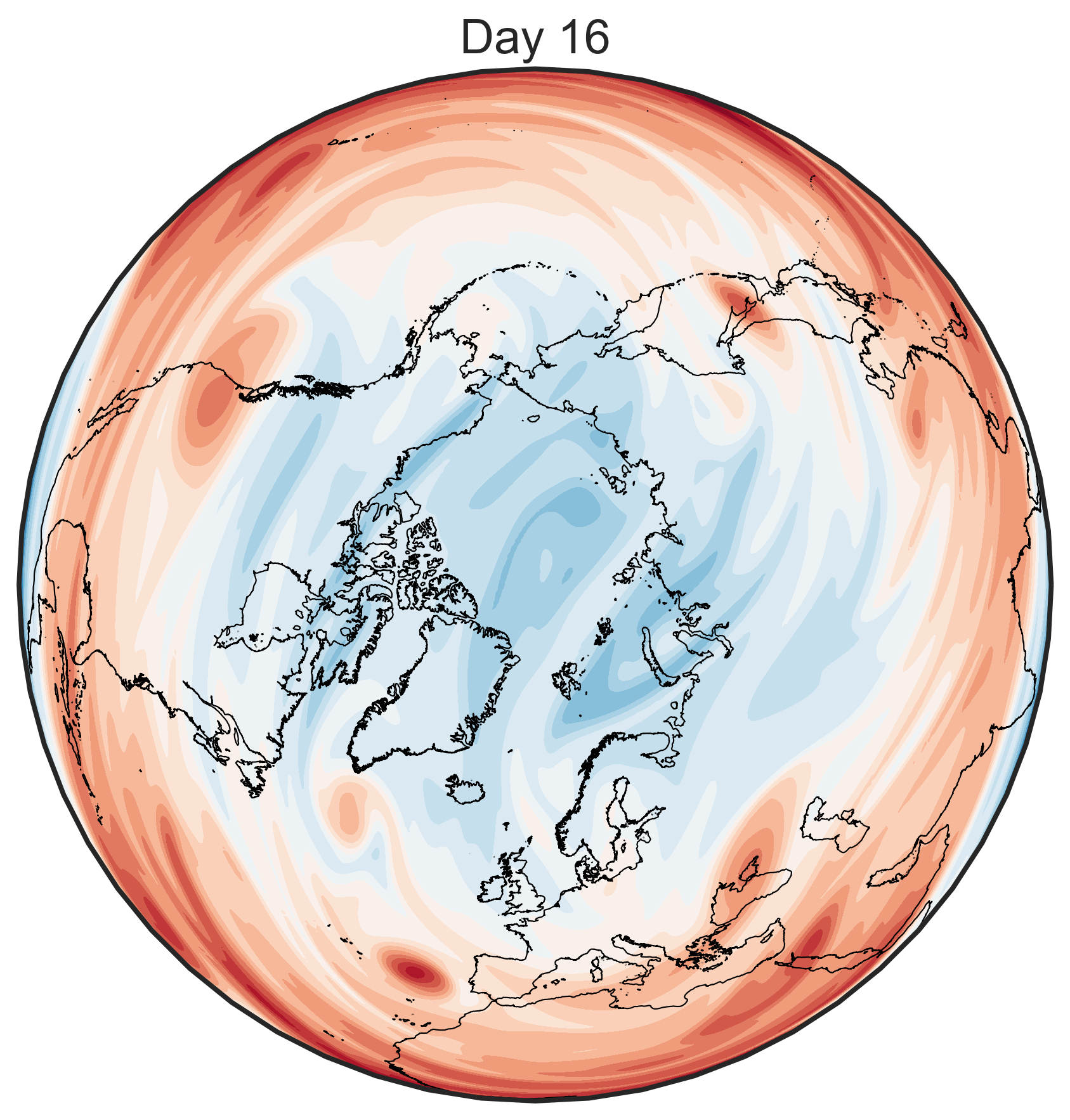}
    \end{subfigure}
    \begin{subfigure}{0.20\textwidth}
        \includegraphics[width=\linewidth]{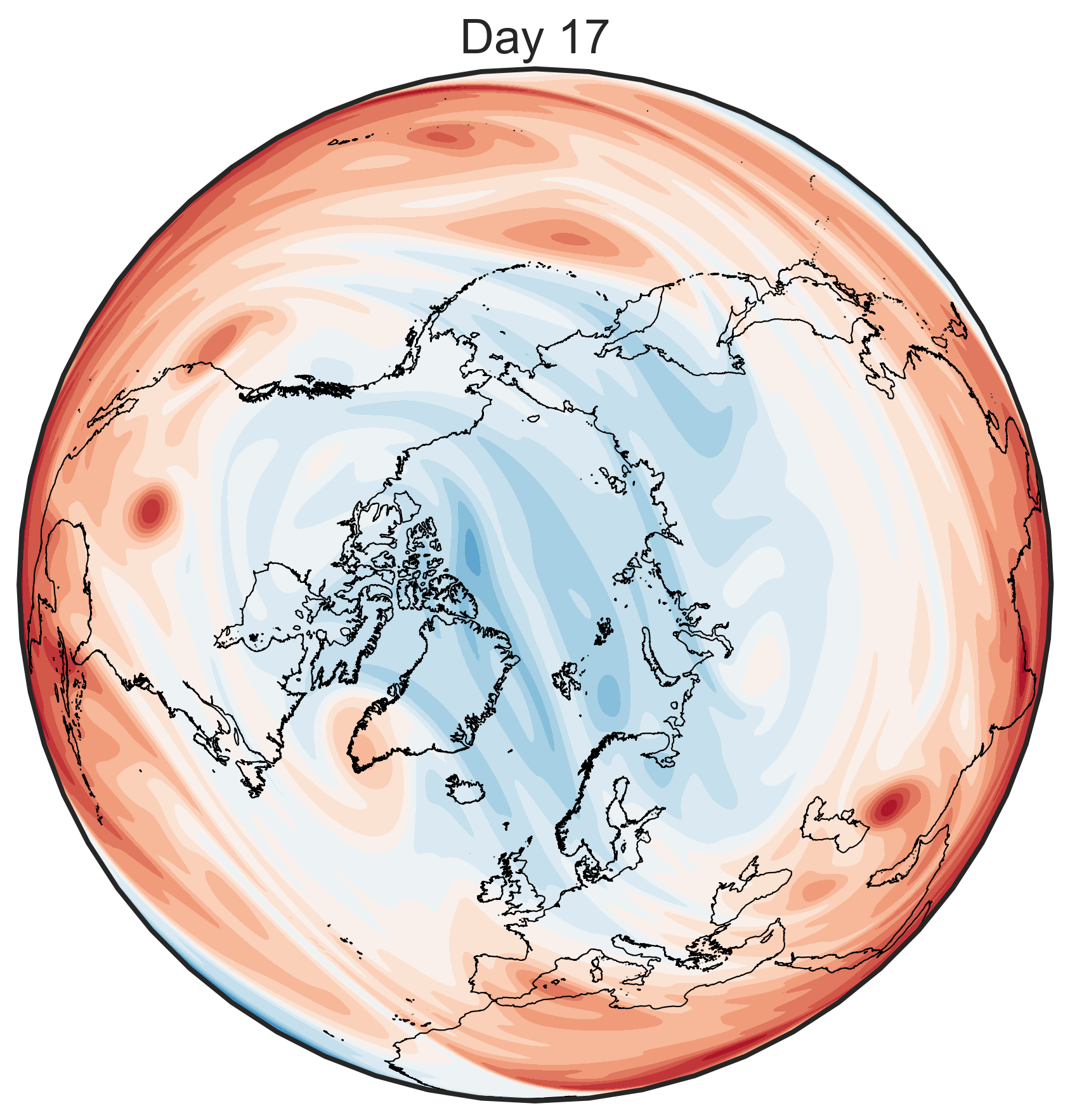}
    \end{subfigure}
    \hspace{0.025\textwidth}
    \begin{subfigure}{0.20\textwidth}
        \includegraphics[width=\linewidth]{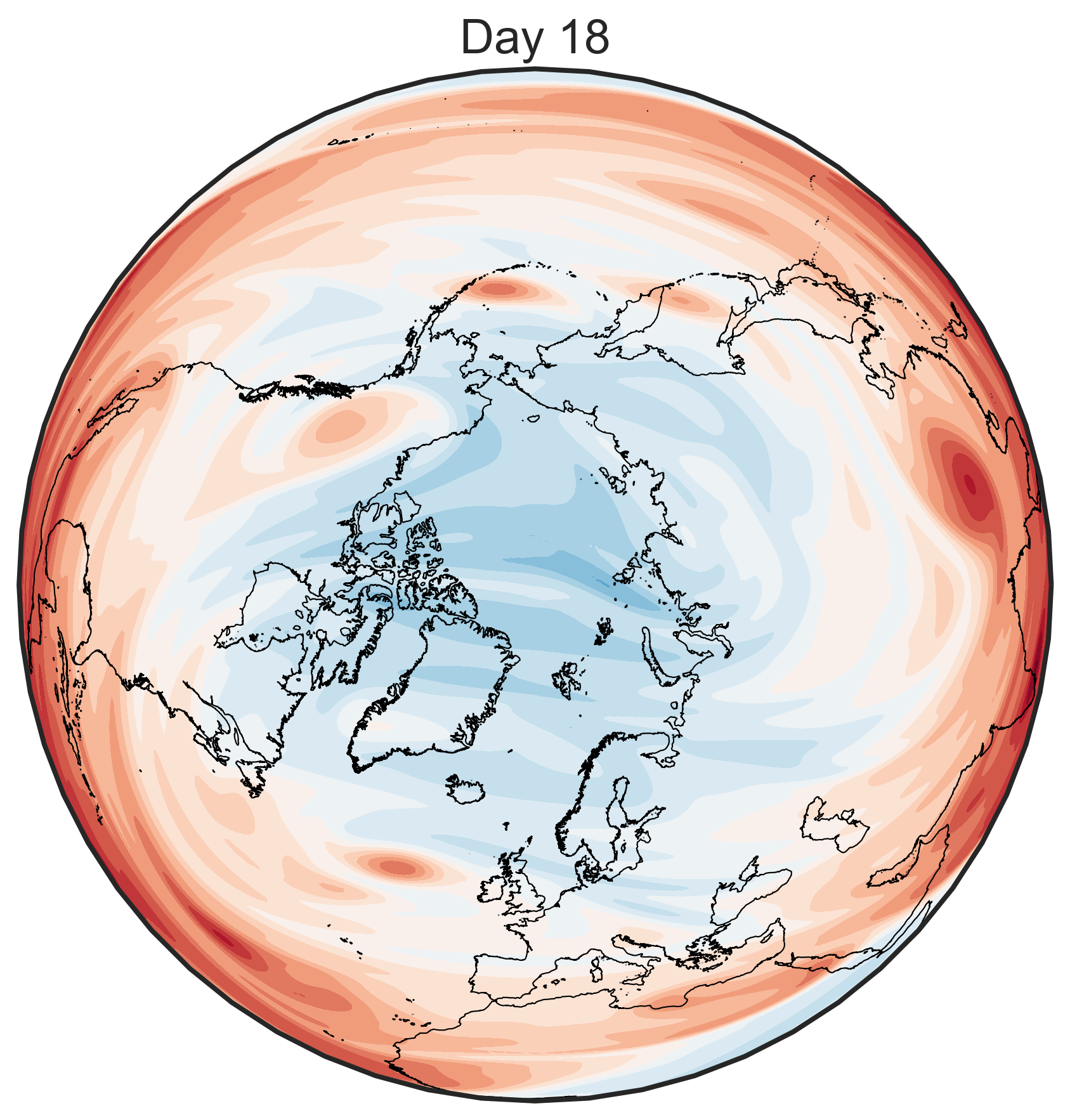}
    \end{subfigure}
    \hspace{0.025\textwidth}
    \begin{subfigure}{0.20\textwidth}
        \includegraphics[width=\linewidth]{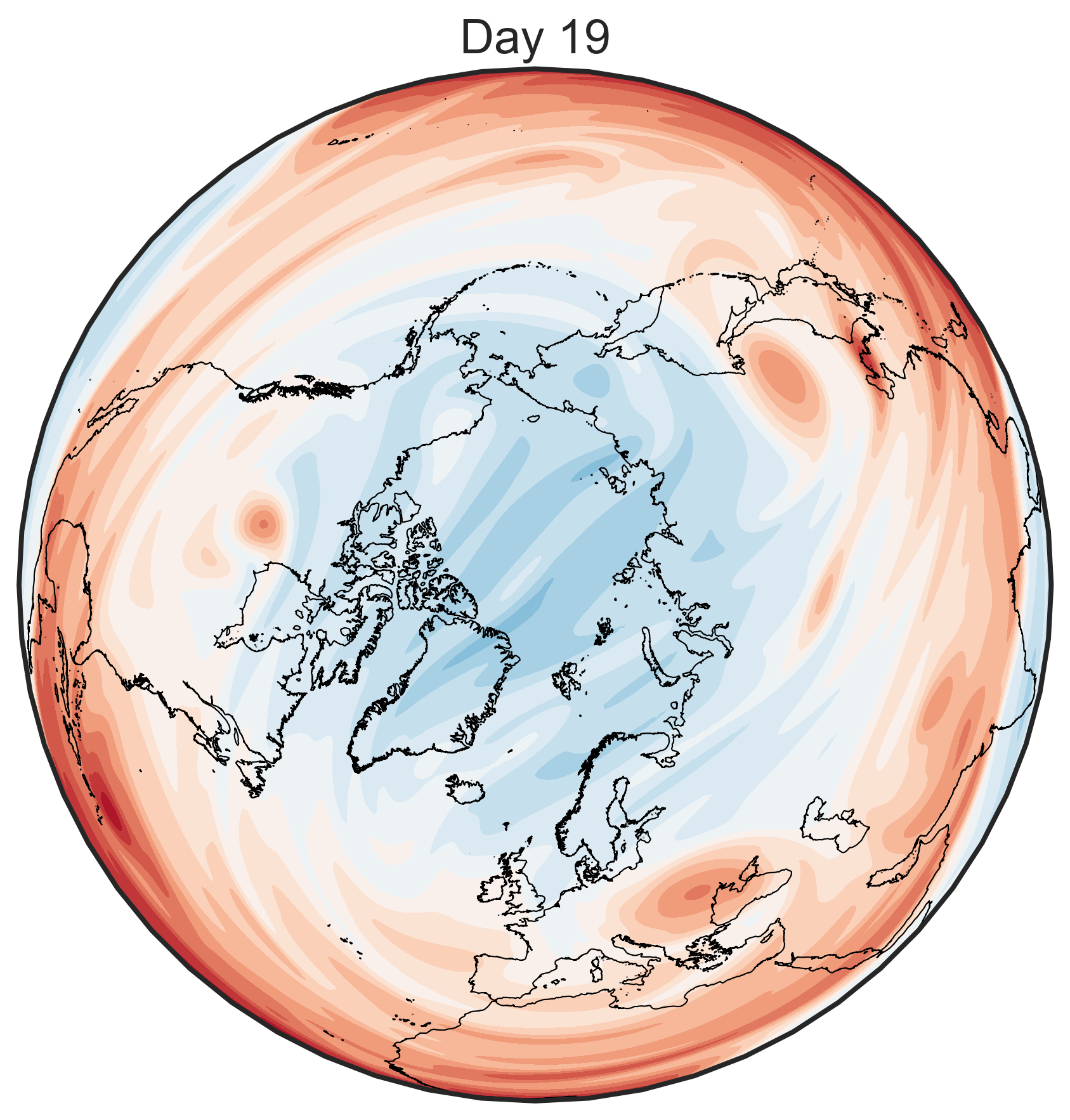}
    \end{subfigure}
    \hspace{0.025\textwidth}
    \begin{subfigure}{0.20\textwidth}
        \includegraphics[width=\linewidth]{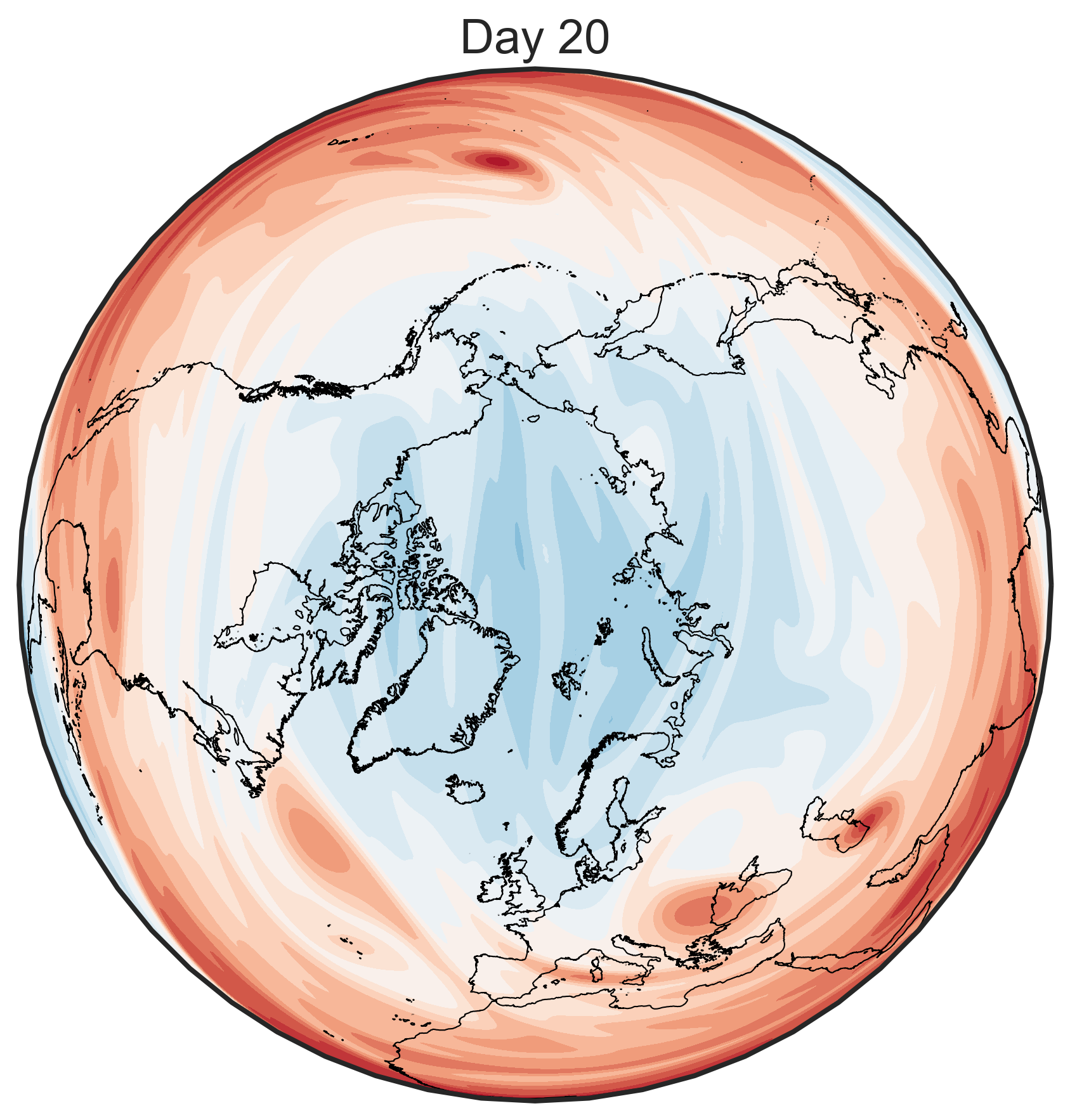}
    \end{subfigure}
    \begin{subfigure}{0.20\textwidth}
        \includegraphics[width=\linewidth]{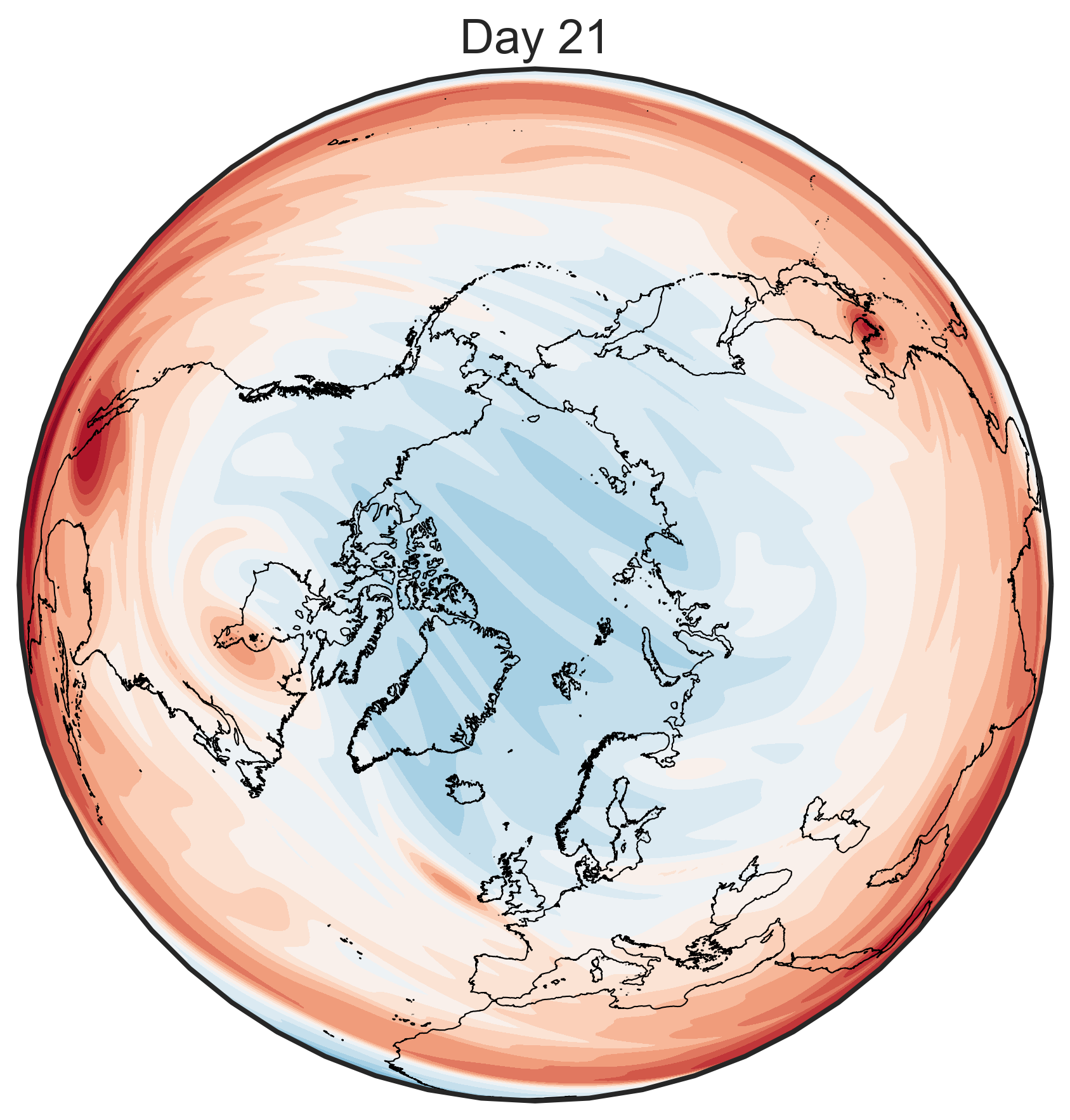}
    \end{subfigure}
    \hspace{0.025\textwidth}
    \begin{subfigure}{0.20\textwidth}
        \includegraphics[width=\linewidth]{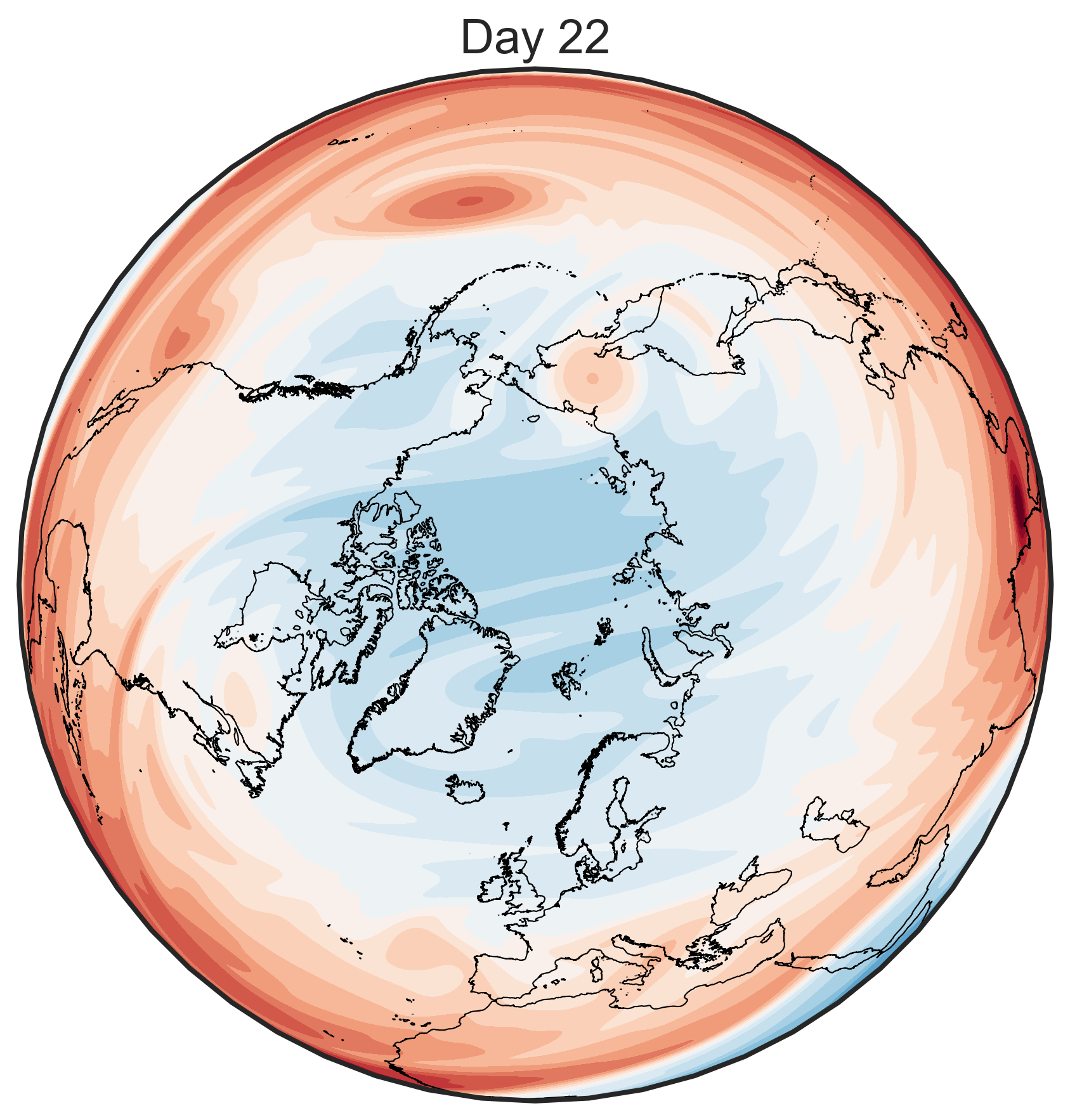}
    \end{subfigure}
    \hspace{0.025\textwidth}
    \begin{subfigure}{0.20\textwidth}
        \includegraphics[width=\linewidth]{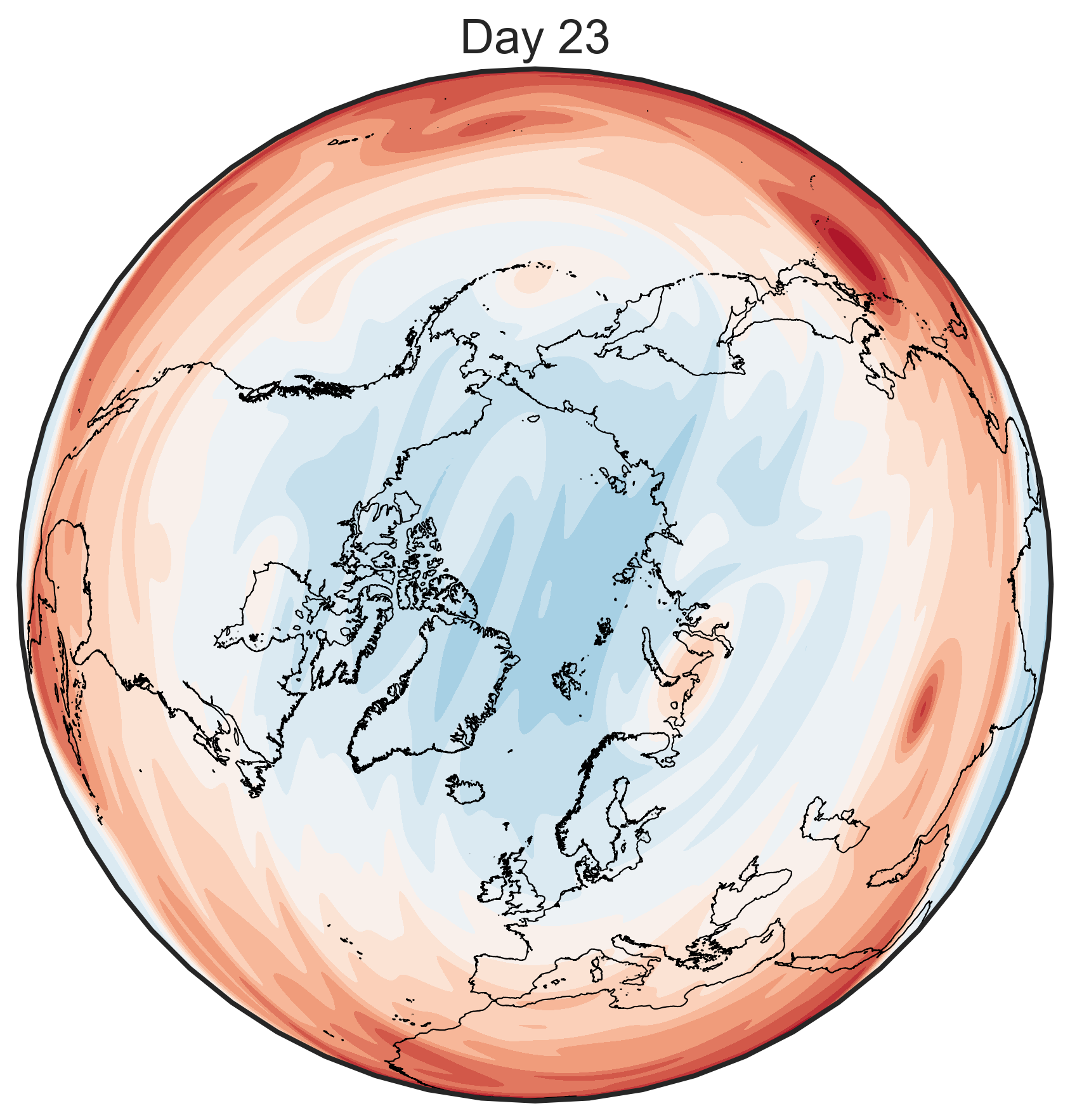}
    \end{subfigure}
    \hspace{0.025\textwidth}
    \begin{subfigure}{0.20\textwidth}
        \includegraphics[width=\linewidth]{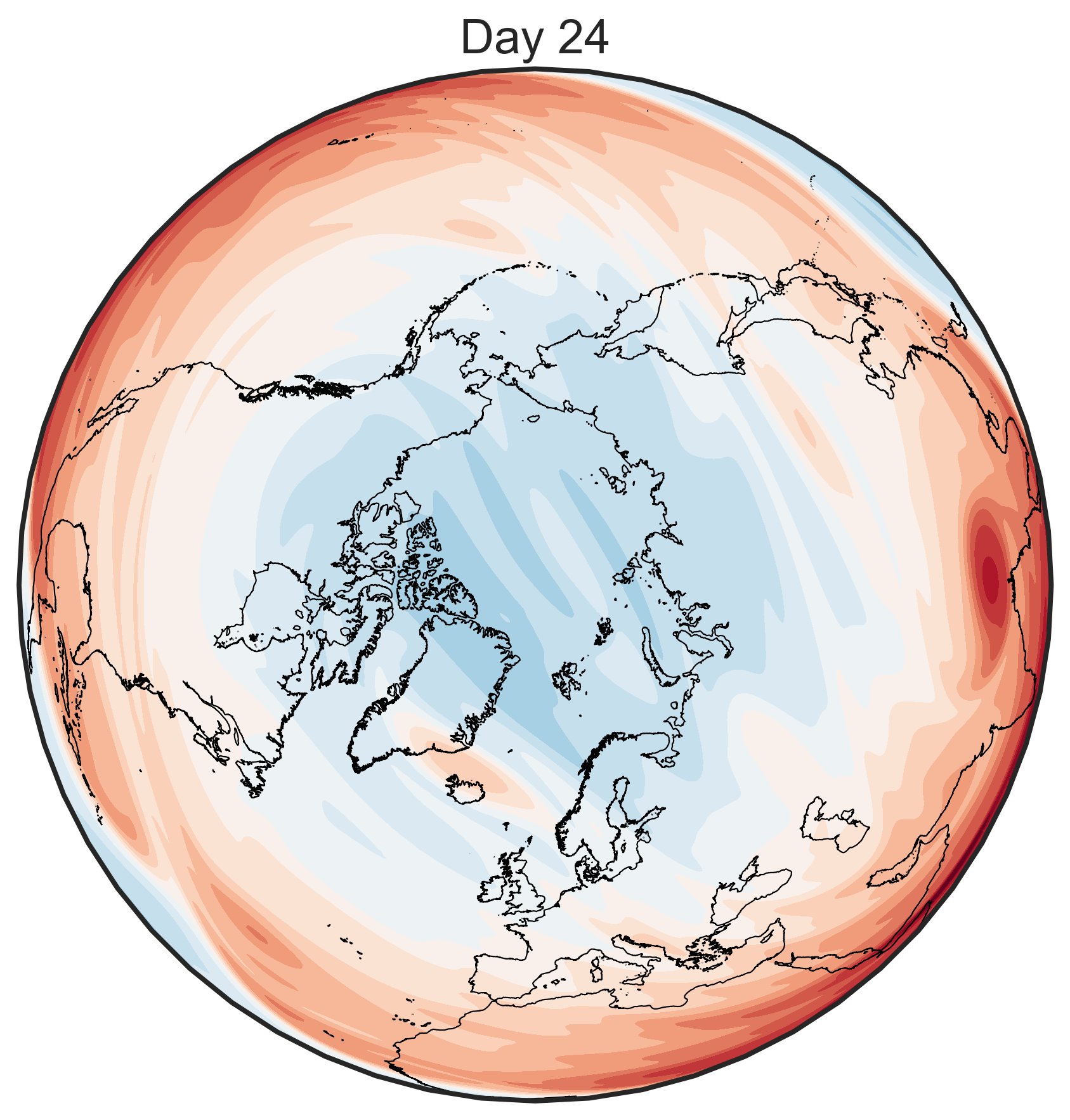}
    \end{subfigure}
    \begin{subfigure}{\textwidth}
        \centering
        \includegraphics[scale=0.25,trim = {0 1.75cm 0 0.31cm}]{cbar.png}
    \end{subfigure}
    \caption{The results of the polar vortex collapse with a wavenumber 2 forcing. The first plot is at day 1, and each plot increments by 1 day. The forcing shuts off at day 15.}
    \label{fig:ss2results}
\end{figure}

We plot the results of these simulations in Figures~\ref{fig:ss1results} and~\ref{fig:ss2results}. We perform these simulations with 163842 particles, $d=6$, and $\theta=0.7$. These two forcings behave as expected, and demonstrate that we can capture different polar vortex collapse morphologies in our simplified model. With the wavenumber 1 forcing, we see that the polar vortex is first displaced off the pole before we see the development of a tripolar vortex structure. Once the forcing shuts off at day 15, we see that the polar vortex slowly starts to recover. With the wavenumber 2 forcing, we see that the polar vortex is very quickly split, and a region of negative relative vorticity forms at the pole. This region of negative vorticity persists through the duration of the model run, and the polar vortex does not recover its initial state, with a polar vortex that instead has anticyclonic behavior. 

\section{Conclusion and Outlook}
\label{sec:conc}

In this work, we presented a novel spherical tree code for performing spherical convolutions in $O(N\log{N})$ time, and applied it to solving the dynamics of the barotropic vorticity equation on a sphere. The spherical tree code uses triangular polynomial interpolation on an iteratively refined icosahedron to perform cluster-particle, particle-cluster, and cluster-cluster interactions. This spherical tree code is kernel independent and is suitable for a wide range of kernels on the sphere, including those that are only known numerically. We demonstrated the improved scaling over direct summation, and also demonstrated the effectiveness of the MPI parallelization with a strong scaling test. We also compared the error behavior with that of the BLTC. 

We then tested the spherical tree code on three different barotropic vorticity equation test cases. The first was the Rossby Haurwitz wave, which is an exact solution of the barotropic vorticity equation, and allowed us to exactly measure the error. Next, we tested it on the Gaussian vortex test case, which allowed us to see that the fast summation worked well with the adaptive mesh refinement scheme. This was also a test case without an exact solution, and our results with fast summation compared well with the results computed using direct summation. Lastly, we tested the polar vortex collapse test case, in which we add a vorticity forcing to collapse the polar vortex. This is a challenging test case with steep vorticity gradients and a long runtime, and we saw that the forcings do behave as predicted, while also leading to a wide range of dynamics, including the formation of tripolar vortex structures. 

We plan on applying this fast summation technique to more challenging problems, such as the shallow water equations (SWE) and the computation of Self Attraction and Loading (SAL), as well as running this tree code on Graphics Processing Units (GPUs) in order to go to even larger problem sizes. The SWE state 
\begin{align}\begin{split}
\frac{D\mathbf{u}}{Dt}&=-\mathbf{f}\times\mathbf{u}-g\nabla h\\
\frac{Dh}{Dt}&=-h\nabla\cdot\mathbf{u}
\end{split}\end{align}
where $\mathbf{u}$ is the horizontal fluid velocity, $\mathbf{f}=f\mathbf{e}_z$ is the Coriolis parameter, $g$ is gravity, and $h$ is the height of the fluid layer. A Lagrangian discretization of the SWE would have multiple kernels, accounting for both the vorticity and the divergence, and we plan on applying our spherical tree code to this problem. The SAL is a term that occurs in ocean modeling \cite{barton2022global, brus2023scalable}. When accounting for tidal motions in the SWE, one has the equation set 
\begin{align}\begin{split}
\frac{D\mathbf{u}}{Dt}&=-\mathbf{f}\times\mathbf{u}-g\nabla(h-\zeta_{\mathrm{EQ}}-\zeta_{\mathrm{SAL}})\\
\frac{Dh}{Dt}&=-h\nabla\cdot\mathbf{u}
\end{split}\end{align}
where we have an additional $\zeta_{\mathrm{EQ}}$ term representing the equilibrium tide due to astronomical forcing, and $\zeta_{\mathrm{SAL}}$, accounting for the self attraction of the mass of the water column, changes in the shape of the Earth due to elastic loading, and subsequent changes to the Earth's gravitational potential due to the changes in the Earth's shape. Not including the SAL in tidal computations can result in amplitude errors of 20\% and significant phase errors. The traditional method of computing SAL has been to compute the spherical harmonic components $\widehat{h}_n^m$ of the sea surface height, and multiply each component by a coefficient based on the Love numbers \cite{farrell1972deformation}. Then, the spherical harmonic components can be summed together to give the $\zeta_{\mathrm{SAL}}$ term. However, by the convolution theorem, we can rewrite this as a convolution with a Green's function, and when discretized, this can be rewritten in a way to take advantage of our spherical tree code. 

\section{Acknowledgements}
\label{sec:acknow}

We would like to thank Peter Bosler (Sandia National Laboratories) and Robert Krasny (University of Michigan) for the helpful discussions. 

This material is based upon work supported by the National Science Foundation
Graduate Research Fellowship under Grant No. DGE-1841052. Any opinion, findings, and conclusions or recommendations expressed in this material are those of the authors(s) and do not necessarily reflect the views of the National Science Foundation. This research was supported in part through computational resources and services provided by Advanced Research Computing at the University of Michigan, Ann Arbor. 

 \appendix
 \section{Adaptive Mesh Refinement}
 \label{sec:amr:appendix}
 
 We briefly describe our AMR scheme. First, because we are using particles from an iteratively refined icosahedron, our particles form a triangulation of the sphere with a tree structure. Then, to perform adaptive mesh refinement, we loop over the leaf triangles, checking if the triangle and the particles at its vertices have a large circulation 
 \begin{equation}A\frac{\zeta_1+\zeta_2+\zeta_3}{3}\geq\varepsilon_1
 \end{equation} where $A$ is the area of the triangle, the $\zeta_i$ are the vorticities of the vertex particles, and $\varepsilon_1$ is some threshold, or if the vorticity variation is large
 \begin{equation}\max\{\zeta_1,\zeta_2,\zeta_3\}-\min\{\zeta_1,\zeta_2,\zeta_3\}\geq\varepsilon_2
 \end{equation} where $\varepsilon_2$ is another parameter. If either inequality is true, we refine the triangle, splitting it into four smaller triangles and introducing new particles. We limit the maximum refinement to three levels. When these criteria are no longer met, we coarsen the grid, but not past the initial resolution. AMR is performed at each time step. The values of epsilon that are used in the Gaussian Vortex test case are stated in Sec.~\ref{sec:results:gv}. 

 \bibliographystyle{elsarticle-num} 
 \bibliography{fast-bve}

\end{document}